%% file: paper.tex
\documentclass[pamq,keywordasfootnote]{ipart}
\usepackage{color,amsmath,amssymb,amsthm,hyperref,graphicx}
\usepackage[normalem]{ulem}

\startlocaldefs
\numberwithin{equation}{section}
\newcommand{\diff}{\mathrm{d}}

\newtheorem*{conj}{Conjecture and Definition}
\newtheorem{conjecture}{Conjecture}
\newtheorem*{define}{Definition}
\newtheorem*{rems}{Remarks}
\newtheorem{convention}[conjecture]{Convention}
\newcommand{\definedas}{\mathrel{\raise.095ex\hbox{\rm :}\mkern-5.2mu=}}
 \newcommand{\asdefined}{\mathrel{=\mkern-5.2mu\raise.095ex\hbox{\rm :}}}

\endlocaldefs

\begin{document}

\begin{frontmatter}

\title[A flow approach to Bartnik's static metric extension conjecture]
  {A flow approach to Bartnik's static metric extension conjecture in axisymmetry}

\begin{aug}
  \author{\fnms{Carla} \snm{Cederbaum}%
    \ead[label=e1]{cederbaum@math.uni-tuebingen.de}},
   \address{Department of Mathematics, University of T\"ubingen,\\ 
     Auf der Morgenstelle 10, 72076 T\"ubingen, Germany\\ 
     \printead{e1}}
   \author{\fnms{Oliver} \snm{Rinne}%
     \ead[label=e2]{oliver.rinne@htw-berlin.de}}
  \address{Faculty 4, HTW Berlin---University of Applied Sciences,\\ 
    Treskowallee 8, 10318 Berlin, Germany, and\\
    Max Planck Institute for Gravitational Physics (Albert Einstein 
    Institute),\\Am M\"uhlenberg 1, 14476 Potsdam, Germany\\
    \printead{e2}}
  \and
  \author{\fnms{Markus} \snm{Strehlau}%
    \ead[label=e3]{markus.strehlau@b-tu.de}}
  \address{Institute of Mathematics, Brandenburgische Technische 
    Universit\"at\\Cottbus-Senftenberg,
    Postfach 10 13 44, 03013 Cottbus, Germany, and\\
    Max Planck Institute for Gravitational Physics (Albert Einstein 
    Institute),\\Am M\"uhlenberg 1, 14476 Potsdam, Germany\\
    \printead{e3}}
\end{aug}

\begin{abstract}
We investigate Bartnik's static metric extension conjecture under the
additional assumption of axisymmetry of both the given Bartnik data
and the desired static extensions. To do so, we suggest a geometric
flow approach, coupled to the Weyl--Papapetrou formalism for
axisymmetric static solutions to the Einstein vacuum equations. The
elliptic Weyl--Papapetrou system becomes a free boundary value problem
in our approach. We study this new flow and the coupled flow--free
boundary value problem numerically and find axisymmetric static
extensions for axi\-symmetric Bartnik data in many situations,
including near round spheres in spatial Schwarzschild of positive
mass.\\

This paper is dedicated to Robert Bartnik on the occasion of his 60th birthday. Happy Birthday, Robert!
 \end{abstract}

\begin{keyword}
\kwd{General Relativity}
\kwd{Axisymmetry}
\kwd{Weyl--Papapetrou Coordinates}
\kwd{Geometric Flow}
\kwd{Free Boundary Value Problem}
\end{keyword}


\end{frontmatter}

\input{paper_sec1}
\input{paper_sec2}
\input{paper_sec3}
\input{paper_sec4}
\input{paper_sec5}
\input{paper_sec6}

\section*{Acknowledgements}
The authors would like to thank 
Michael Anderson, 
Olaf Baake,
Armando J.~Cabrera Pacheco, 
Friederike Dittberner, 
Georgios Doulis,
Leon Escobar,
Gerhard Huisken, 
Claus Kiefer,
Heiko Kr\"oner,
Elena M\"ader-Baumdicker,
Claudio Paganini,
and Christian Schell 
for stimulating discussions and for asking helpful questions.

CC is indebted to the Baden-W\"urttemberg Stiftung for the financial
support of this research project by the Eliteprogramme for Postdocs. 
Work of CC is supported by the Institutional Strategy of the
University of T\"ubingen (Deutsche Forschungsgemeinschaft, ZUK 63). 
The early stages of OR's work on this project were supported by a
Heisenberg Fellowship and grant RI 2246/2 of the Deutsche
Forschungsgemeinschaft.
MS ist indebted to the Graduate Research School (GRS) of the
Brandenburg University of Technology Cottbus--Senftenberg for
financial support.

The authors would like to extend thanks to the Albert Einstein
Institute for technical support and for allowing us to collaborate in
a stimulating environment.

\input{paper.bbl}
\vfill

\end{document}

%% file: paper_sec1.tex
\section{Introduction}
\label{s:intro}
In~\cite{Bartnik89}, Robert Bartnik introduced within the theory of General Relativity a new notion of quasi-local ``mass'' or ``capacity'' for bounded spatial regions in an initial data set in a given spacetime. This definition, which is now referred to as the \emph{Bartnik mass}, is given as an infimum over the ADM masses of all ``admissible'' asymptotically flat initial data set extensions of the given bounded region -- with no reference to the spacetime in which the region is contained to begin with. Bartnik then conjectured that the infimum should be attained by a stationary, vacuum, asymptotically flat initial data set that attaches to the given bounded region in a suitably regular manner. This leads to the related question of whether or not such stationary, vacuum, asymptotically flat initial data sets extending the given region in a suitably regular manner will generically exist. This question of existence of stationary ``extensions'' of bounded spatial regions remains open until today to the best knowledge of the authors.

More is known when one restricts to the time-symmetric (or Riemannian) case, as we will do here. In the time-symmetric context, a \emph{bounded spatial region} is described by a smooth compact Riemannian $3$-manifold $(\Omega,\breve{\gamma})$ with non-empty boundary $\partial\Omega\neq\emptyset$. For simplicity and definiteness, we will assume that this boundary is diffeomorphic to $\mathbb{S}^{2}$, $\partial\Omega\approx\mathbb{S}^{2}$, has positive Gaussian curvature $K>0$, and positive mean curvature $H>0$ (with respect to the outward pointing unit normal)\footnote{Our convention for the mean curvature is such that the round spheres of radius $r$ in Euclidean $3$-space will have mean curvature $H=\frac{2}{r}$ with respect to the outward unit normal.}. Furthermore, we assume that the scalar curvature $\operatorname{R}_{\breve{\gamma}}\geq0$ is non-negative, or in other words that the (Riemannian) dominant energy condition is satisfied. In the time-symmetric setting, the question of existence of stationary extensions reduces to a question that is known as \emph{Bartnik's static metric extension conjecture}: Given a bounded spatial region $(\Omega,\breve{\gamma})$ as described above, does there always exist an asymptotically flat Riemannian $3$-manifold $(\mathcal{M},\gamma)$, called the \emph{static metric extension}, such that $(\Omega,\breve{\gamma})\hookrightarrow(\mathcal{M},\gamma)$ isometrically, $(\mathcal{M},\gamma)$ is smooth except possibly across $\partial\Omega$, and is (standard) static vacuum in the sense that there exists a smooth \emph{lapse function} $N\colon \mathcal{M}\setminus\Omega\to\mathbb{R}^{+}$ with $N\to1$ suitably fast in the asymptotic end, so that the \emph{static vacuum Einstein equations}
\begin{align}\label{e:SVE1}
N\operatorname{Ric}_{\gamma}&=\nabla_{\gamma}^{2}N\\\label{e:SVE2}
\triangle_{\gamma}N&=0
\end{align}
hold on $\mathcal{M}\setminus\Omega$. Here,
$\operatorname{Ric}_{\gamma}$ denotes the Ricci curvature tensor of 
$\gamma$, and $\nabla_{\gamma}^{2}$ and $\triangle_{\gamma}$ denote the Hessian and the Laplacian with respect to $\gamma$, respectively. Note that the static vacuum Einstein equations~\eqref{e:SVE1},~\eqref{e:SVE2} imply scalar flatness of $(\mathcal{M}\setminus{\Omega},\gamma)$,  $\operatorname{R}_{\gamma}=0$, such that the (Riemannian) dominant energy condition is automatically
satisfied in the extension $(\mathcal{M},\gamma)$, away from
$\partial\Omega$. Furthermore, one requests that
$(\mathcal{M},\gamma)$ be regular enough across $\partial\Omega$
so that the scalar curvature of $\gamma$ can be assumed to be distributionally non-negative. 

Depending on the precise definition of Bartnik mass one uses, additional conditions will need to be requested of the static extension in order to connect the static metric extension problem to the search of a minimizer of Bartnik's quasi-local mass in the time-symmetric context. One such condition would be that $\partial\Omega$ needs to be area outer minimizing in $(\mathcal{M}\setminus\Omega,{\gamma})$ or that there shall be no minimal surfaces in $(\mathcal{M}\setminus\Omega,{\gamma})$ (homologous to $\partial\Omega$).\\

It has become customary to study the following simplified ``boundary version'' of Bartnik's static metric extension conjecture, replacing the isometric embedding condition with suitable regularity across $\partial\Omega$ by a boundary condition compatible with the distributional non-negativity condition on the scalar curvature. This is the conjecture we will address in this paper.
\begin{conj}[Bartnik's static metric extension conjecture, boundary version]\label{conj:Bartnik}
Let $(\Sigma\approx\mathbb{S}^{2},g)$ be a smooth Riemannian $2$-manifold with positive Gaussian curvature, and let $H\colon\Sigma\to\mathbb{R}^{+}$ be a smooth positive function. The tuple $(\Sigma,g,H)$ is called \emph{Bartnik data}. Then there conjecturally is a smooth Riemannian $3$-manifold $(\mathcal{M},\gamma)$ with boundary $\partial \mathcal{M}$ and a smooth, positive \emph{lapse function} $N\colon \overline{\mathcal{M}}\to\mathbb{R}^{+}$ such that the \emph{static system} $(\mathcal{M},\gamma,N)$
\begin{enumerate}
\item satisfies the \emph{static vacuum Einstein equations}
\begin{align*}
N\operatorname{Ric}_{\gamma}&=\nabla_{\gamma}^{2}N\\
\triangle_{\gamma}N&=0,
\end{align*}\item is \emph{asymptotically flat}, i.e.~there exists a smooth diffeomorphism $\varphi=(x^{i})\colon \mathcal{M}\setminus\mathcal{K}\to\mathbb{R}^{3}\setminus B$, with $B\subset\mathbb{R}^{3}$ some bounded, open ball, $\mathcal{K}\subset \mathcal{M}$ is compact, and
\begin{align*}
(\varphi_{*}\gamma)_{ij}&=\delta_{ij}+\mathcal{O}_{2}(r^{-1}),\\
\varphi_{*}N&=\phantom{_{ij}}1+\mathcal{O}_{2}(r^{-1})
\end{align*}
as $r\definedas\sqrt{(x^{1})^{2}+(x^{2})^{2}+(x^{3})^{2}}\to\infty$, where $\delta_{ij}$ denotes the Euclidean metric on $\mathbb{R}^{3}\setminus B$,
\item and has inner boundary $(\partial \mathcal{M},\gamma\vert_{\partial \mathcal{M}})$ isometric to $(\Sigma,g)$ with induced mean curvature $H$ with respect to the unit normal pointing to the asymptotic end in $\mathcal{M}$. 
\end{enumerate}
If $(\mathcal{M},\gamma,N)$ exists, we call it a \emph{static metric extension of $(\Sigma,g,H)$}.
\end{conj}

\begin{rems}
Let us make the following remarks.
\begin{itemize}
\item We do not request any ``outward minimizing property'' nor any ``no minimal surfaces condition''. We do, however, consistently with either of those assumptions, assume that the lapse function $N$ be positive. 
\item We do not explicitly request that there be a fill-in $(\Omega,\gamma)$ with boundary $\partial\Omega\approx\Sigma$, such that $(\partial\Omega,\breve{\gamma}\vert_{\partial\Omega})$ is isometric to $(\Sigma,g)$ and has mean curvature $H$. For a more thorough discussion on fill-ins, see e.g.~\cite{Jauregui:2011fu}. 
\item We do not make any claim about uniqueness of the static metric extension $(\mathcal{M},\gamma,N)$.
\end{itemize}
\end{rems}

Clearly, if $(\mathcal{M},\gamma,N)$ \emph{is} a static system satisfying (1), (2) above and is such that $\partial \mathcal{M}$ is diffeomorphic to $\mathbb{S}^{2}$, has positive Gaussian curvature $K>0$, and positive mean curvature $H>0$ with respect to the unit normal pointing to the asymptotically flat end, then its induced Bartnik data $(\partial \mathcal{M},g,H)$ naturally possess the static metric extension $(\mathcal{M},\gamma,N)$. Of course, we do not know if this is the only static metric extension of $(\partial \mathcal{M},g,H)$. Neither do we know of an explicit method of reconstructing $(\mathcal{M},\gamma,N)$ from the Bartnik data $(\partial \mathcal{M},g,H)$ in general. 

The static metric extension conjecture becomes much simpler, and indeed a priori resolved, once one restricts one's attention to \emph{spherically symmetric Bartnik data $(\Sigma,g,H)$}\label{p:sphericalsymmetry}, i.e.~to the case where $(\Sigma,g)$ is round, meaning isometric to $(\mathbb{S}^{2},R^{2}\sigma)$ with some radius $R>0$ and $\sigma$ denoting the canonical metric on the unit sphere, and $H>0$ is a constant. Such spherically symmetric Bartnik data are always extended by the well-known \emph{Schwarzschild static system of mass $M_*$}, more precisely by $(\mathcal{M},\gamma,N)$ given by
\begin{align}
\begin{split}\label{e:schwarzschild}
\mathcal{M}&=\;\left(R,\infty\right)\;\times\;\mathbb{S}^{2},\\
\gamma&=N^{-2}dr^{2}+r^{2}\sigma,\\
N&=N(r)=\sqrt{1-\frac{2M_*}{r}},
\end{split}
\end{align} 
where the mass can be picked as the ``Hawking mass'' of the Bartnik data,
\begin{align}\label{Smass}
M_*&\definedas m_\mathrm{H}=\frac{R}{2}\left(1-\frac{H^{2}R^{2}}{4}\right),
\end{align}
see below. From~\eqref{Smass}, one recovers the Euclidean case $H=\frac{2}{R}$, where $N\equiv1$ and $M_*=0$. It is well-known~\cite{Schwarzschild} that the Schwarzschild static systems are the only spherically symmetric solutions of the static vacuum Einstein equations~\eqref{e:SVE1},~\eqref{e:SVE2}. A simple computation shows that the mass $M_*$ computed in~\eqref{Smass} is the only one that induces the given Bartnik data. In this sense, one can say that static metric extensions are ``unique in the category of spherically symmetric extensions''. Again, we do not know in general if there will be an additional, non-spherically symmetric static metric extension of given spherically symmetric Bartnik data.

Using a subtle implicit function theorem argument, Miao~\cite{Miao2003} showed that, given Bartnik data $(\Sigma,g,H)$ that are close to Euclidean unit round sphere Bartnik data $(\mathbb{S}^{2},g=\sigma,H=2)$ in a suitable Sobolev norm, and that possess a certain $\mathbb{Z}_{2}\times\mathbb{Z}_{2}\times\mathbb{Z}_{2}$-symmetry, there exists a static metric extension close to the Euclidean static system $(\left(1,\infty\right)\times\mathbb{S}^{2},\gamma=dr^{2}+r^{2}\sigma,N=1)$ in a suitably weighted Sobolev norm. This result easily generalises to Bartnik data near round sphere data of arbitrary radius~\cite{Piubello}. Later, this result was generalised by Anderson~\cite{Anderson} for general perturbations of the flat background. A related result by the first author can heuristically be stated as saying that the number of (functional) degrees of freedom of Bartnik data coincides with the number of (functional) degrees of freedom of asymptotically flat static vacuum systems, see~\cite[Sec. 3.4]{Cederbaum2011} for details. Shortly thereafter, Anderson--Khuri~\cite{Anderson_2013} addressed the question of degrees of freedom using methods from functional analysis and showed that the correspondence between static vacuum solutions and Bartnik boundary data is Fredholm of index $0$. In other words, the linearisation of this correspondence has finite-dimensional kernel and co-kernel at every point.

In this paper, we will address Bartnik's static metric extension conjecture in the form stated above under the additional assumption that both the Bartnik data $(\Sigma,g,H)$ and the desired static metric extensions $(\mathcal{M},\gamma,N)$ be ``compatibly axisymmetric'' in the following sense.

\begin{define}[Axisymmetric Bartnik data and extensions]\label{def:axisymm}
Let $(\Sigma,g,H)$ be Bartnik data. We say that $(\Sigma,g,H)$ is \emph{axisymmetric} if there is a Killing vector field $X$ on $(\Sigma,g)$, $\mathcal{L}_{X}g=0$, with closed orbits that keeps the mean curvature $H$ invariant in the sense that
\begin{align}
X(H)=0.
\end{align}
Now let $(\mathcal{M},\gamma,N)$ be a static metric extension of Bartnik data $(\Sigma,g,H)$ that are axisymmetric with respect to some field $X$ on $\Sigma$. We say that $(\mathcal{M},\gamma,N)$ \emph{is a compatibly axisymmetric static metric extension of $(\Sigma,g,H)$} or more sloppily \emph{is axisymmetric} if $X$ extends to a smooth Killing vector field $\widehat{X}$ of $(\mathcal{M},\gamma)$, $\mathcal{L}_{\hat X} \gamma = 0$, with closed orbits that keeps $N$ invariant in the sense that 
\begin{align}
\widehat{X}(N)=0.
\end{align}
\end{define}

Naturally, the spherically symmetric situation discussed above is a special case of this axisymmetric setup. Again, we only address the question of whether, given axisymmetric Bartnik data, there exists a compatibly axisymmetric static metric extension, and make no assertions about uniqueness nor about (in)existence of non-axisymmetric extensions. We thus address the following conjecture which speaks about a smaller category but voices a stronger expectation.

\begin{conjecture}[Bartnik's static metric extension conjecture in axisymmetry]\label{conj:Bartnikaxisymm}
Let $(\Sigma,g,H)$ be axisymmetric Bartnik data. Then there conjecturally exists a compatibly axisymmetric static metric extension $(\mathcal{M},\gamma,N)$.
\end{conjecture}
\newpage
To the best knowledge of the authors, the notion of axisymmetry has not yet been considered before in this context. We do not know of any reasons derived from the original derivation of Bartnik's static metric extension conjecture that leads to the expectation that axisymmetric Bartnik data should possess axisymmetric extensions in general. This may in fact be an interesting question to study in its own right. However, as we will see in Section~\ref{s:formulation}, the axisymmetric setup suggested here allows for the introduction of ideas and tools that are very different from those that have been used in the generic scenario and that allow us to achieve at least some positive numerical results.

In practice, we will make an additional assumption of compatible reflection symmetry of the Bartnik data and the desired static metric extensions in order to simplify our analysis of Conjecture~\ref{conj:Bartnikaxisymm} across a plane orthogonal to the axis of rotation, see Convention~\ref{conv:reflection}. Altogether, our symmetry assumptions imply the $\mathbb{Z}_{2}\times\mathbb{Z}_{2}\times\mathbb{Z}_{2}$-symmetry assumption in Miao's work~\cite{Miao2003}. Again, we do not know of any reason why a static metric extension of reflection symmetric Bartnik data should a priori be reflection symmetric.

\paragraph*{Strategy.} To discuss Conjecture~\ref{conj:Bartnikaxisymm}, we will draw on the work of Weyl
\cite{Weyl1917} and Papapetrou~\cite{Papapetrou} and adopt global quasi-isotropic coordinates $(r,\theta,\varphi)$ adapted to the
axisymmetry of the desired static metric extension $(\mathcal{M},\gamma,N)$. In these so-called Weyl--Papapetrou coordinates, the static vacuum equations~\eqref{e:SVE1},~\eqref{e:SVE2} will take on a particularly simple form~\eqref{e:Ueqn},~\eqref{e:Veqns} which we will call the Weyl--Papapetrou equations. In particular, the coordinate $\varphi$ corresponding to the axial
Killing vector field $\widehat{X}=\partial_{\varphi}$ will drop out and the equations as well as the remaining two scalar function
variables $(U,V)$ subject to the Weyl--Papapetrou equations will be stated in a symmetry-reduced form in the $(r,\theta)$-half-plane orthogonal to the axis of symmetry. Due to the compatible axisymmetry,
the Bartnik data will symmetry-reduce to a curve~$\Gamma$ in this
half-plane which will take the role of a free boundary for the
Weyl--Papapetrou equations as we do not know the values of the
coordinates $(r,\theta)$ along $\Gamma$ a priori.

We will then approach Conjecture~\ref{conj:Bartnikaxisymm} as follows:
Interpret given axisymmetric Bartnik data $(\Sigma,g,H)$ as a free
boundary curve $\Gamma$ in the half-plane orthogonal to the axis of
symmetry of a potential compatibly axisymmetric static metric
extension. The metric $g$ then prescribes
Dirichlet boundary values for the free field $U$ along $\Gamma$ (and $V$ can be obtained from $U$ by integration). 
The mean curvature $H$ takes the role of a consistency condition to ensure that the free boundary is located correctly.
The asymptotic decay conditions prescribe boundary values at infinity, while smoothness requirements provide additional compatibility conditions along the axis. \emph{If we knew} the position of the curve $\Gamma$ in the Weyl--Papapetrou coordinates of the desired extension it would in fact be straightforward to solve the Weyl--Papapetrou equations and compute the solution $(U,V)$, allowing us to reconstruct the desired compatibly axisymmetric static metric extension~$(\mathcal{M},\gamma,N)$.

However, we do of course \emph{not} know $\Gamma$ in Weyl--Papapetrou
coordinates a priori. We thus adopt a strategy that couples solving the Weyl--Papapetrou equations to a geometric flow of the boundary curve: We guess an
initial curve $\Gamma$ in Weyl--Papapetrou coordinates as well as an
initial solution $(U,V)$ of the Weyl--Papapetrou equations compatible
with the asymptotic decay conditions and the compatibility
requirements along the axis. For these initial guesses, we do
\emph{not} request that the geometry of the Bartnik data be consistent
with the inner boundary values of $(U,V)$. Then, we deform $\Gamma$ by
a geometric flow in the background geometry of the half-plane 
initially
given by the initial $(U,V)$; the chosen flow naturally depends on the
geometry of the given Bartnik data. 
At the same time, we compute the boundary values for $(U,V)$ induced 
by the geometry of the given Bartnik data on the current (flowing) 
boundary curve and solve the Weyl--Papapetrou equations to 
  update the fields $(U, V)$. The flow is chosen such that the flowing curve
$\Gamma_{t}$ shall ideally approach the ``true'' position of the
boundary curve in case this true position exists in the 
constantly updated
background described by $(U,V)$.

\paragraph*{The paper is structured as follows:} We will first remind the reader very briefly of some notions from Mathematical Relativity that will be used to formulate our approach and to analyse the numerical results. Then, in Section~\ref{s:formulation}, we formulate our approach to the axisymmetric version of Bartnik's static metric extension conjecture~\ref{conj:Bartnikaxisymm}, including the definition of the geometric flow we use and the derivation of the free boundary value Weyl--Papapetrou problem. In Section~\ref{s:analysis.existence}, we study some analytic and geometric properties of the flow in a flat background. Then, in Section~\ref{s:nummethod}, we introduce and describe our numerical schemes for the geometric flow and the free boundary value Weyl--Papapetrou system. We present our numerical results in Section~\ref{s:numresults} and discuss them in Section~\ref{s:discussion}.

\subsection{Some notions from Mathematical Relativity}\label{s:masses}
In the numerical analysis of the geometric flow and the combined free
boundary value approach we suggest, we will use several notions of
total and quasi-local mass to gain some insight into the numerical
solutions. We will briefly discuss these notions and their relevant
properties here, adjusted to the context of asymptotically flat
solutions $(\mathcal{M},\gamma,N)$ of the static vacuum Einstein equations
\eqref{e:SVE1},~\eqref{e:SVE2}.

In this context, the (total) 
\emph{Arnowitt--Deser--Misner (ADM)} mass~\cite{ADM}, $m_\mathrm{ADM}$,  is given by
\begin{align}\label{def:mADM}
m_\mathrm{ADM}&\definedas \frac{1}{16\pi}\lim_{r\to\infty}\int_{\mathbb{S}^{2}_{r}} \left(\gamma_{ii,j}-\gamma_{ij,i}\right)\frac{x^{j}}{r}\,dA,
\end{align}
where $dA$ denotes the (Euclidean) area element on
$\mathbb{S}^{2}_{r}$, $\gamma_{ij}$ is short for
$(\varphi_{*}\gamma)_{ij}$, 
$r^2 = \delta_{ij} x^i x^j$,
and commas denote partial (coordinate) derivatives. The ADM mass is well-defined as the scalar curvature vanishes due to the static vacuum Einstein equations~\eqref{e:SVE1},~\eqref{e:SVE2} and because of the asymptotic decay conditions we assume, see~\cite{Bartnik86,Chrus2}. 

We will numerically compute two quasi-local masses, namely the Hawking mass $m_\mathrm{H}$~\cite{Hawking} already alluded to above and the pseudo-Newtonian mass $m_\mathrm{PN}$~\cite{Cederbaum2011}, which is genuinely only defined in the static realm. The \emph{Hawking mass} of Bartnik data $(\Sigma,g,H)$ is given by 
\begin{align}\label{def:mH}
m_{H}&\definedas \frac{\sqrt{\vert\Sigma\vert_{g}}}{16\pi}\left(1-\frac{1}{16\pi}\int_{\Sigma}H^{2}\,dA_{g}\right).
\end{align}
Here, $\vert\Sigma\vert_{g}$ denotes the area of $\Sigma$ and $dA_{g}$ the area element with respect to~$g$. It follows from Huisken--Ilmanen's proof of the Penrose inequality~\cite{HI} that
\begin{align}\label{genPI}
m_\mathrm{H}&\leq m_\mathrm{ADM}
\end{align}
holds for the Hawking mass of any Bartnik data $(\Sigma,g,H)$ sitting
inside an asymptotically flat static system $(\mathcal{M},\gamma,N)$ of ADM mass
$m_\mathrm{ADM}$ in an \emph{area outer minimizing} way, meaning that
any $2$-surface $\widetilde{\Sigma}\hookrightarrow(\mathcal{M},\gamma)$
homologous to $\Sigma$ with induced metric $\widetilde{g}$ will have
area at least as big as that of $\Sigma$,
$\vert\Sigma\vert_{g}\leq\vert\widetilde{\Sigma}\vert_{\widetilde{g}}$. We
will make use of this \emph{generalised Penrose inequality}
\eqref{genPI} to check consistency of our numerical results, see
Sections~\ref{s:numresults} and~\ref{s:discussion}.

The \emph{pseudo-Newtonian mass} $m_\mathrm{PN}$ of Bartnik data $(\Sigma,g,H)\hookrightarrow(\mathcal{M},\gamma,N)$ in an asymptotically flat static system $(\mathcal{M},\gamma,N)$ is defined in~\cite{Cederbaum2011} as
\begin{align}\label{def:mPN}
m_\mathrm{PN}&\definedas \frac{1}{4\pi}\int_{\Sigma}\nu(N)\,dA_{g},
\end{align}
where $\nu$ denotes the unit normal to $\Sigma$ in $(\mathcal{M},\gamma)$
pointing to the asymptotically flat end. A straightforward computation
shows that the notion of pseudo-Newtonian mass in fact
coincides with that of Komar mass~\cite{Komar}. Moreover, it follows
from the divergence theorem that the pseudo-Newtonian mass
$m_\mathrm{PN}$ of Bartnik data
$(\Sigma,g,H)\hookrightarrow(\mathcal{M},\gamma,N)$ indeed coincides with the
ADM mass $m_\mathrm{ADM}$ of the surrounding static system
$(\mathcal{M},\gamma,N)$ if $(\mathcal{M},\gamma,N)$ solves the static vacuum Einstein
equations~\eqref{e:SVE1},~\eqref{e:SVE2}:
\begin{align}\label{e:massidentity}
m_\mathrm{PN}&=m_\mathrm{ADM},
\end{align}
see~\cite[Chapt.~4]{Cederbaum2011}. This fact will also be used to
check consistency of our numerical results, see
Sections~\ref{s:numresults} and~\ref{s:discussion}.

%% file: paper_sec2.tex
\section{Formulation of the problem}
\label{s:formulation}

From this section onwards, we will overline all functions and tensor fields corresponding to or induced by prescribed Bartnik data in order to distinguish them from fields of the same geometric kind that we are flowing or otherwise computing. For example, Bartnik data themselves will from now on be denoted by~$(\Sigma,\overline{g},\overline{H})$.

Recall that a given axisymmetric surface $(\Sigma\approx\mathbb{S}^{2},\overline{g})$ can be rewritten in terms of an arclength parametrisation of its rotational profile as
\begin{align}\label{e:2dmetricintrinsic}
\overline{g} &= d\tau^2 + \overline{\lambda}^2 d\phi^2,
\end{align}
where $\tau\in[0,\overline{L}]$ denotes the arclength coupling parameter,
$\phi\in[0,2\pi)$ the angle of rotation, $\overline{L}$ the total length of the rotation profile, and $\overline{\lambda}=\overline{\lambda}(\tau)$ is a function induced by $\overline{g}$ determining the intrinsic geometry of $(\Sigma,\overline{g})$. Accordingly, a given function $\overline{H}\colon\Sigma\to\mathbb{R}$, can be understood as a function $\overline{H}\colon[0,\overline{L}]\to\mathbb{R}$, slightly abusing notation. We will pursue this perspective for Bartnik data $(\Sigma,\overline{g},\overline{H})$ throughout the remainder of this work.

\subsection{Axisymmetric static systems}\label{s:formulation.saveinstein}
Let us now consider the axisymmetric version of Bartnik's static metric extension conjecture, Conjecture~\ref{conj:Bartnikaxisymm}, using the global ansatz 
\begin{align}\label{e:metric}
  \gamma &= e^{-2U} [ e^{2V} (dr^2 + r^2 d\theta^2) + r^2 \sin^2\theta \, d\phi^2 ]\\\label{e:lapse}
  N&= e^{U}
\end{align}
in global Weyl--Papapetrou coordinates $(r,\theta,\phi)$ for the
axisymmetric metric $\gamma$ and lapse function $N$ of a static system
$(\mathcal{M},\gamma,N)$, where $\partial_\phi=\widehat{X}$ denotes the axial
Killing vector field. This ansatz goes back to Weyl~\cite{Weyl1917}
and Papapetrou~\cite{Papapetrou}. We identify the manifold $\mathcal{M}$ with a
domain $\Omega\times[0,2\pi)$ with free inner boundary of the
coordinate range $\left(\mathbb{R}^{+}\times[0,\pi]\right)\times[0,2\pi)$.
Here, the free functions $U=U(r,\theta)$ and $V=V(r,\theta)$,
$(r,\theta)\in\Omega\subset\mathbb{R}^{+}\times[0,\pi]$, denote smooth
real valued functions that contain all the geometric information of the given static system. Using this ansatz, we obtain the following standard formula for the \emph{length $\lambda$} of the axial Killing vector field $\partial_\phi=\widehat{X}$\begin{align}\label{e:lambda}
  \lambda &= \sqrt{\gamma_{\phi\phi}} = e^{-U} r \sin\theta.
\end{align}

The static vacuum Einstein equations~\eqref{e:SVE1},~\eqref{e:SVE2} for a metric and lapse of the form~\eqref{e:metric},~\eqref{e:lapse} reduce to the \emph{Weyl--Papapetrou equations}
\begin{align}  \label{e:Ueqn}
  \Delta_{\delta}\, U &= U_{,rr} + \frac{2}{r} U_{,r} + \frac{1}{r^2} 
    (U_{,\theta\theta} + \cot \theta \, U_{,\theta}) = 0,\\  
   \begin{split}\label{e:Veqns}
  V_{,r} &= r \sin^2\theta \, U_{,r}^2 + 2 \sin\theta\cos\theta\, U_{,r} U_{,\theta}
  - \frac{\sin^2\theta}{r} U_{,\theta}^2,\\ 
  V_{,\theta} &= -r^2\sin\theta\cos\theta\, U_{,r}^2 + 2r\sin^2\theta\, U_{,r} U_{,\theta}
  + \sin\theta\cos\theta \, U_{,\theta}^2.
  \end{split}
\end{align}
Here, $\Delta_{\delta}$ denotes the $3$-dimensional Euclidean Laplacian in spherical polar coordinates. Observe that the first equation,~\eqref{e:Ueqn}, is a standard (Euclidean) Laplace equation for $U$ and thus linear second-order elliptic. It is decoupled from the second set of equations and in particular does not depend on $V$. The second set of equations,~\eqref{e:Veqns}, gives the first partial derivatives of $V$ in terms of (first partial derivatives of) $U$.

To incorporate the asymptotic flatness condition imposed on static metric extensions, we furthermore assume that the Weyl--Papapetrou coordinates are consistent with the asymptotic flatness assumptions in the sense that the asymptotic coordinates $(x^{i})$ can be chosen such that they are the Cartesian coordinates corresponding to the spherical polar coordinates $(r,\theta,\phi)$. This leads to the asymptotic decay conditions
\begin{align}\label{e:decay}
U,V&=\mathcal{O}_{2}(r^{-1}) \text{ as }r\to\infty.
\end{align}

Now assume we are given a smooth, $2$-dimensional, compatibly axisymmetric surface $(\Sigma,\overline{g})$ isometrically sitting inside a static system $(\mathcal{M},\gamma,N)$ of the form~\eqref{e:metric},~\eqref{e:lapse} that satisfies~\eqref{e:Ueqn},~\eqref{e:Veqns}. Clearly, $\Sigma$ can be described in Weyl--Papapetrou coordinates by a curve $\Gamma=(r,\theta)\colon I\to \mathbb{R}^{+}\times[0,\pi]$, where $I=[0,b]$ is some interval. For smoothness reasons, the curve $\Gamma$ has to stay away from the axis $\lbrace{z\definedas r\cos\theta=0\rbrace}$ except at the endpoints, where it needs to be horizontal. This horizontality condition is equivalent to requesting that $\Gamma$ satisfies the boundary conditions
\begin{align}
\begin{split}\label{e:smooth}
\theta(0)&=0,\quad \theta(b)=\pi,\\
r'(0)&=0,\quad r'(b)=0.
\end{split}
\end{align}

If $\Gamma=\Gamma(\tau)=(r(\tau),\theta(\tau))$ is parametrised by arclength on $I=[0,\overline{L}]$ with arclength parameter $\tau$ and total length $\overline{L}$ determined by $(\Sigma,\overline{g})$, the metric can  --- by compatibility of the intrinsic axisymmetry and the axisymmetry of the surrounding static system --- be written as
\begin{align}\label{e:2dmetric}
\overline{g} &= d\tau^2 + \lambda^2 d\phi^2,
\end{align}
 with $\lambda=\lambda(\tau)$, or equivalently, in view of~\eqref{e:2dmetricintrinsic},
\begin{align}\label{e:lambdaequal}
\overline{\lambda}(\tau)=\lambda(\tau),
\end{align}
which
ensures that the embedding $(\Sigma,\overline{g})\hookrightarrow(\mathcal{M},\gamma,N)$ is indeed isometric, for $\tau\in[0,\overline{L}]$, $\phi\in[0,2\pi)$, see~\eqref{e:2dmetricintrinsic}. The restriction of $U$ to the surface $\Sigma$ can then be expressed as 
\begin{align}  \label{e:Ubc}
  U \circ \Gamma (\tau) =-\ln \frac{\overline{\lambda}(\tau)}
   {r(\tau) \sin\theta(\tau)}
\end{align}
for $\tau\in[0,\overline{L}]$ via~\eqref{e:lambda}. The condition that the surface $(\Sigma,\overline{g})\hookrightarrow(\mathcal{M},\gamma,N)$ be isometrically embedded, or equivalently that the curve parameter $\tau$ indeed be the arclength parameter, can also be stated as
\begin{align}\label{e:embedding}
  \ell^2 &\definedas e^{2(V-U)\circ\Gamma} (r'^2 + r^2\theta'^2) \equiv 1
\end{align}
on $[0,\overline{L}]$, where a prime denotes a derivative w.r.t.~$\tau$. For notational simplicity, we will from now on drop $\circ\Gamma$ after $U$, $V$, etc. in expressions such as the one in~\eqref{e:embedding} and hope that no confusion will arise from this.

Taking a $\tau$-derivative of~\eqref{e:embedding} and using~\eqref{e:Veqns} to replace derivatives of $V$ with derivatives of $U$, we obtain the identity
\begin{align}  
\begin{split}\label{e:C}
  0 = C \definedas  \ell^{-2} \ell' =\;& e^{2(V-U)} \ell^{-3} (r' r'' + r r' \theta'^2 + r^2 \theta' \theta'')    \\
     &+ \ell^2 \Big[ - r' U_{,r} - \theta' U_{,\theta}\\
     &\qquad + (r \sin\theta\, U_{,r}^2 - \frac{\sin\theta}{r} U_{,\theta}^2 )
       (r'\sin\theta - r \theta' \cos\theta)\\
     &\qquad + 2\sin\theta\, U_{,r} U_{,\theta} 
     (r'\cos\theta + r \theta' \sin\theta) \Big].
\end{split}
\end{align}
We will use this identity in the definition of the geometric flow below.

Finally, a straightforward computation shows that the induced mean curvature~$H$ of $(\Sigma,\overline{g})\hookrightarrow(\mathcal{M},\gamma,N)$ with respect to the unit normal pointing to the asymptotically flat end can be expressed as a function $H=H(\tau)$ by 
\begin{align}
\begin{split}\label{e:H}
  H =\;& e^{2(V-U)} 
    \ell^{-3} (-r r'' \theta' + 2 r'^2 \theta' + r r' \theta'' + r^2 \theta'^3)\\
    & + \ell^{-1} \Big[ - \frac{r'}{r} \cot\theta + \theta' 
    + 2\left(\frac{r'}{r} U_{,\theta} - r \theta' U_{,r}\right)\\
    & \qquad + \left(r\sin\theta\, U_{,r}^2 
      - \frac{\sin\theta}{r} U_{,\theta}^2\right)
      (r'\cos\theta + r \theta'\sin\theta)\\
    & \qquad - 2 \sin\theta\, U_{,r} U_{,\theta} 
    (r'\sin\theta - r \theta'\cos\theta)  \Big].
 \end{split}
\end{align}
This expression is valid even if $\tau$ is a parameter different from arclength.

\subsection{Static metric extensions in Weyl--Papapetrou form}\label{s:formulation.metricext}
We can thus rephrase Conjecture~\ref{conj:Bartnikaxisymm} as follows.
\begin{conjecture}[Bartnik's static metric extension conjecture in axisymmetry in Weyl--Papapetrou form]\label{conj:Bartnikaxisymm2}
Let $(\Sigma,\overline{g},\overline{H})$ be axisymmetric Bartnik data. Then there conjecturally exists a domain $\Omega\subset\mathbb{R}^{+}\times[0,\pi]$ \emph{containing the asymptotically flat end} in the sense that $(r_{0},\infty)\times[0,\pi]\subseteq\Omega$ for some $r_{0}>0$, and a solution $(U,V)$ of the Weyl--Papapetrou equations
\begin{align*}
  \Delta_{\delta}\, U &= U_{,rr} + \frac{2}{r} U_{,r} + \frac{1}{r^2} 
    (U_{,\theta\theta} + \cot \theta \, U_{,\theta}) = 0,\\  
   \begin{split}
  V_{,r} &= r \sin^2\theta \, U_{,r}^2 + 2 \sin\theta\cos\theta\, U_{,r} U_{,\theta}
  - \frac{\sin^2\theta}{r} U_{,\theta}^2,\\ 
  V_{,\theta} &= -r^2\sin\theta\cos\theta\, U_{,r}^2 + 2r\sin^2\theta\, U_{,r} U_{,\theta}
  + \sin\theta\cos\theta \, U_{,\theta}^2
  \end{split}
\end{align*}
 on $\Omega$ with asymptotic decay 
 \begin{align*}
U,V&=\mathcal{O}_{2}(r^{-1})\text{ as }r\to\infty
\end{align*}
such that the inner boundary $\partial\Omega$ of the domain $\Omega$ can be written as a curve $\Gamma=(r,\theta)\colon[0,\overline{L}]\to\mathbb{R}^{+}\times[0,\pi]$, $\partial\Omega=\Gamma([0,\overline{L}])$, parametrised by arclength,
\begin{align*}
\ell&=e^{2(U-V)}(r'^{2}+r^{2}\theta'^{2})\equiv1,
\end{align*}
satisfying the boundary conditions
\begin{align*}
\theta(0)&=0,\quad \theta(\overline{L})=\pi,\\
r'(0)&=0,\quad r'(\overline{L})=0,
\end{align*}
staying away from the axis in the sense that $\theta(\tau)\in(0,\pi)$ for $\tau\in(0,\overline{L})$, and inducing the Bartnik data 
\begin{align}\label{e:metricWP}
\overline{g}=d\tau^{2}+\overline{\lambda}^{2}d\phi^{2}&=d\tau^{2}+\lambda^{2}d\phi^{2},\\\label{e:HWP}
\overline{H}&=H
\end{align}
on $\partial\Omega$, where $\lambda$ and $H$ are given by
\begin{align}\label{e:overlinelambdaWP}
\lambda&\definedas e^{-U}r\sin\theta,\\
\begin{split}\label{e:overlineHWP}
H&\definedas e^{2(V-U)} 
    \ell^{-3} (-r r'' \theta' + 2 r'^2 \theta' + r r' \theta'' + r^2 \theta'^3)\\
    & + \ell^{-1} \Big[ - \frac{r'}{r} \cot\theta + \theta' 
    + 2\left(\frac{r'}{r} U_{,\theta} - r \theta' U_{,r}\right)\\
    & \qquad + \left(r\sin\theta\, U_{,r}^2 
      - \frac{\sin\theta}{r} U_{,\theta}^2\right)
      (r'\cos\theta + r \theta'\sin\theta)\\
    & \qquad - 2 \sin\theta\, U_{,r} U_{,\theta} 
    (r'\sin\theta - r \theta'\cos\theta)  \Big].
    \end{split}
\end{align}
We will call $([0,\overline{L}],\overline{\lambda},\overline{H})$ \emph{(Weyl--Papapetrou) Bartnik data} and $(\Omega,U,V)$ \emph{(Weyl--Papapetrou) static metric extensions} for simplicity.
\end{conjecture}

Condition~\eqref{e:metricWP} combined with~\eqref{e:overlinelambdaWP}
can be considered as giving Dirichlet boundary values for $U$ along
$\Gamma$, closing the first Weyl--Papapetrou equation~\eqref{e:Ueqn}. 
  Once the unique solution for $U$ is found, the second Weyl--Papapetrou 
  equations~\eqref{e:Veqns} together with the asymptotic condition 
  $V\to 0$ as $r\to\infty$ uniquely determine $V$.

Note that a resolution of this conjecture will indeed resolve Conjecture~\ref{conj:Bartnikaxisymm} as discussed above. However, Conjecture~\ref{conj:Bartnikaxisymm2} is slightly stronger than Conjecture~\ref{conj:Bartnikaxisymm} as we have assumed existence of global Weyl--Papapetrou coordinates which furthermore need to be compatible with the asymptotic flatness assumptions in the derivation of Conjecture~\ref{conj:Bartnikaxisymm2}.

\begin{convention}[Reflection symmetry]\label{conv:reflection}
From now on, we will assume in addition that the Bartnik data $(\Sigma,\overline{g},\overline{H})$ and static metric extensions we consider are \emph{(compatibly) reflection symmetric} in the following sense: once cast in Weyl--Papapetrou coordinates, we request that the free boundary curve $\Gamma$ and thus the domain $\Omega$ satisfy
\begin{align}\label{e:reflection}
\Gamma(\overline{L}-\tau)&=\Gamma(\tau)
\end{align}
for all $\tau\in[0,\overline{L}]$. In terms of the polar coordinates $(r,\theta)$ along $\Gamma$, this~reads
\begin{align}
\begin{split}\label{e:reflpolars}
r(\overline{L}-\tau)&=r(\tau),\\
\theta(\overline{L}-\tau)&=\pi-\theta(\tau)
\end{split}
\end{align}
for $\tau\in[0,\overline{L}]$. We will refer to condition~\eqref{e:reflection} or~\eqref{e:reflpolars} as \emph{reflection symmetry of $\Gamma$}. This condition is compatible with all other conditions stated in Conjecture~\ref{conj:Bartnikaxisymm2}.
\end{convention}

Assuming reflection symmetry will allow us to exploit parity 
and restrict to one half of the computational domain
in the numerical scheme we will describe in Section
\ref{s:nummethod}. Furthermore, assuming reflection symmetry ensures
that the geometric flow we will devise in the next section cannot
``slide up'' the axis of symmetry even though the flow equation~\eqref{e:flow} 
 is invariant under   translations of the curve along the axis. There would of course be other solutions to this sliding issue such as fixing the center of mass, but we prefer to work with symmetry, here. Note that this condition in fact implies the $\mathbb{Z}_{2}\times\mathbb{Z}_{2}\times\mathbb{Z}_{2}$-symmetry assumed in~\cite{Miao2003}.


\subsection{A geometric flow}
\label{s:formulation.flow}
Before we write down the geometric flow we suggest for studying Conjecture~\ref{conj:Bartnikaxisymm2} under the additional assumption of reflection symmetry explained in Convention~\ref{conv:reflection}, let us first introduce some helpful notation. First of all, we will now switch to abstract index notation and write the boundary curve $\Gamma\colon[0,\overline{L}]\to\mathbb{R}^{+}\times[0,\pi]$ as $\Gamma\asdefined x^a$, $a=1,2$. We do \emph{not} assume here that $\Gamma$ is parametrised by arclength, but will stick with $[0,\overline{L}]$ for the domain of definition of $\Gamma$ nevertheless. The unit tangent $t^{a}$ and outward unit normal $n^{a}$ (pointing to the asymptotic end) to the curve $\Gamma$ in the domain $\Omega$ bounded by $\partial\Omega=\Gamma([0,\overline{L}])$ with geometry induced by $(U,V)$ as described above can then be computed to be
\begin{align}\label{e:tangent}
  t^a &= \ell^{-1} (r', \theta'),\\\label{e:normal}
  n^a &= \ell^{-1} \left(r\theta', -\frac{r'}{r}\right),
\end{align}
recalling the definition of $\ell$ given in~\eqref{e:embedding}. If $\Gamma$ is parametrised by arclength, of course $\ell\equiv1$ and~\eqref{e:tangent} and~\eqref{e:normal} simplify accordingly.

Consider now a one-parameter family of curves $\Gamma_t(\tau)= (r_t(\tau), \theta_t(\tau))$ with $\tau \in [0, \overline{L}]$, \emph{not} necessarily parametrised by arclength, with flow ``time'' parameter $t\in[0,T)$ for some $T>0$, with $T=\infty$ allowed in principle. Using abstract index notation, the curves $\Gamma_t \asdefined x^a_t$ will be evolved by the novel geometric curve flow
\begin{align} \label{e:flow}
  \frac{d x_t^a}{dt} &= -(H_t - \overline{H}) n_t^a + C_t t_t^a 
    +\kappa \pi \left(\frac{1}{L_t} - \frac{1}{\overline{L}}\right) n_t^a,
\end{align}
where $t^a_t$ now denotes the unit tangent and $n^a_t$ the outward unit
normal to $\Gamma_t$ with respect to the geometry induced on 
the domain $\Omega_{t}$ and its boundary $\partial\Omega_{t}=\Gamma_{t}([0,\overline{L}])$)
by $(U,V)$, see~\eqref{e:tangent} and~\eqref{e:normal}, $\kappa$ is a coupling parameter, $L_t$ denotes the \emph{length of the curve $\Gamma_{t}$} given by
\begin{align}
  \label{e:length}
  L_t \definedas  \int_0^{\overline{L}} \ell_t \, d\tau,
\end{align}
with $\ell_{t}$ as defined in~\eqref{e:embedding}, and $C_{t}$ is defined as in~\eqref{e:C}, where both $\ell_{t}$ and $C_{t}$ are now computed for $\Gamma_{t}$.

The first term in~\eqref{e:flow} moves the curve in the normal direction by its mean curvature or rather by the difference between its actual mean curvature and the desired mean curvature from the Bartnik data. The second term moves the points of the curve tangentially along the curve in order to
drive the parametrisation of the curve to the desired one in terms of 
 arclength, which corresponds to the isometric embedding of the Bartnik data. The
precise choice of this term is made such that 
\eqref{e:flow} has parabolic character, see Section~\ref{s:analysis.existence}.
The third term counteracts the tendency of the mean curvature flow to shrink 
any curve to a point, and instead drives the curve length $L_{t}$ to the target value
$\overline{L}$. Our analysis in Sections~\ref{s:analysis.spherical} and~\ref{s:analysis.linstab} will reveal that we need the coupling parameter $\kappa > 2$. 

At each instant of $t$, we evaluate $U$ on $\Gamma_t$ using~\eqref{e:Ubc}
with the prescribed embedding term function $\overline{\lambda}(\tau)$.
With this Dirichlet boundary condition for $U$ on $\Gamma_t$ and asymptotic
condition $U\rightarrow 0$ as $r\rightarrow \infty$, we solve
the Laplace equation~\eqref{e:Ueqn} for $U$ in the exterior and then 
determine $V$ by integrating~\eqref{e:Veqns} with asymptotic condition 
$V \rightarrow 0$ as $r\rightarrow \infty$.
Now that we know $U$ and $V$ in the exterior, we can evaluate all the terms
on the right-hand side of~\eqref{e:flow} on $\Gamma_t$, in particular the 
normal derivatives of $U$. We will give arguments in favour of
parabolicity of the symbol of~\eqref{e:flow} and of short-time
existence of solutions to~\eqref{e:flow} in Section~\ref{s:analysis.existence}.

As required in~\eqref{e:smooth} for a single curve, we must impose the boundary conditions
\begin{align}
\begin{split}
\theta_{t}(0)&=0,\quad \theta_{t}(\overline{L})=\pi,\\
r'_{t}(0)&=0,\quad r'_{t}(\overline{L})=0.
\end{split}
\end{align}
for all times $t$ along the geometric flow~\eqref{e:flow}. Imposing these boundary conditions at the initial time $t=0$ gives rise
to compatibility conditions, namely
\begin{align}
  \frac{dr_t'}{dt} \doteq 0, \quad \frac{d\theta_t}{dt} \doteq 0,
\end{align}
at $t=0$, where $\doteq$ denotes equality at $\tau = 0$ and $\tau = \overline{L}$.
We have shown (using computer algebra) that these are satisfied provided that
\begin{align}
  U_{,\theta} \doteq V_{,\theta} \doteq \overline{H}' \doteq 
  r' \doteq \left( \frac{r'}{\sin\theta} \right)' \doteq r''' \doteq
  \theta'' \doteq 0.
\end{align}
These conditions follow from elementary flatness on the axis and are automatically enforced
by the expansions~\eqref{e:spectralexpansion} we use in the code.

Clearly, if $(\Omega,U,V)$ corresponds to a (Weyl--Papapetrou) static metric extension of given (Weyl--Papapetrou) Bartnik data $([0,\overline{L}],\overline{\lambda},\overline{H})$ parametrised by arclength, the geometric flow~\eqref{e:flow} will be \emph{stationary}, i.e.~$\frac{dx^{a}_{t}}{dt}\equiv0$ for all $t\geq0$. This observation is independent of the numerical value of $\kappa$. However, as we will discuss in Section~\ref{s:stationary}, free boundary positions with correctly induced Bartnik data are \emph{not} the only stationary states of the flow, even in a fixed background.

%% file: paper_sec3.tex
\section{Analysis of the geometric flow}
\label{s:analysis}

In order to gain more insight into the novel geometric
flow~\eqref{e:flow} we couple to the Weyl--Papapetrou
equations~\eqref{e:Ueqn} and~\eqref{e:Veqns} in our numerical
analysis, we will now study its properties in some restricted
scenarios such as in spherical symmetry and/or in a fixed
background. First, in Section~\ref{s:analysis.existence}, 
we analyse the symbol of the geometric flow equation~\eqref{e:flow} 
and find that it is parabolic, which suggests short-time existence of solutions. Next, in Section~\ref{s:analysis.spherical}, we study the geometric flow~\eqref{e:flow}, and briefly the coupled system, in spherical symmetry (Euclidean and Schwarzschildean backgrounds). In particular, we will discuss the chosen threshold for the coupling parameter~$\kappa$ there. In Section~\ref{s:analysis.linstab}, we will linearise the geometric flow in a fixed Euclidean background around a coordinate circle and study the behaviour of the linearised flow. Finally, in Section~\ref{s:stationary}, we will briefly discuss the occurrence of rather unintended stationary states of the geometric flow~\eqref{e:flow} and the coupled flow-Weyl--Papapetrou system.

It would of course be desirable both to study the flow in other fixed backgrounds and to analyse the full coupled system. We will not pursue these ideas here as they would lead too far for this first treatment of Bartnik's conjecture in axisymmetry.


\subsection{Short-time existence}
\label{s:analysis.existence}

Equation~\eqref{e:flow} is a system of parabolic partial differential
equations in $t$ and $\tau$ for $x^a_t(\tau)$.
More precisely,
\[ \partial_t \left( \begin{array}{c}r_t(\tau) \\\theta_t(\tau) 
                     \end{array} \right)
  = \ell^{-2} \left( \begin{array}{c}r_t''(\tau) \\\theta_t''(\tau) \end{array}
  \right) + \ldots,
\]
where the omitted terms do not contain any second derivatives;
note also from \eqref{e:embedding} that $\ell$ contains up to first 
$\tau$-derivatives of $r_t(\tau)$ and $\theta_t(\tau)$.
Thus we have a manifestly strongly parabolic quasi-linear system of
second order. The principal part decouples into two separate heat equations
  thanks to a cancellation between the first two terms in
  \eqref{e:flow}; to make this work, we introduced a factor
  $\ell^{-2}$ in the definition of $C$ in \eqref{e:C}.

If we disregard the last term in \eqref{e:flow}, which because of
\eqref{e:length} contains a $\tau$-integral of the
unkowns $x_t^a(\tau)$, then standard theorems (e.g.~Theorem 7.2 in
\cite{TaylorBook}) imply that for smooth initial data, a unique smooth
solution exists for a finite time.

Given that an integral term of a very similar form occurs in
area preserving curve shortening flow~\cite{Gage1986} and other
constrained curve flows~\cite{Dittberner2017}, we expect that
short-time existence results for such flows~\cite{Huisken1987,Pihan1998}
will carry over to our flow as well, even though we have not proven
such a theorem (which is complicated by the fact that the functions 
$U, V$ occurring in the flow equation~\eqref{e:flow} are determined by
elliptic equations in the coupled case).

A further indication is obtained in Section
  \ref{s:analysis.linstab}, where we show that when the flow equation
  is linearised about a circle in a Euclidean background, then a unique
  solution to the initial-boundary value problem exists for all future times.

We cannot make any statements about global existence
of the nonlinear coupled flow, although our
numerical experiments indicate that in a variety of situations bounded
smooth solutions exist for infinite flow time.


\subsection{Flowing in spherical symmetry}
\label{s:analysis.spherical}
In this section, we will restrict our attention to the spherically symmetric case, i.e.~to the case where the axisymmetric static metric extensions characterised by $(U,V)$ are indeed spherically symmetric and the free boundary curves evolving under the geometric flow~\eqref{e:flow} are circles in coordinates adjusted to the spherical symmetry and thus represent symmetry reductions of the orbital spheres of the spherically symmetric static metric extensions. This of course corresponds to prescribing Bartnik data that are spherically symmetric or in other words that have a round metric and constant positive mean curvature, see also the discussion on page~\pageref{p:sphericalsymmetry}.

We will first study evolving circles in a Euclidean background and give a very brief insight into the coupled system with prescribed Bartnik data corresponding to a centred coordinate sphere in a Euclidean background in Section~\ref{s:Euclidean}. Then, we will study the evolution of circles corresponding to centred round spheres in a fixed Schwarzschild background in Section~\ref{s:Schwarzschild} and investigate the stability of the flow in a fixed Schwarzschild background.

\subsubsection{Euclidean case}\label{s:Euclidean}
Let us first look at the geometric flow~\eqref{e:flow} in a fixed Euclidean background, corresponding to $U\equiv V\equiv0$. We will study the evolution of coordinate circles in Weyl--Papapetrou coordinates, which indeed coincide with standard (Euclidean) polar coordinates and represent symmetry reductions of Euclidean coordinate spheres centred at the origin. We will furthermore assume that these circles are parametrised proportionally to arclength. Our Bartnik data are then uniquely determined by prescribing their coordinate sphere radius $\overline{R}$, and we have $\overline{L}=\pi\overline{R}$.
 
In Weyl--Papapetrou coordinates, circles evolving under~\eqref{e:flow} are described by
\begin{align}
  r_t(\tau) \equiv R_t,\quad  \theta_t(\tau) = \frac{\pi \tau}{\overline{L}},
\end{align}
and we obtain by direct computation from the formulas in Section~\ref{s:formulation} (or, equivalently, by directly computing all the geometric notions for coordinate spheres in Euclidean space) that
\begin{align}
  H_t = \frac{2}{R_t},\quad   \overline{H}_t \equiv\frac{2}{\overline{R}},\quad
  n_t^a = (1, 0),\quad  C_t = 0,
\end{align}
so that the geometric flow~\eqref{e:flow} reduces to the ODE
\begin{align}\label{e:ODEEuclidean}
  \frac{dR_t}{dt} = (\kappa-2) \left( \frac{1}{R_t} - \frac{1}{\overline{R}} 
  \right).
\end{align}
Provided $\kappa > 2$, the circle of radius $R_{t}=\overline{R}$ will be a stationary state of~\eqref{e:ODEEuclidean}, while circles with $R_t > \overline{R}$ will shrink and circles with $R_t < \overline{R}$ will expand with $R_t \rightarrow \overline{R}$
as $t \rightarrow \infty$, as desired, showing global stability of the (unique) stationary state $R_{t}\equiv\overline{R}$ (but see Section~\ref{s:stationary}). For $\kappa\leq2$, the flow does not produce the desired behaviour. This observation gives rise to the restriction $\kappa>2$ we introduced in Section~\ref{s:formulation}.

Now, still looking at Weyl--Papapetrou coordinate circles parametrised proportionally to arclength, let us look at the coupled flow--Weyl--Papapetrou system, meaning that we continuously solve for $U$ und $V$ outside the flowing circles. Again, we consider Bartnik data corresponding to a \emph{target} circle of radius $\overline{R}$. The corresponding \emph{target functions} are
\begin{align}
  \overline{\lambda}(\tau) &= \frac{\overline{L}}{\pi} \sin \left( \frac{\tau \pi}{\overline{L}} 
  \right),\\
  \overline{H}(\tau) &\equiv \frac{2\pi}{\overline{L}}.
\end{align}
We start the coupled flow--Weyl--Papapetrou system with $U_{0}\equiv V_{0}\equiv0$, i.e.~with fields corresponding to Euclidean space. During the flow, we assume the flowing curve $\Gamma_{t}=(r_{t},\theta_{t})$ remains a circle with a $t$-dependent
radius $R_{t}$, so that
\begin{align}
  r_t(\tau) \equiv \frac{L_t}{\pi},\quad  \theta_t(\tau) = \frac{\pi\tau}{\overline{L}}.
\end{align}
Hence, the free boundary condition for $U=U_{t}$,~\eqref{e:Ubc}, evaluates to 
\begin{equation}
  U_t \circ \Gamma_t (\tau) =-\ln \frac{\overline{\lambda}(\tau)}
  {r_t(\tau) \sin\theta_t(\tau)} = - \ln \frac{\overline{L}}{L_t}
\end{equation}
along the evolving free boundary curve, i.e.~$U_{t}$ is constant along $\Gamma_{t}$, although changing in time, and, in fact \emph{differs} from the Euclidean $U\equiv0$. Thus, inserting this Dirichlet boundary condition for $U_{t}$ into the Weyl--Papapetrou equation~\eqref{e:Ueqn} for $U=U_{t}$ and imposing the asymptotic flatness condition $U_{t}\rightarrow 0$ as $r\rightarrow \infty$, we will thus get a spherically symmetric solution $U_{t}$ which is non-zero for finite time $t$, and thus corresponds to a non-Euclidean metric (once we have also solved for $V_{t}$). If $L_{t}\to\overline{L}$ as $t\to T$ for some $T\leq\infty$, approaching a stationary state as one may expect, one sees that indeed $U_{t}\to 0$ and indeed also $V_{t}\to0$ as $t\to T$ (see, however, Section~\ref{s:stationary}). We indeed observe this (temporary) deviation from Euclidean space for Euclidean Bartnik data and $U_{0}\equiv V_{0}\equiv0$ numerically in much more general situations, see Section~\ref{s:numresults}.

\subsubsection{Schwarzschild case}\label{s:Schwarzschild}
Let us now look at the geometric flow~\eqref{e:flow} in a fixed Schwarzschild background of mass $M>0$, cf.~\eqref{e:schwarzschild}. To do so, we will need to express the Schwarzschild background in terms of potentials $(U,V)$ in Weyl--Papapetrou coordinates, see~\cite{Griffiths2009}. In cylindrical Weyl--Papapetrou coordinates 
$\rho=r\sin\theta$ and $z=r\cos\theta$, the fields $U$ and $V$ are given by
\begin{align}
 \begin{split} \label{e:SchwarzschildUV}
  U &= \frac{1}{2}  \ln \frac{R_+ + R_- - 2M}
    {R_+ + R_- + 2M},\\[1ex]
  V &= \frac{1}{2}  \ln \frac{(R_+ + R_-)^2 - 4 M^2}{4 R_+ R_-},
\end{split}\end{align}
where
\begin{align}
  R_\pm &= \sqrt{\rho^2 + (z \pm M)^2}.
\end{align}
In cylindrical Weyl--Papapetrou coordinates, the Schwarzschild metric becomes singular on the piece of the axis of rotation given by $\lbrace{\rho=0, |z|<M\rbrace}$ which corresponds to the black hole horizon. The usual Schwarzschild coordinates $(r_S, \theta_S,\varphi_S=\varphi)$ (used in~\eqref{e:schwarzschild} without the index $S$) are related to the Weyl--Papapetrou coordinates via the coordinate transformations
\begin{align}
\begin{split}\label{e:transfo}
\rho &\asdefined M \sqrt{(x^2-1)(1-y^2)}, \\
z &\asdefined M x y,\\[1ex]
  r_S &\definedas M(x+1),\\
   \cos\theta_S &\definedas y.
  \end{split}
\end{align}

In order to gain more insight into the geometric flow~\eqref{e:flow}, now in a fixed Schwarzschild background, i.e.~with $(U,V)$ as in~\eqref{e:SchwarzschildUV}, let us consider flowing Schwarzschild coordinate circles parametrised proportionally to arclength, corresponding to centred Schwarzschild coordinate spheres from a $3$-dimensional viewpoint. Note that it is not immediately obvious but indeed follows from the geometric nature of~\eqref{e:flow} that such Schwarzschild coordinate circle solutions to~\eqref{e:flow} exist; in contrast to the Euclidean case discussed in Section~\ref{s:Euclidean}, they cannot be written as coordinate circles in Weyl--Papapetrou coordinates. 

In order to derive the ODE for the circle radius, we perform the following computations. Let us first look at a single circle of Schwarzschild radius $r_S=R$ in our fixed Schwarzschild background of mass $M$. Written as a curve in Schwarzschild coordinates $(r_S,\theta_S)$ which is parametrised proportionally to arclength on some interval $\left[0,\overline{L}\right]$, we find $r_S(\tau)\equiv R$ and $\theta_S(\tau)=\frac{\tau}{\overline{L}}$ as in the Euclidean case. Also, abbreviating $\overline{R}\definedas\frac{\overline{L}}{\pi}$, one computes 
\begin{align}\label{e:L}
L&=\pi R,\\\label{e:l}
 \ell(\tau)&\equiv\frac{R}{\overline{R}},\\\label{e:H(R)}
 H(\tau)&=\frac{2}{R}\,\sqrt{1-\frac{2M}{R}}.
\end{align}
Performing the canonical transformation into the cylindrical Weyl--Papapetrou coordinates~\eqref{e:transfo} and then changing back into standard polar Weyl--Papapetrou coordinates in which~\eqref{e:flow} is written, we obtain
\begin{align}\label{e:r(tau)}
r(\tau)&=M\,\sqrt{\left(\frac{R}{M}-1\right)^2-\sin^2\left(\frac{\tau}{\overline{R}}\right)},\\\label{e:theta(tau)}
\theta(\tau)&=\arctan\left(\frac{\sqrt{1-\frac{2M}{R}}}{1-\frac{M}{R}}\,\tan\left(\frac{\tau}{\overline{R}}\right)\right).
\end{align}
Now consider a flow of circles solving~\eqref{e:flow} with radius $(r_{S})_{t}(\tau)\equiv R_t$. The geometric flow equation~\eqref{e:flow} then reduces to the system
\begin{align}\label{e:rcircle}
  \frac{d r_t}{dt} &= \left(-(H_t - \overline{H}) + \kappa \pi \left(\frac{1}{L_t} - \frac{1}{\overline{L}}\right) \right)\ell_t^{-1}\,r_t\,\theta'_t,\\\label{e:thetacircle}
  \frac{d \theta_t}{dt} &= \left(-(H_t - \overline{H}) + \kappa \pi \left(\frac{1}{L_t} - \frac{1}{\overline{L}}\right) \right)\ell_t^{-1}\,\left(-\frac{r'_t}{r_t}\right).
\end{align}
Recall that $'$ denotes a derivative with respect to $\tau$. Here, $\overline{H}=\frac{2}{\overline{R}}\,\sqrt{1-\frac{2M}{\overline{R}}}$ and of course $\overline{L}=\pi\overline{R}$. Evaluating~\eqref{e:L}--\eqref{e:theta(tau)} along the flow and plugging them
 into~\eqref{e:rcircle} and~\eqref{e:thetacircle}, both of these flow equations reduce to the following ODE for the radius $R_t$ of the flowing circles
 \begin{align}
 \begin{split}\label{e:SchwarzODE}
 \frac{dR_t}{dt}&=\left[-\left(\frac{2}{R_t}\sqrt{1-\frac{2M}{R_t}}-\frac{2}{\overline{R}}\sqrt{1-\frac{2M}{\overline{R}}}\right)+\kappa\left(\frac{1}{R_t}-\frac{1}{\overline{R}}\right)\right]\\
 &\qquad\times\sqrt{1-\frac{2M}{R_t}}.
\end{split}
\end{align}

This ODE is consistent with the Euclidean case discussed above, where \eqref{e:ODEEuclidean} arises from~\eqref{e:SchwarzODE} by setting $M=0$, as expected. Clearly, this ODE has a stationary solution $R_t\equiv \overline{R}$ such that the prescribed Bartnik data or in other words the target circle $r_{S}=\overline{R}$ is indeed a stationary state of the flow. We will now perform a direct computation to show that there are no other stationary states of this ODE for $\kappa>2$, recalling that the Schwarzschild radial coordinate needs to remain larger than the black hole radius, $r_{S}>2M$. The same computation will demonstrate that circles with $R_t > \overline{R}$ will shrink and circles with $2M<R_t < \overline{R}$ will expand with $R_t \rightarrow \overline{R}$ as $t \rightarrow \infty$, as in the Euclidean case and as desired, provided that $\kappa>2$. In particular, this computation will show that the stationary state $R_{t}\equiv\overline{R}$ is globally stable as a stationary point of~\eqref{e:SchwarzODE} as in the Euclidean case discussed in Section~\ref{s:Euclidean} (but see Section~\ref{s:stationary}). To prove these claims, let us set $\rho\definedas\frac{R_{t}}{\overline{R}}$ and rewrite~\eqref{e:SchwarzODE} as
\begin{align*}
\frac{R_{t}}{\sqrt{1-\frac{2M}{R_{t}}}\left(\kappa-2\sqrt{1-\frac{2M}{\overline{R}}}\right)}\times\frac{dR_t}{dt}&=\frac{\kappa-2\sqrt{1-\frac{2M}{\rho\overline{R}}}}{\kappa-2\sqrt{1-\frac{2M}{\overline{R}}}}-\rho\asdefined f_{M}^{\kappa}(\rho).
\end{align*}
This gives $f_{M}^{\kappa}(1)=0$ and $\frac{df_{M}^{\kappa}(\rho)}{d\rho}<0$ as long as $\rho>\frac{2M}{\overline{R}}$, which asserts both that the only zero of $f_{M}^{\kappa}$ is $\rho=1$, proving that $R_{t}\equiv\overline{R}$ is the only stationary state of~\eqref{e:SchwarzODE}, and that the sign of  $\frac{dR_t}{dt}$ is negative for $R_{t}>\overline{R}$ and positive for $R_{t}<\overline{R}$, as needed to assert global stability of the stationary state $R_{t}\equiv\overline{R}$ as a solution of~\eqref{e:SchwarzODE} (but see Section~\ref{s:stationary}).

\subsection{Linear stability analysis in a Euclidean background}
\label{s:analysis.linstab}

In the above considerations, we have analysed the stability of the stationary state $R_{t}\equiv\overline{R}$ in the spherically symmetric ODE setting and found global stability both in a fixed Euclidean and a fixed Schwarzschild background. Complementing this analysis, we will now linearise the geometric flow~\eqref{e:flow} around such a stationary circle, but only in a Euclidean background, leaving the computationally more involved Schwarzschild case for future work. As the background is Euclidean, we have $U\equiv V\equiv0$ and the stationary circle in Weyl--Papapetrou coordinates, parametrised by arc\-length, is given by $r(\tau)\equiv\overline{R}$, $\theta(\tau)=\frac{\tau}{\overline{R}}$,  as above. Recall this circle corresponds to Bartnik data consisting of a round sphere of radius $\overline{R}$ in Euclidean space centred at the origin. We recall that
\begin{align}
\begin{split}\label{e:stationarycircle}
 \overline{L}=\pi\overline{R},&\quad 
 \overline{H}\equiv\frac{2}{\overline{R}},\\
 \overline{C}\equiv0,&\quad \overline{\ell}\equiv1.
 \end{split}
 \end{align}
To compute the linearised flow equations, set
\begin{align}\label{def:phi}
r_t(\tau)&\asdefined \overline{R}+\varepsilon\mu_t(\tau),\\\label{def:psi}
\theta_t(\tau)&\asdefined \frac{\tau}{\overline{R}}+\varepsilon\psi_t(\tau),
\end{align}
for smooth families of functions $\mu_{t},\psi_{t}\colon[0,\overline{L}]\to\mathbb{R}$, $t\in[0,T)$, and $\varepsilon$ a small parameter. In order to sustain the boundary conditions~\eqref{e:smooth}, we need to ask that
\begin{align}\label{e:bdryphi}
 \mu'_t(0)&=\mu'_t(\overline{L})=0,\\\label{e:bdrypsi}
 \psi_t(0)&=\psi_t(\overline{L})=0
\end{align}
for all $t\in[0,T)$. Moreover, recalling our Convention~\ref{conv:reflection} on reflection symmetry, we need to ask in addition that
\begin{align}\label{e:musymm}
\mu_{t}(\overline{L}-\tau)&=\mu_{t}(\tau),\\\label{e:psisymm}
\psi_{t}(\overline{L}-\tau)&=-\psi_{t}(\tau)
\end{align}
for $\tau\in[0,\overline{L}]$ and all $t\in[0,T)$. As a consequence of~\eqref{e:C},~\eqref{e:H} as well as~\eqref{e:stationarycircle} and~\eqref{e:bdrypsi}, we find the following linearised system
 \begin{align}\label{e:phi_t}
    \frac{d\mu_t}{dt}&=\mu''_t+\frac{\mu'_t}{\overline{R}}\cot\left(\frac{\tau}{\overline{R}}\right)+\frac{2\mu_t}{\overline{R}^2}-\frac{\kappa}{\pi\overline{R}^3}\int_0^{\overline{L}}\mu_t(\tau)\,d\tau,\\\label{e:psi_t}
    \frac{d\psi_t}{dt}&=\frac{\mu'_t}{\overline{R}^2}+\psi''_t.
\end{align}

By linearisation, the system~\eqref{e:phi_t}, \eqref{e:psi_t} is a system of linear parabolic equations of second order. Observe that~\eqref{e:phi_t} decouples from~\eqref{e:psi_t} as it does
not contain $\psi_{t}$. Equation~\eqref{e:psi_t} is an inhomogeneous linear 
heat equation for $\psi_t$, once $\mu_t$ has been computed. 
We will see below that, for given initial data, the
system~\eqref{e:phi_t}, \eqref{e:psi_t} has a unique solution for all
(future) times.
To solve~\eqref{e:phi_t}, we make an ansatz of separation of variables
\begin{align}
 \mu_t(\tau)=a(t)\,b(\tau).
\end{align}
Dividing as usual by $ab$, we rearrange~\eqref{e:phi_t} to
\begin{align}\label{e:sep}
 \frac{\dot{a}}{a}&=\frac{b''}{b}+\frac{b'}{b\,\overline{R}}\cot\left(\frac{\tau}{\overline{R}}\right)+\frac{2}{\overline{R}^2}-\frac{\kappa}{b\,\pi\overline{R}^3}\int_0^{\overline{L}}b(\tau)\,d\tau\asdefined \frac{\alpha+2}{\overline{R}^2},
\end{align}
where $\dot{a}=\frac{da}{dt}$ and $b'=\frac{db}{d\tau}$ as before, and $\alpha$ denotes a real parameter. Equation~\eqref{e:sep} can then immediately be seen to possess the unique solution
\begin{align}
 a(t)&=a_\alpha\,\exp\left(\frac{\alpha+2}{\overline{R}^2}\,t\right)
\end{align}
for some $a_\alpha\in\mathbb{R}$ for $a$. In order to show linear stability of the stationary state $r\equiv\overline{R}$ of the geometric flow~\eqref{e:flow}, we need to show that only modes with $\alpha+2<0$ prevail, or in other words that $\alpha\geq-2$ implies $b\equiv0$. 

So let us study solutions to~\eqref{e:sep} for $b$, 
\begin{align}\label{e:b}
 b''+\frac{b'}{\overline{R}}\cot\left(\frac{\tau}{\overline{R}}\right)-\frac{\kappa}{\pi\overline{R}^3}\int_0^{\overline{L}}b(\tau)\,d\tau&=\frac{\alpha}{\overline{R}^2} \,b.
\end{align}
From the boundary conditions~\eqref{e:bdryphi}, we find the boundary conditions
\begin{align}\label{e:bbounds}
 b'(0)=b'(\overline{L})=0
\end{align}
for~\eqref{e:b}. To simplify this equation consistently with its boundary conditions, we perform the following transformation of variables
\begin{align}
\begin{split}\label{e:trafo}
x&\definedas \cos\left(\frac{\tau}{\overline{R}}\right),\\
\widehat{b}(x)&\definedas b(\tau).
\end{split}
\end{align}
This transformation allows us to rewrite~\eqref{e:b} as
\begin{align}\label{e:bwidehat}
(1-x^{2})\,\widehat{b}\,''(x)-2x\,\widehat{b}\,'(x)-\alpha\,\widehat{b}(x)&=\frac{\kappa}{\pi}\int_{-1}^{1}\frac{\widehat{b}(y)}{\sqrt{1-y^{2}}}\,dy
\end{align}
for $x\in[-1,1]$, where, slightly abusing notation, we denote $x$-derivatives by $'$ as well. The transformation~\eqref{e:trafo} is only allowed because of the boundary conditions~\eqref{e:bbounds}, it would otherwise be degenerate at the endpoints of $[-1,1]$. There are no prescribed boundary values for $\widehat{b}$. By reflection symmetry~\eqref{e:musymm}, 
\begin{align}\label{e:widehatbeven}
\widehat{b}(-x)=\widehat{b}(x)\text{ for all }x\in[-1,1].
\end{align}

To analyse solutions to~\eqref{e:bwidehat}, let us first treat the
case that $\widehat{b}$ is constant, 
$\widehat{b}(x)\equiv:\lambda$, 
which takes a special role because of the integral term on the right hand side of \eqref{e:bwidehat}. Plugging this ansatz into \eqref{e:bwidehat}, we obtain 
\begin{align}
-\alpha\lambda&=\frac{\kappa\lambda}{\pi}\int_{-1}^{1}\frac{1}{\sqrt{1-y^{2}}}\,dy=\kappa\lambda
\end{align}
so that $\alpha=-\kappa$, or $\lambda=0$ and hence $\widehat{b}\equiv0$. Thus, the only spatially constant mode of $\mu_{t}$ decays in time exponentially with rate $\frac{2-\kappa}{\overline{R}^{2}}$.

Now that we have understood the constant case, let us analyse the full set of solutions of \eqref{e:bwidehat} by distinguishing the cases $\alpha=0$ and $\alpha\neq0$. If $\alpha=0$, the $x$-derivative of~\eqref{e:bwidehat} gives
\begin{align}
\left((1-x^{2})\,\widehat{b}\,'(x)\right)''&=0,
\end{align}
so that there exist constants $\lambda,\rho\in\mathbb{R}$ such that
\begin{align}\label{e:integrable}
(1-x^{2})\,\widehat{b}\,'(x)&=\lambda x+\rho\text{ on }[-1,1].
\end{align}
Plugging $x=\pm1$ into~\eqref{e:integrable} leads to $\lambda=\rho=0$
or in other words $\widehat{b}$ must be constant on $[-1,1]$, which, by the above discussion, immediately implies $\widehat{b}\equiv0$, recalling that $\kappa>2$ (cf.~Section~\ref{s:Euclidean}). 
Thus $\alpha=0$ is excluded or in other words there is no mode in $\mu_{t}$ which grows exponentially with rate $\frac{2}{\overline{R}^{2}}$.

Let us now discuss the case $\alpha\neq0$. The right hand side of~\eqref{e:bwidehat} is manifestly constant. This observation can be rephrased as saying that
\begin{align}
(1-x^{2})\,\widehat{b}''(x)-2x\,\widehat{b}'(x)-\alpha\,\widehat{b}(x)&=\lambda
\end{align}
for some $\lambda\in\mathbb{R}$. By setting
\begin{align}\label{e:btwiddle}
\widetilde{b}(x)&\definedas \widehat{b}(x)+\frac{\lambda}{\alpha},
\end{align}
we find the Legendre differential equation
\begin{align}\label{e:Legendre}
(1-x^{2})\,\widetilde{b}''(x)-2x\,\widetilde{b}'(x)-\alpha\,\widetilde{b}(x)&=0
\end{align}
for $\widetilde{b}$ on $[-1,1]$. It is well known that the only solutions $\widetilde{b}$ to~\eqref{e:Legendre} extending continuously to $[-1,1]$ are the Legendre polynomials $P_{l}$, with $l\in\mathbb{N}_{0}$, where $\alpha=-l(l+1)$. We have already handled (excluded) $l=\alpha=0$. For $l=1$, $\alpha=-2$, we find the translational mode $\widetilde{b}(x)=Bx$ with $B\in\mathbb{R}$ and thus $\widehat{b}(x)=Bx-\frac{\lambda}{\alpha}$. On the other hand, we know that
\begin{align}
\lambda&=\frac{\kappa}{\pi}\int_{-1}^{1}\frac{\widehat{b}(y)}{\sqrt{1-y^{2}}}\,dy=\frac{\kappa}{\pi}\left(0-\frac{\lambda}{\alpha}\pi\right),
\end{align}
which implies $\lambda=0$ or $\kappa=-\alpha=2$. As $\kappa>2$, we
deduce $\lambda=0$. Moreover, $B=0$ by~\eqref{e:widehatbeven}. But
this leads to $\widehat{b}\equiv0$ in case $l=1$ or 
equivalently $\alpha=-2$. In other words, there is no constant-in-time mode in $\mu_{t}$. For all $l\geq2$ however, $\alpha=-l(l+1)<-2$, so that all modes discovered via the separation of variables performed above are decaying. Now recall that the Legendre polynomial $P_{l}$ is even if $l\in\mathbb{N}_{0}$ is even, and odd if $l$ is odd, and that the set of even Legendre polynomials $\lbrace{P_{2n}\,\vert\,n\in\mathbb{N}_{0}\rbrace}$ is a complete orthonormal system for the separable Hilbert space $\lbrace{\widehat{u}\in L^{2}[-1,1]\,\vert\,\widehat{u}\text{ even}\rbrace}$. Reversing the transformations~\eqref{e:trafo} and~\eqref{e:btwiddle}, this allows us to conclude that every smooth $u\colon[0,\overline{L}]\to\mathbb{R}$ satisfying the boundary conditions~\eqref{e:bdryphi} and the symmetry condition~\eqref{e:musymm} can be expanded as 
\begin{align}\label{e:bgeneral}
u(\tau)&=B_{0}+\sum_{n=1}^{\infty} B_{n} b_{2n}(\tau)
\end{align} 
for all $\tau\in[0,\overline{L}]$ with suitable coefficients
$B_{0}, B_{n}\in\mathbb{R}$, where
\begin{align}\label{e:bl} 
b_{l}(\tau)&\definedas P_{l}\left(\cos\left(\frac{\tau}{\overline{R}}\right)\right)+\frac{\lambda_{l}}{l(l+1)},\\\label{e:lambdal}
\lambda_{l}&\definedas \frac{l(l+1)\kappa}{(l(l+1)-\kappa)\pi}\int_{-1}^{1}\frac{P_{l}(y)}{\sqrt{1-y^{2}}}\,dy
\end{align}
for $\tau\in[0,\overline{L}]$ and for each $l\in\mathbb{N}_{>0}$. Hence, in view of Parseval's identity for complete orthonormal systems, we can conclude that the separation of variables performed above has found the most general solution of~\eqref{e:phi_t} consistent with the boundary conditions~\eqref{e:bdryphi} and the symmetry condition~\eqref{e:musymm} which consequently must be of the form
\begin{align}\label{e:mugeneral}
\mu_{t}(\tau)&=A_{0}\exp\left(\frac{2-\kappa}{\overline{R}^{2}}t\right)+\sum_{n=1}^{\infty}A_{n}\exp\left(\frac{2-2n(2n+1)}{\overline{R}^{2}}\,t\right)b_{2n}(\tau),
\end{align}
for all $t>0$ and all $\tau\in[0,\overline{L}]$, with suitable coefficients $A_{0}, A_{n}\in\mathbb{R}$. This finishes the linear stability argument for $\mu_{t}$.

Quantitatively speaking, we found that the most slowly decaying mode of $\mu_{t}$ decays at least as fast as $\exp\left(-\frac{4}{\overline{R}^{2}}t\right)$ or as $\exp\left(\frac{2-\kappa}{\overline{R}^{2}}t\right)$, depending on the choice of $\kappa>2$ --- whichever of those modes decays more slowly.

The general form~\eqref{e:mugeneral} of a solution of \eqref{e:phi_t} allows us to establish existence and uniqueness of solutions of~\eqref{e:phi_t} for all times $t>0$ for given initial data $\mu_{0}\colon[0,\overline{L}]\to\mathbb{R}$. Indeed, let $\mu_{0}$ be smooth initial data for~\eqref{e:phi_t}, obeying the boundary conditions
\begin{align}\label{e:bdryphi0}
\mu_{0}'(0)&=\mu_{0}'(\overline{L})=0
\end{align}
as in~\eqref{e:bdryphi} and the symmetry condition
\begin{align}\label{e:mu0symm}
\mu_{0}(\overline{L}-\tau)&=\mu_{0}(\tau)
\end{align}
for all $\tau\in[0,\overline{L}]$ as in~\eqref{e:musymm}. As asserted above (see~\eqref{e:bgeneral} and the text above it), we can expand $\mu_{0}$ uniquely as
\begin{align}
\mu_{0}(\tau)&=A_{0}+\sum_{n=1}^{\infty}A_{n}b_{2n}(\tau)
\end{align}
for $\tau\in[0,\overline{L}]$ with coefficients $A_{n}\in\mathbb{R}$
for $n\in\mathbb{N}_{0}$. Thus, \eqref{e:mugeneral} is the unique
solution of~\eqref{e:phi_t} with initial data $\mu_{0}$ and automatically exists for all future times $t>0$.

Let us now study~\eqref{e:psi_t} to get information on $\psi_{t}$,
exploiting that we have already asserted that $\mu_{t}$ is
decaying and of the form~\eqref{e:mugeneral}. By linearity of~\eqref{e:psi_t}, we can separately
investigate the behaviour of the solutions to the homogeneous 
heat equation
\begin{align}\label{e:homogeneous}
\frac{d\psi_t}{dt}&=\psi''_t
\end{align}
and that of special solutions to the inhomogeneous system~\eqref{e:psi_t}, inserting the individual modes of $\mu_{t}$.

First, considering the homogeneous system~\eqref{e:homogeneous}, using the separation of variables
\begin{align}\label{e:sep2ansatz}
\psi_{t}(\tau)&=c(t)d(\tau),
\end{align}
dividing as usual by $cd$, and rearranging~\eqref{e:homogeneous}, we find
\begin{align}\label{e:sep2}
\frac{\dot{c}}{c}&=\frac{d''}{d}\asdefined -\frac{\beta}{\overline{R}^{2}}
\end{align}
for some real parameter $\beta$. As above, \eqref{e:sep2} can then immediately be seen to possess the unique solution
\begin{align}
 c(t)&=c_\beta\,\exp\left(-\frac{\beta}{\overline{R}^2}\,t\right)
\end{align}
for some $c_\beta\in\mathbb{R}$ for $c$. Thus, in order to show decay of the solutions to the homogeneous part of~\eqref{e:psi_t}, i.e. to~\eqref{e:homogeneous}, we need to ensure that $\beta>0$. 

Let us first study the case $\beta=0$. In this case, \eqref{e:sep2} implies that there are constants $\lambda$, $\rho$ for which $d$ satisfies
\begin{align}
d(\tau)&=\lambda\tau+\rho
\end{align}
for all $\tau\in[0,\overline{L}]$. The boundary conditions~\eqref{e:bdrypsi}, or, alternatively, the reflection symmetry condition~\eqref{e:psisymm}, tell us that $d(0)=d(\overline{L})=0$ so that $\rho=\lambda=0$ and hence $d\equiv0$ in this case. This rules out the case $\beta=0$.  

Next, let us study the case $\beta<0$. In this case, $d$ can be written as 
\begin{align}
d(\tau)&=\lambda\sinh\left(\frac{\sqrt{-\beta}}{\overline{R}}\,\tau\right)+\rho\cosh\left(\frac{\sqrt{-\beta}}{\overline{R}}\,\tau\right)
\end{align}
for some constants $\lambda$, $\rho$. The boundary condition $d(0)=0$ that follows from~\eqref{e:bdrypsi} gives $\rho=0$. The reflection symmetry condition $d(\overline{L}-\tau)=-d(\tau)$ that follows from~\eqref{e:psisymm}, together with the hyperbolic addition theorem, then implies that $\cosh\left(\pi\sqrt{-\beta}\right)=1$ if $\lambda\neq0$, recalling $\overline{L}=\pi\overline{R}$. This, however, is of course excluded for $\beta<0$ so that $\lambda=0$ and thus again $d\equiv0$, which rules out the case $\beta<0$ and indeed asserts decay of all modes of the solution to the homogeneous system~\eqref{e:homogeneous}.

In order to get a better idea of how fast the slowest mode of the solution to the homogeneous system~\eqref{e:homogeneous} is decaying, let us briefly also consider the case $\beta>0$. Arguing as before, we find
\begin{align}\label{e:d}
d(\tau)&=\lambda\sin\left(\frac{\sqrt{\beta}}{\overline{R}}\,\tau\right)+\rho\cos\left(\frac{\sqrt{\beta}}{\overline{R}}\,\tau\right)
\end{align}
for some constants $\lambda$, $\rho$, and the boundary and reflection
symmetry requirements, together with the trigonometric addition
theorem, tell us that $\rho=0$ and that  $\cos\left(\pi\sqrt{\beta}\right)=1$.
This leads to $\beta=4n^{2}$ for some $n\in\mathbb{N}_{>0}$. Hence, the most slowly decaying mode of the solution to the homogeneous equation~\eqref{e:homogeneous} decays at least as fast as $\exp\left(-\frac{4}{\overline{R}^{2}}t\right)$.

It follows from the above considerations that the general solution to the homogeneous equation~\eqref{e:homogeneous}, $\psi_{t}^{\text{hom}}$, reads
\begin{align}\label{e:hompsi}
\psi_{t}^{\text{hom}}(\tau)&=\sum_{n=1}^{\infty}C_{n}\exp\left(-\frac{4n^{2}}{\overline{R}^2}\,t\right)\sin\left(\frac{2n}{\overline{R}}\,\tau\right)
\end{align}
for suitable constants $C_{n}\in\mathbb{R}$.

We now turn our attention to the inhomogeneous system~\eqref{e:psi_t}, inserting a fixed mode solution
\begin{align}
\mu_{t}(\tau)&=A_{n}\exp\left(\frac{2-2n(2n+1)}{\overline{R}^{2}}\,t\right)b_{2n}(\tau)
\end{align}
for some $n\in\mathbb{N}_{>0}$, where $b_{2n}$ is given by~\eqref{e:bl}, and $A_{n}\in\mathbb{R}$. This leads to the equation
\begin{align}
\frac{d\psi_t}{dt}&=\frac{A_{n}}{\overline{R}^2}\exp\left(\frac{2-2n(2n+1)}{\overline{R}^{2}}\,t\right)\,b_{2n}'+\psi''_t.
\end{align}
Making the same 
ansatz~\eqref{e:sep2ansatz} of separation of variables as before (but multiplied by $A_{n}$), we find
\begin{align}\label{e:inhomoseparated}
\dot{c}_{2n}(t)d_{2n}(\tau)&=\frac{1}{\overline{R}^2}\exp\left(\frac{2-2n(2n+1)}{\overline{R}^{2}}\,t\right)\,b_{2n}'(\tau)+c_{2n}(t)d''_{2n}(\tau).
\end{align}
One special solution $c_{2n}$, $d_{2n}$ of~\eqref{e:inhomoseparated} is then given by
\begin{align}\label{e:overc}
c_{2n}(t)&=\exp\left(\frac{2-2n(2n+1)}{\overline{R}^{2}}\,t\right),
\end{align}
with $d_{2n}$ being any solution of
\begin{align}\label{e:overd}
\frac{2-2n(2n+1)}{\overline{R}^{2}}\,d_{2n}(\tau)-d_{2n}''(\tau)&=\frac{1}{\overline{R}^2}\,b_{2n}'(\tau).
\end{align}
In particular, $\psi_{t}(\tau)=A_{n}c_{2n}(t)d_{2n}(\tau)$ decays for $t\to\infty$ because $n\geq1$. We still need to discuss the inhomogeneous system~\eqref{e:psi_t} corresponding to the mode solution
$\mu_{t}(\tau)=A_{0}\exp\left(\frac{2-\kappa}{\overline{R}^{2}}\,t\right)$. However, as this is spatially constant, it does not contribute to the inhomogeneity of the heat equation~\eqref{e:psi_t}.
\newpage
Combining this with what we found out about the solutions of the
homogeneous equation~\eqref{e:homogeneous}, we can conclude that all
solutions $\psi_{t}$ of the linearised flow equation~\eqref{e:psi_t}
decay. In particular, we found that the most slowly decaying mode of $\psi_{t}$ decays at least as fast as $\exp\left(-\frac{4}{\overline{R}^{2}}t\right)$. 

We have thus asserted full linear stability of the geometric flow~\eqref{e:flow} in a Euclidean background $U\equiv V\equiv0$ around a Euclidean coordinate circle, $r\equiv\overline{R}$, parametrised by arclength. We found that both $\mu_{t}$ and $\psi_{t}$ decay at least as fast as $\exp\left(-\frac{4}{\overline{R}^{2}}t\right)$ or $\exp\left(\frac{2-\kappa}{\overline{R}^{2}}t\right)$ as $t\to\infty$.
In light of the implicit function theorem, this is a strong indication of fully non-linear stability of the geometric flow~\eqref{e:flow} near such circles, as the decay rate is bounded away from zero. Such a non-linear analysis would lead too far here, and we leave it for future work.

We will now establish existence and uniqueness of solutions to~\eqref{e:psi_t} for all times $t>0$ for given initial data $\psi_{0}\colon[0,\overline{L}]\to\mathbb{R}$ and given $\mu_{t}$. Indeed, let $\psi_{0}$ be smooth initial data for~\eqref{e:psi_t}, obeying the boundary~\eqref{e:bdrypsi} and symmetry conditions~\eqref{e:psisymm} for $t=0$, and let $\mu_{t}$ be given as an expansion as in~\eqref{e:mugeneral}. For $n\geq1$, let $c_{2n}$ be as given in \eqref{e:overc} and let $d_{2n}$ be a fixed special solution of~\eqref{e:overd}. Together, these give a solution $\psi_{t}^{\text{inhom}}$ of the inhomogeneous heat equation~\eqref{e:psi_t}, consistent with~\eqref{e:bdrypsi} and  \eqref{e:psisymm},~namely
\begin{align}
\psi_{t}^{\text{inhom}}(\tau)&=\sum_{n=1}^{\infty} A_{n} \exp\left(\frac{2-2n(2n+1)}{\overline{R}^{2}}\,t\right) d_{2n}(\tau)
\end{align}
for $t>0$, $\tau\in[0,\overline{L}]$. In order to determine the constants $C_{n}$ in~\eqref{e:hompsi},
i.e. the homogeneous part of the solution of~\eqref{e:psi_t}, we
proceed as follows. Expand the smooth function
$\psi_{0}-\psi^{\text{inhom}}_{t}\vert_{t=0}$ in terms of the complete
orthonormal system
$\lbrace{\sin\left(\frac{2n}{\overline{R}}\,\tau\right)\,\vert\,n\in\mathbb{N}_{>0}\rbrace}$
of the separable Hilbert space 
$\lbrace{u\in L^{2}[0,\overline{L}]\,\vert\,u(\overline{L}-\tau)=u(\tau)\text{ for all }\tau\in[0,\overline{L}],\;u(0)=u(\overline{L})=0\rbrace}$,
\begin{align}
\psi_{0}(\tau)-\psi^{\text{inhom}}_{t}\vert_{t=0}(\tau)&=\sum_{n=1}^{\infty}C_{n}\sin\left(\frac{2n}{\overline{R}}\,\tau\right)
\end{align}
for $\tau\in[0,\overline{L}]$, with suitable constants $C_{n}\in\mathbb{R}$. Then, the unique solution of~\eqref{e:psi_t} with initial data $\psi_{0}$ and inhomogeneity $\mu_{t}$ is given by
\begin{align}
\psi_{t}(\tau)&=\sum_{n=1}^{\infty}\exp\left(\frac{2-2n(2n+1)}{\overline{R}^{2}}\,t\right) \left\lbrace C_{n}\sin\left(\frac{2n}{\overline{R}}\,\tau\right)+A_{n}d_{2n}(\tau)\right\rbrace
\end{align}
for $t>0$ and $\tau\in[0,\overline{L}]$.

\subsection{Stationary states}\label{s:stationary}
If the flow~\eqref{e:flow} approaches a stationary state, 
$\frac{dx^a_t}{dt} \to 0$ as $t \to \infty$, we immediately know from 
orthonormality of $\lbrace{t^{a}_{t},n^{a}_{t}\rbrace}$ that 
$C_{t}\to 0$, so that the curve $\Gamma_{t}=x^{a}_{t}$ is
asymptotically parametrised proportionally to
arclength by construction (cf.~\eqref{e:C}). As the curve is given on
the interval $[0,\overline{L}]$, there now are two cases: Either we
indeed have $L_\infty := \lim_{t\to\infty }L_{t}=\overline{L}$ 
so that the stationary point
condition together with $\ell\equiv1$ implies 
$H_\infty := \lim_{t\to\infty} H_t = \overline{H}$, so
that by~\eqref{e:embedding} and~\eqref{e:H}, we have indeed found a
free boundary location which induces the correct Bartnik data. In case
$L_\infty\neq\overline{L}$, the functions $H_\infty$
and $\overline{H}$ must differ by an additive constant determined by the size 
of the coupling parameter $\kappa$ and the difference 
$L_\infty-\overline{L}$. 
We cannot rule out that our flow might aproach such unintended
stationary states.
In all the numerical simulations presented in Section
\ref{s:numresults}, we have ensured that $L_T$ agrees with $\overline{L}$
to within numerical error at the final time $T$ of the simulation.
It would be an interesting question for future research to
investigate why the flow generically does seem to approach the desired
stationary state satisfying $L_t \to \overline{L}$ as $t\to\infty$, and if
there any situations when different stationary states are reached.

%% file: paper_sec4.tex
\section{Numerical methods}
\label{s:nummethod}


\subsection{Flow}
\label{s:nummethod.flow}

We solve the flow equation \eqref{e:flow} numerically using a pseudo-spectral method based on Fourier expansions. We will focus on outlining the details of our method; for background material on
pseudo-spectral methods the reader is referred to textbooks such as~\cite{Fornberg1996,Boyd2001}.

All functions of the curve parameter $\tau\in[0,\overline{L}]$ encountered in our implementation
belong to two classes: \emph{even functions} $f(\tau)$ satisfying $f'(0) = f'(\overline{L}) = 0$ and \emph{odd functions} $g(\tau$) satisfying $g(0) = g(\overline{L}) = 0$. In particular, due to the assumed reflection symmetry of the curves $\Gamma_{t}$, see Convention~\ref{conv:reflection}, the spherical polar coordinate $r(\tau)$ is even and the \emph{modified angular coordinate} $\widehat{\theta}(\tau) \definedas \theta(\tau) - \pi\tau/\overline{L}$ is odd.

We expand these functions in truncated Fourier series,
\begin{equation}
  \label{e:spectralexpansion}
  f = \sum_{n=0}^N \widetilde f_n \cos \frac{n\pi\tau}{\overline{L}}, \qquad
  g = \sum_{n=1}^{N-1} \widetilde g_n \sin \frac{n\pi\tau}{\overline{L}}.
\end{equation}
The expansion coefficients $\{\widetilde f_n\}_{n=0}^N$ and $\{\widetilde g_n\}_{n=1}^{N-1}$ constitute one representation of our numerical approximations to $f$ and $g$.

Since the geometric flow \eqref{e:flow} is non-linear, we will need to evaluate non-linear terms numerically. In a pseudo-spectral method, this is done pointwise at a set of collocation points $\tau_j$, which we take to be $\tau_j = j\overline{L}/N$, $0\leqslant j \leqslant N$. We denote by $f_j \definedas f(\tau_j)$ the values of an even function $f$ at these collocation points, similarly for an odd function $g$.

The expansion coefficients and point values are obviously related by 
\begin{align}
  f_j = \sum_{n=0}^N A_{jn} \widetilde f_n, \qquad 
  g_j = \sum_{n=1}^{N-1} B_{jn} \widetilde g_n
\end{align}
with 
\begin{align}
\begin{split}
  A_{jn} &= \cos\frac{jn\pi}{N}, \quad 0\leqslant j,n \leqslant N, \\
  B_{jn} &= \sin\frac{jn\pi}{N}, \quad 1\leqslant j,n \leqslant N-1.
  \end{split}
\end{align}
The inverse transformation is found to be
\begin{align}
  \widetilde f_j = \sum_{n=0}^N (A^{-1})_{jn} f_n, \qquad 
  \widetilde g_j = \sum_{n=1}^{N-1} (B^{-1})_{jn} g_n
\end{align}
with
\begin{align}
\begin{split}
  (A^{-1})_{nj} &= \frac{2 A_{jn}}{N (1 + \delta_{n0})(1 + \delta_{j0})    (1 + \delta_{nN})(1 + \delta_{jN})},\quad 0\leqslant j,n \leqslant N,\\
  (B^{-1})_{nj} &= \qquad\qquad\qquad\quad \frac{2}{N} B_{jn},\qquad\qquad\qquad\quad 1\leqslant j,n \leqslant N-1.
\end{split}
\end{align}

Derivatives of functions can be computed analytically to the given order of the expansion using the known derivatives of the basis functions in \eqref{e:spectralexpansion}. The derivative of an even function $f$ is an odd function $g = f'$ with expansion coefficients
\begin{align*}
  \widetilde g_n = \sum_{m=0}^N C_{nm} \widetilde f_m, \qquad 
  C_{nm} = -\frac{n\pi}{\overline{L}} \delta_{nm}, \quad 
  1\leqslant n \leqslant N-1, \; 0\leqslant m \leqslant N.
\end{align*}
\newpage
Similarly, the derivative $f = g'$ of an odd function $g$ is an even function $f$ computed as
\begin{align*}
  \widetilde f_n = \sum_{m=1}^{N-1} D_{nm} \widetilde g_m, \qquad 
  D_{nm} = \frac{n\pi}{\overline{L}} \delta_{nm}, \quad 
  0\leqslant n \leqslant N, \; 1\leqslant m \leqslant N-1.
\end{align*}

In the code, we find it convenient to represent functions by their point values. Derivatives can be computed directly in point space by combining the transformations discussed above as follows:
\begin{align}
  f'_k &= \sum_{m=1}^{N-1} \sum_{n=0}^N \sum_{j=0}^N 
           B_{km} C_{mn} (A^{-1})_{nj} f_j 
           \asdefined  \sum_{j=0}^N \mathcal{C}^{(1)}_{kj} f_j, \\
  f''_k &= \sum_{p=0}^N \sum_{m=1}^{N-1} \sum_{n=0}^N \sum_{j=0}^N
            A_{kp} D_{pm} C_{mn} (A^{-1})_{nj} f_j 
            \asdefined  \sum_{j=0}^N \mathcal{C}^{(2)}_{kj} f_j,\\
  g'_k &= \sum_{m=0}^N \sum_{n=1}^{N-1} \sum_{j=1}^{N-1} 
           A_{km} D_{mn} (B^{-1})_{nj} g_j 
           \asdefined  \sum_{j=1}^{N-1} \mathcal{D}^{(1)}_{kj} g_j,\\
  g''_k &= \sum_{p=1}^{N-1} \sum_{m=0}^N \sum_{n=1}^{N-1} \sum_{j=1}^{N-1}
            B_{kp} C_{pm} D_{mn} (B^{-1})_{nj} g_j 
            \asdefined  \sum_{j=1}^{N-1} \mathcal{D}^{(2)}_{kj} g_j.
\end{align}
The differentiation matrices $\mathcal{C}^{(1)}$ and $\mathcal{D}^{(2)}$
can be computed once and for all before the simulation starts.

In particular, for the functions $r(\tau)$ and $\theta(\tau)$ representing the curve $\Gamma_{t}$, we have
\begin{align}
  &r'_k = \sum_{j=0}^N \mathcal{C}^{(1)}_{kj} r_j, \qquad 
     r''_k = \sum_{j=0}^N \mathcal{C}^{(2)}_{kj} r_j, \\
  &\theta'_k = \sum_{j=1}^{N-1} \mathcal{D}^{(1)}_{kj} 
     \left(\theta_j - \frac{\pi\tau_j}{\overline{L}} \right) 
     + \frac{\pi}{\overline{L}}, \qquad
  \theta''_k = \sum_{j=1}^{N-1}  \mathcal{D}^{(2)}_{kj} 
     \left(\theta_j - \frac{\pi\tau_j}{\overline{L}} \right).
\end{align}
Occasionally, we will need to divide two odd functions $g$ and $h$ by each other, which results in an even function $f = g/h$. This is done pointwise for $1\leqslant j \leqslant N-1$. At $j=0$ and $j=N$, where the quotient is ill-defined, we apply L'Hospital's rule
\begin{align}
  f_{0} = \left(\frac{g}{h}\right)_{0} = \left(\frac{g'}{h'}\right)_{0}, \quad
  f_{N} = \left(\frac{g}{h}\right)_{N} = \left(\frac{g'}{h'}\right)_{N}.
\end{align}

In order to compute integrals such as \eqref{e:L}, we note that all even 
expansion functions in \eqref{e:spectralexpansion} vanish when 
integrated from $0$ to $\overline{L}$ except for the constant mode $n=0$, 
which integrates to $\overline{L}$.
Thus we have
\begin{align}
  \int_0^{\overline{L}} f(\tau) \diff \tau = \overline{L} \widetilde f_0 
\end{align}
for an even function $f$.

The flow equation \eqref{e:flow} is stepped forward in time using the
Euler forward method
\begin{equation}
  \frac{dx}{dt} = F[x] \; \rightarrow \; x^{n+1} = x^n + \Delta t \, F[x^n]. 
\end{equation}
This method was chosen because of its low computational cost and because the
numerical accuracy of the solution as a function of flow time $t$ does
not matter much to us since we are mainly interested in the asymptotic
behaviour as $t\to\infty$.

For a parabolic equation like \eqref{e:flow}, numerical stability requires
\begin{align*}
\Delta t \leqslant c (\Delta \tau)^2
\end{align*}
for some constant $c$ (the Courant--Friedrichs--Lewy condition~\cite{Courant1928}), where $\Delta \tau = \overline{L}/N$ is the spatial grid spacing. For the simulations presented in Section~\ref{s:numresults}, we will typically set $c=0.1$ and use between $N=30$ and $N=75$ collocation points.

As is typical for pseudo-spectral methods applied to non-linear
partial differential equations, aliasing errors introduce
high-frequency errors that may lead to numerical instabilities.
We address this problem by applying the 2/3 filtering rule~\cite{Boyd2001},
whereby the top third of the spectral coefficients are set to zero
after evaluating the non-linear terms.

In order to construct initial data for the flow \eqref{e:flow}, we often specify the
coordinate location of the initial curve as a function $r(\theta)$,
i.e., the curve parameter is preliminarily taken to be $s=\theta$.
Arclength $\tau$ is then computed according to
\begin{equation}
  \tau(s) = \int_0^s \sqrt{g(\Gamma_t'(s), \Gamma_t'(s))} \, \diff s,
\end{equation}
where $g$ now denotes the Riemannian metric on the boundary surface $\Sigma$ corresponding to the boundary curve $\Gamma_{t}$, see \eqref{e:2dmetric}.

In order to reparametrise the curve by arclength, the function $\tau(s)$ now needs to be inverted numerically in order to obtain $s(\tau)$ and
thus $x^a(s(\tau)) = (r(s(\tau)), \theta(s(\tau)))$.
In practice, this is done by interpolating onto the equidistant
collocation points $\tau_j$.
We have found this procedure to introduce inaccuracies that cause the
``embedding term'' $C$ defined in \eqref{e:C} (which should vanish for
a curve parametrised by arclength) to be unacceptably large. In order to deal with this problem, we apply a few (typically $1000$)
steps of the simplified flow equation
\begin{equation}
  \frac{\diff x_t^a}{\diff t} = C_t t_t^a 
\end{equation}
before the actual simulation starts.
This smoothing procedure is designed to drive $C$ to zero.

The coupling parameter $\kappa$ in the flow equation \eqref{e:flow} is
typically chosen to be $\kappa=4$ unless otherwise noted, see Section~\ref{s:analysis.spherical} for a discussion of the numerical value of $\kappa$.


\subsection{Weyl--Papapetrou equations}
\label{s:nummethod.einstein}

At each timestep of the flow equation \eqref{e:flow}, the
Weyl--Papapetrou equations \eqref{e:Ueqn} and \eqref{e:Veqns}
must be solved for the metric fields $U$ and $V$.
We follow Weyl~\cite{Weyl1917} (see also~\cite{Griffiths2009}) and
expand the solution in spherical harmonics, which in axisymmetry
reduce to the Legendre polynomials $P_n$.

The general solution to the Laplace equation \eqref{e:Ueqn} with asymptotic
boundary condition $U \rightarrow 0$ as $r\rightarrow \infty$ is
\begin{align}
  \label{e:LegendreU}
  U &= - \sum_{n=0}^\infty a_n r^{-(n+1)} P_n (\cos\theta).
\end{align}
From any such solution $U$, the solution to \eqref{e:Veqns} for $V$ 
with $V\rightarrow 0$ as $r\rightarrow \infty$ is obtained as in 
\cite{Griffiths2009} by
\begin{align}
  \label{e:LegendreV}
  V = -\sum_{k=0}^\infty \sum_{l=0}^\infty a_k a_l \frac{(k+1)(l+1)}{(k+l+2)}
  \frac{(P_k P_l - P_{k+1} P_{l+1})}{r^{k+l+2}}.
\end{align}

The coefficients $a_n$ are determined by the inner Dirichlet boundary
conditions \eqref{e:metricWP} and \eqref{e:overlinelambdaWP}.
Numerically, we compute them by truncating the sum in
\eqref{e:LegendreU} at $n=N$ and evaluating the equation at
the $N+1$ points on the curve $(r(\tau_i), \theta(\tau_i))$, 
$0 \leqslant i \leqslant N$, where $\tau_i$ are the collocation points.
This results in an $(N+1)\times(N+1)$ linear system of equations for
the coefficients $a_n$ which we solve using a standard direct linear
solver (\texttt{numpy.linalg.solve} in 
Python\footnote{\url{https://www.python.org},
\url{http://www.numpy.org}, \url{https://scipy.org},
\url{https://matplotlib.org}}, 
which implements the
LAPACK\footnote{\url{http://www.netlib.org/lapack/}} 
routine \texttt{gesv}).

From \eqref{e:LegendreU}, one expects that this numerical procedure
becomes unstable when the curve contains points with radius $r\lesssim
1$, which is indeed what we observe---particularly the higher $a_n$
become very large. 
The problem can be alleviated by only solving \eqref{e:LegendreU}
for the lowest few $a_n$ in the least squares sense 
(we use the routine \texttt{numpy.linalg.lstsq}).
Typically, only the lowest $2/5$ of the coefficients are solved for.
We switch from the full linear solver to the least squares method as
soon as the radius $r$ of one point on the curve becomes smaller than $1.7$.

Once the coefficients $a_n$ in \eqref{e:LegendreU} have been
determined, $V$ is approximated by again truncating the sums in
\eqref{e:LegendreV} at $k,l=N$.
With these numerical approximations to $U$ and $V$, we are also able to
compute all their derivatives at the collocation points on the curve,
as needed for the flow equation \eqref{e:flow}.

The code has been written in Python using the libraries NumPy, SciPy,
and Matplotlib.
A graphical user interface allows the parameters to be
specified and displays plots of various quantities during the flow
(see Section~\ref{s:numresults} for snapshots).
Typical CPU times on a laptop for the simulations with fixed metric
(Section~\ref{s:numresults.fixed}) are one to two minutes, for the
simulations with evolved metric (Section~\ref{s:numresults.evolved})
about ten minutes.

%% file: paper_sec5.tex
\section{Numerical results}
\label{s:numresults}

In this section, we present our numerical evolutions of the flow
equation~\eqref{e:flow} derived in Section~\ref{s:formulation.flow}.

With the exception of the perturbed data considered at the end of
Section~\ref{s:numresults.evolved}, we 
will compute all (Weyl--Papapetrou) Bartnik data 
$([0,\overline{L}],\overline{\lambda},\overline{H})$ 
as well as the initial data for the flow from known asymptotically Euclidean, static, axisymmetric, vacuum solutions in
Weyl--Papapetrou form.
These are taken from the Zipoy--Voorhees family (which includes the
Schwarz\-schild family and Euclidean space) and the
Curzon--Chazy family of solutions.
We present the form of these solutions here; for further details and
the physical interpretation, we refer the reader to~\cite{Griffiths2009}.
\newpage
In the \emph{Zipoy--Voorhees} (or \emph{$\gamma$-metric}) family of solutions, 
the metric functions $U$ and $V$ are given in cylindrical Weyl--Papapetrou coordinates 
$\rho=r\sin\theta$ and $z=r\cos\theta$ by
\begin{align}
 \begin{split} \label{e:zipoy-voorhees}
  U &= \frac{1}{2} \delta \ln \frac{R_+ + R_- - 2M/\delta}
    {R_+ + R_- + 2M/\delta},\\[1ex]
  V &= \frac{1}{2} \delta^2 \ln \frac{(R_+ + R_-)^2 - 4 (M/\delta)^2}{4 R_+ R_-},
\end{split}\end{align}
where
\begin{equation}
  R_\pm = \sqrt{\rho^2 + (z \pm M/\delta)^2},
\end{equation}
and $M,\delta > 0$ are real parameters.
This metric becomes singular at $\rho=0$, $|z|<M/\delta$.
When $\delta\neq 1$, this line segment is a naked curvature singularity
(i.e., there is no event horizon).
However for $\delta=1$ and $M>0$, \eqref{e:zipoy-voorhees} reduces to the 
Schwarzschild solution representing a static spherically
symmetric black hole of mass $M$, and the line segment
$\rho=0$, $|z|<M$ corresponds to the event horizon (which is
\emph{not} a curvature singularity), see also Section~\ref{s:Schwarzschild}. Finally, if we set $\delta = 1$ and $M = 0$, the Weyl--Papapetrou system described by \eqref{e:zipoy-voorhees}
reduces to the Euclidean space, $U=V=0$.

A different family of Weyl--Papapetrou solutions is given by the
\emph{Curzon--Chazy} family,
\begin{align}
  U = -\frac{M}{r}, \qquad V = -\frac{M^2 \sin^2\theta}{2 r^2}.
\end{align}
This solution has a naked curvature singularity at $r=0$.


\subsection{Fixed metric}
\label{s:numresults.fixed}

We begin our analysis of the coupled flow~\eqref{e:flow} and Weyl--Papapetrou system \eqref{e:Ueqn}, \eqref{e:Veqns} by holding the metric fixed, i.e.~we do not solve for the
functions $U$ and $V$ during the flow.
We construct examples of Bartnik data by prescribing the coordinate location 
$\overline{x}^{a}(\tau)$ of the target curve $\Gamma$.
Along this curve, arclength $\tau$ is computed, $0<\tau<\overline{L}$, using
\eqref{e:embedding}.
The function $\overline{\lambda}(\tau)$ is computed from \eqref{e:Ubc} 
and $\overline{H}(\tau)$ is obtained by evaluating \eqref{e:H} on the target curve $\Gamma$
with the given metric functions $U, V$.
For the initial data of the flow, we specify the coordinates of a
different curve. 
Here we have the option to either choose the curve parameter
arbitrarily or to parametrise the curve proportionally to its arclength
in the given metric. (The range of the parameter $\tau$ is $0<\tau<\overline{L}$ at all times along the flow, where $\overline{L}$ is
computed from the target curve as described above.)
The aim now is to check if the flow correctly ``finds'' the specified
target curve in the Weyl--Papapetrou coordinate half-plane and to study how it is approached.
The numerical resolution is taken to be $N=75$ collocation points in
this subsection.

\begin{figure}[t]
  \centering
  \includegraphics[width=0.49\textwidth]{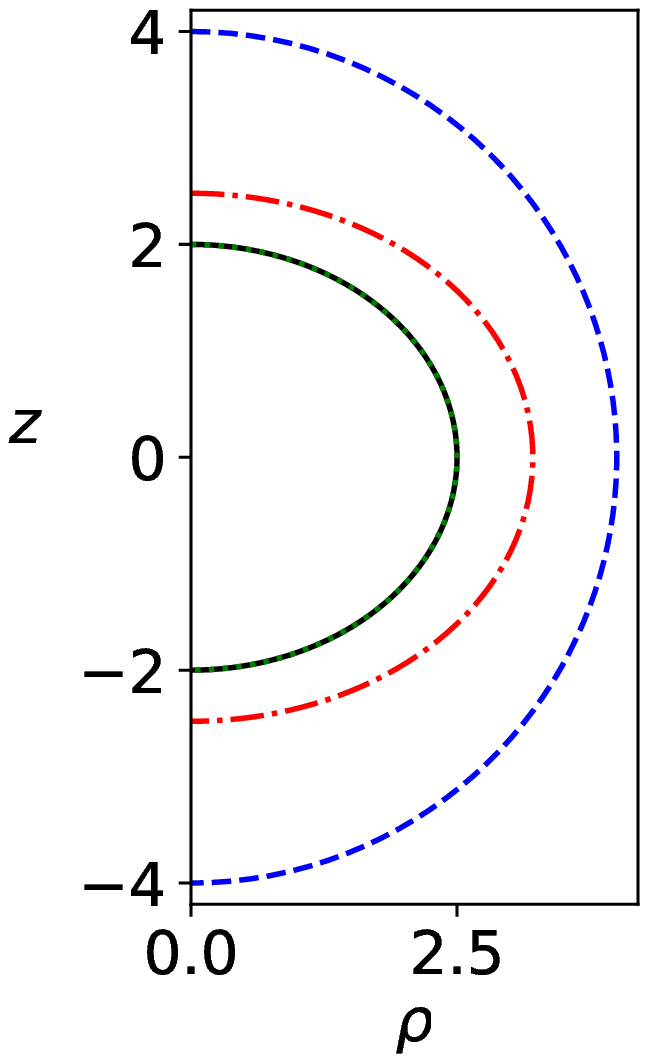}\\
  \includegraphics[width=0.49\textwidth]{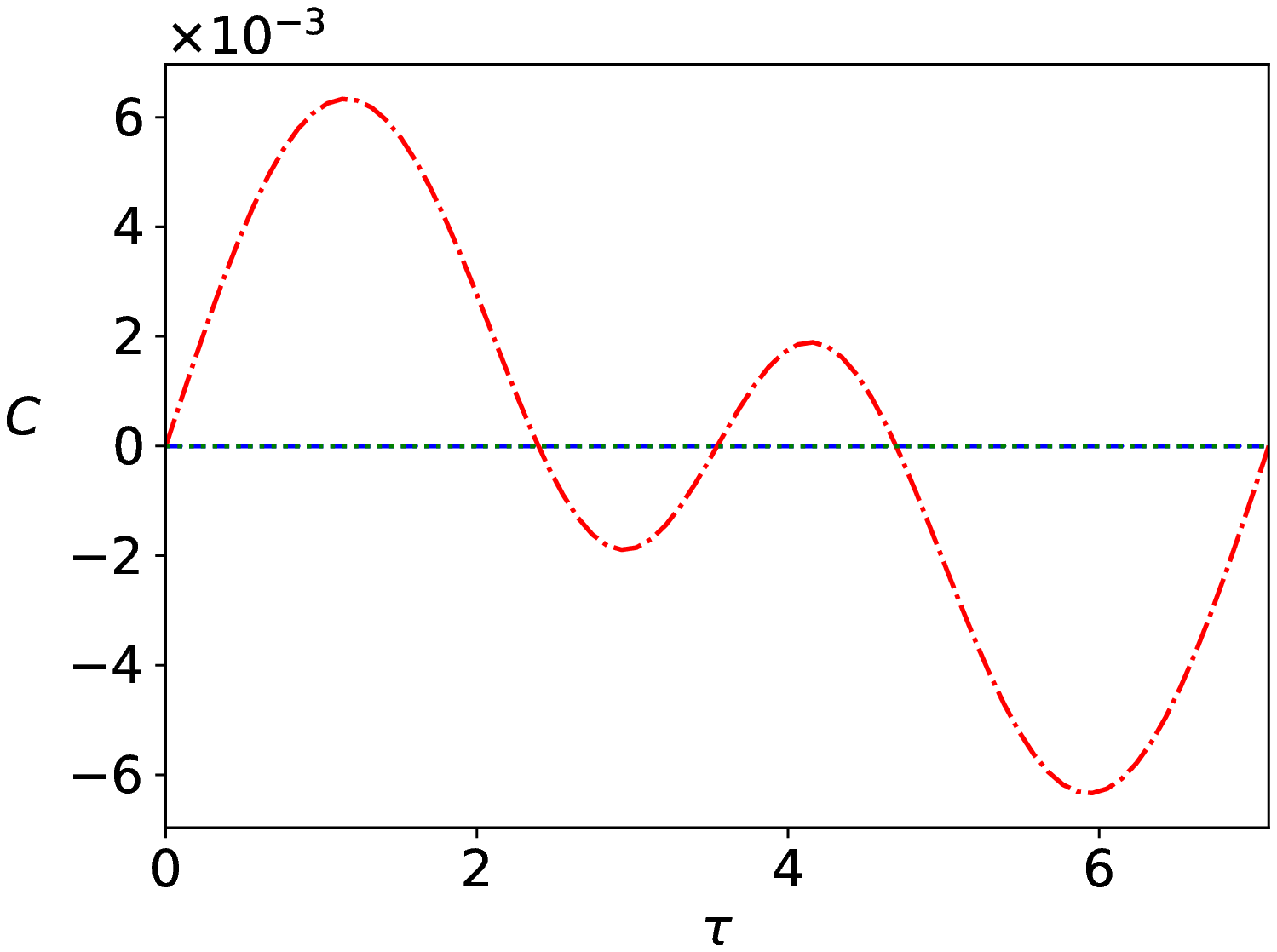}
  \includegraphics[width=0.49\textwidth]{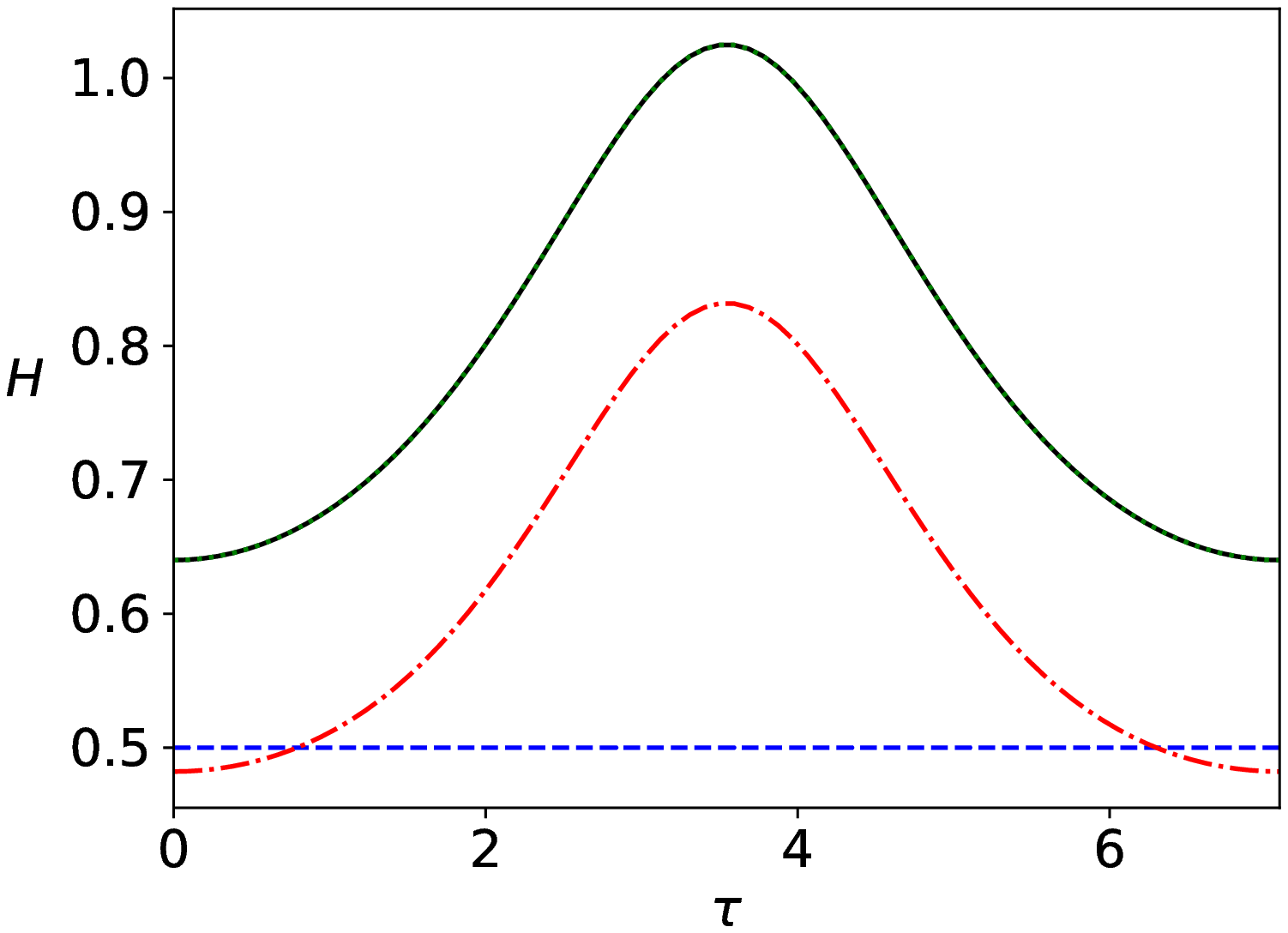}\\
  \includegraphics[width=0.49\textwidth]{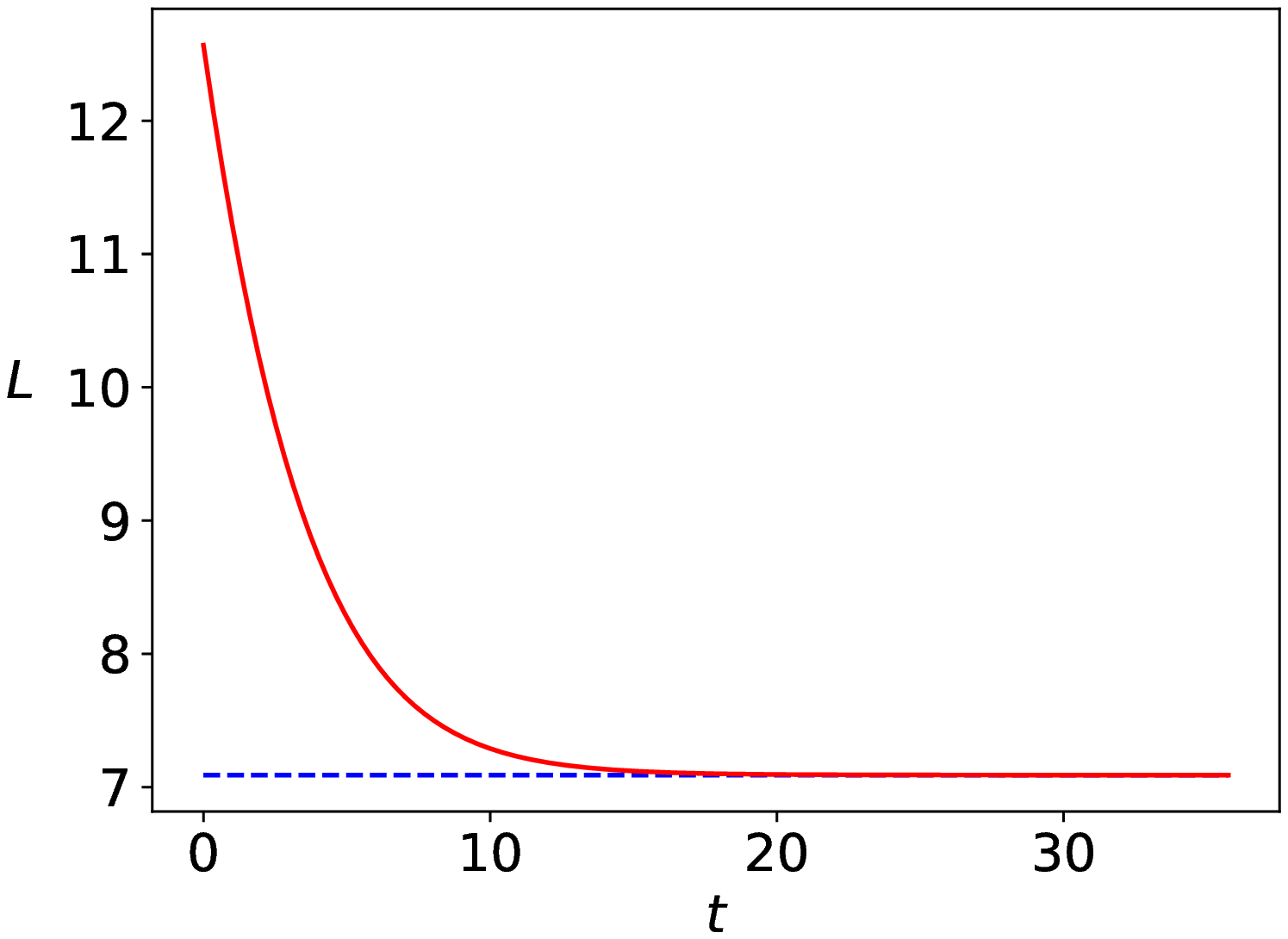}
  \includegraphics[width=0.49\textwidth]{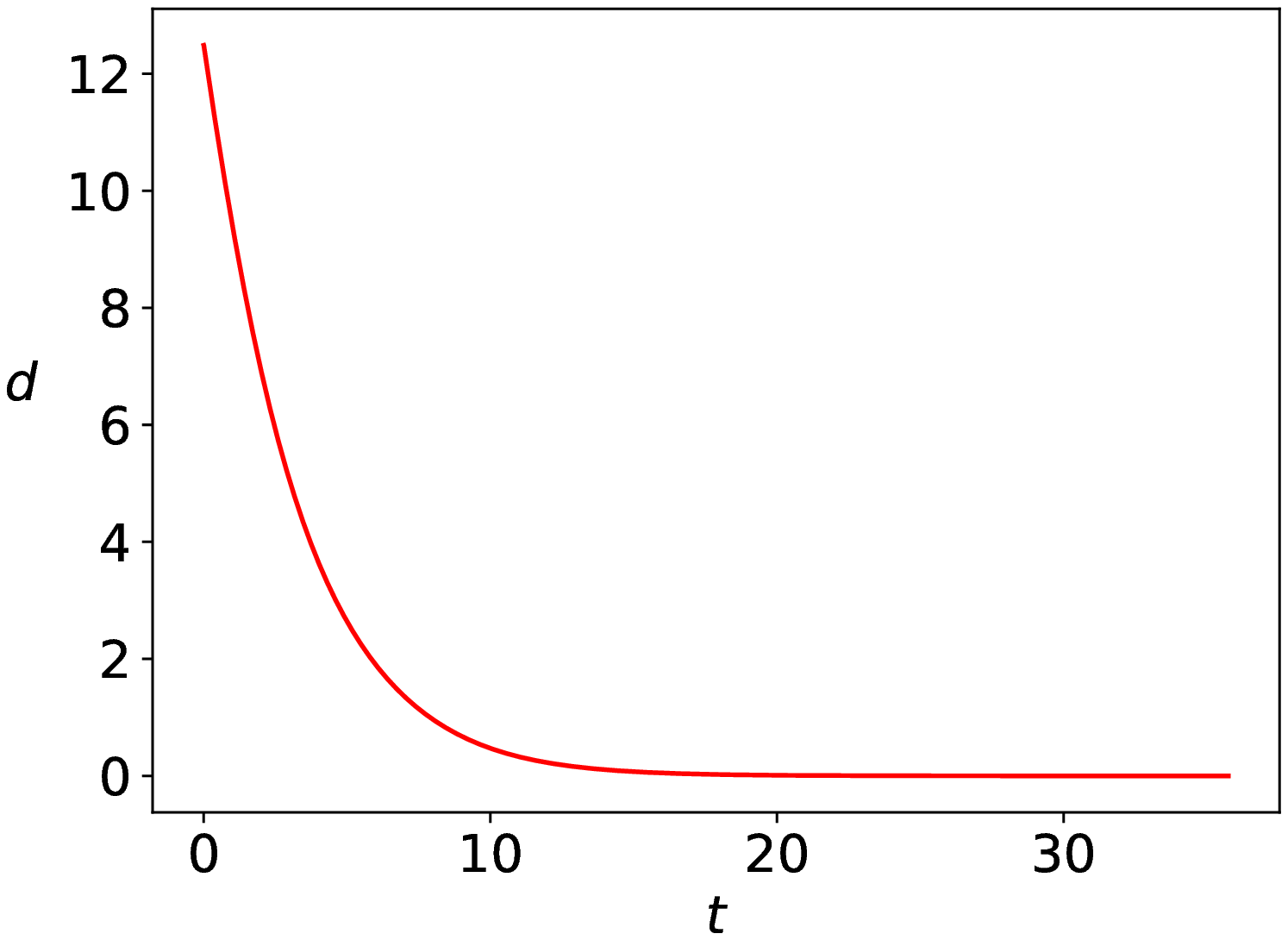}\\
  \caption{Flow with fixed Euclidean background metric.
    The initial curve is taken to be a circle in Weyl--Papapetrou 
    (= Euclidean cylindrical) coordinates 
    (radius $r_0 = 4$, parametrised by arclength) 
    and the target curve an ellipse (semi-major axes 
    $\overline{\rho} = 2.5, \overline{z} = 2$).
    Shown are the coordinate location of the curve (top panel), the
    embedding term $C$ and the mean curvature $H$ (middle panels)
    at flow times $t=0$ (dashed blue), $t=3.6$ (dash-dotted red)
    and $t=35.8$ (dotted green).
    The target curve is plotted in solid black (indistinguishable from
    the dotted green curve here).
    In the bottom panels, we plot the total curve length $L$ 
    (with the target length $\overline{L}$ shown in dashed blue)
    and the Euclidean distance $d$ to the target curve \eqref{e:distance}
    as functions of flow time $t$.
  }
  \label{f:static_Euclidean_circle_to_ellipse}
\end{figure}

First we take the metric to be Euclidean, $U=V=0$. Note that the Weyl--Papapetrou coordinates coincide with the standard 
  Euclidean coordinates in this case. In Figure~\ref{f:static_Euclidean_circle_to_ellipse}, we show the evolution
of an initial circle in Weyl--Papapetrou coordinates to a target curve given by an
ellipse in Weyl--Papapetrou coordinates. 
As expected, the curve approaches the target curve asymptotically as
flow time $t\to\infty$.
The embedding term $C$ defined in \eqref{e:C} vanishes initially (as
the initial circle is parametrised by arclength), departs from zero
during the flow but returns to zero asymptotically, as it should.
The mean curvature $H$ is constant initially and 
approaches its non-constant target profile asymptotically.
The total arclength $L$ (cf.~\eqref{e:length}) approaches its target
value $\overline{L}$ asymptotically.

In order to quantify the approach of the flowing curve $x_t^a(\tau)$ to its
target $\overline{x}^a(\tau)$, we define a distance function
\begin{equation}
  \label{e:distance}
  d(t) := \int_0^{\overline{L}} \lVert x_t^a(\tau) - \overline{x}^a(\tau) \rVert \, \diff\tau,
\end{equation}
where $\lVert \cdot \rVert$ refers to the Euclidean 2-norm.
It is observed to decrease monotonically 
(Figure~\ref{f:static_Euclidean_circle_to_ellipse}). 
We carried out a parameter search over initial and final ellipses
  with various combinations of semi-major axes ranging 
  from $0.1$ to $2.0$, including perturbations of these shapes at the
  $10\%$ level.
  In all cases that led to stable numerical evolutions, $d(t)$ was found
  to be monotonically decreasing.
  We have not been able to prove this in general but our numerical
  results suggest that the distance function \eqref{e:distance}
  might be an interesting quantity to study further in this setting of
  a fixed Euclidean background metric.

\begin{figure}[t]
  \centering
  \includegraphics[width=0.49\textwidth]{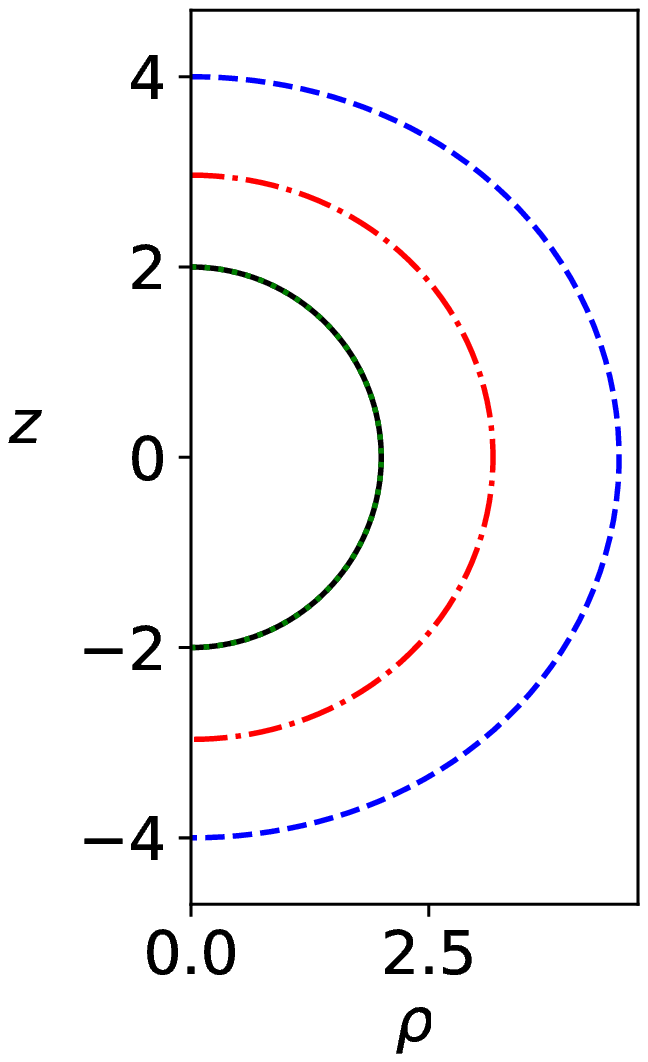}\\
  \includegraphics[width=0.49\textwidth]{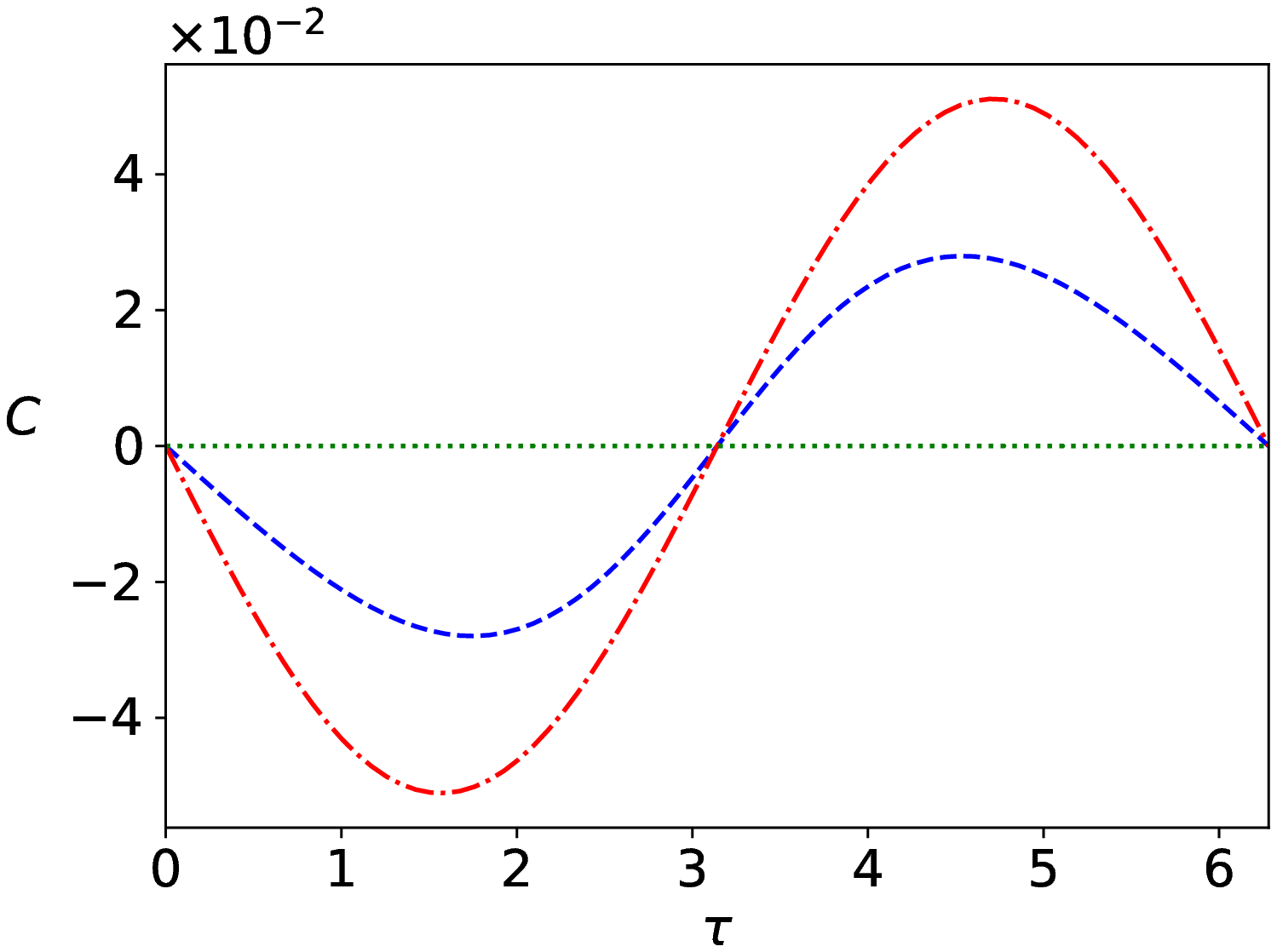}
  \includegraphics[width=0.49\textwidth]{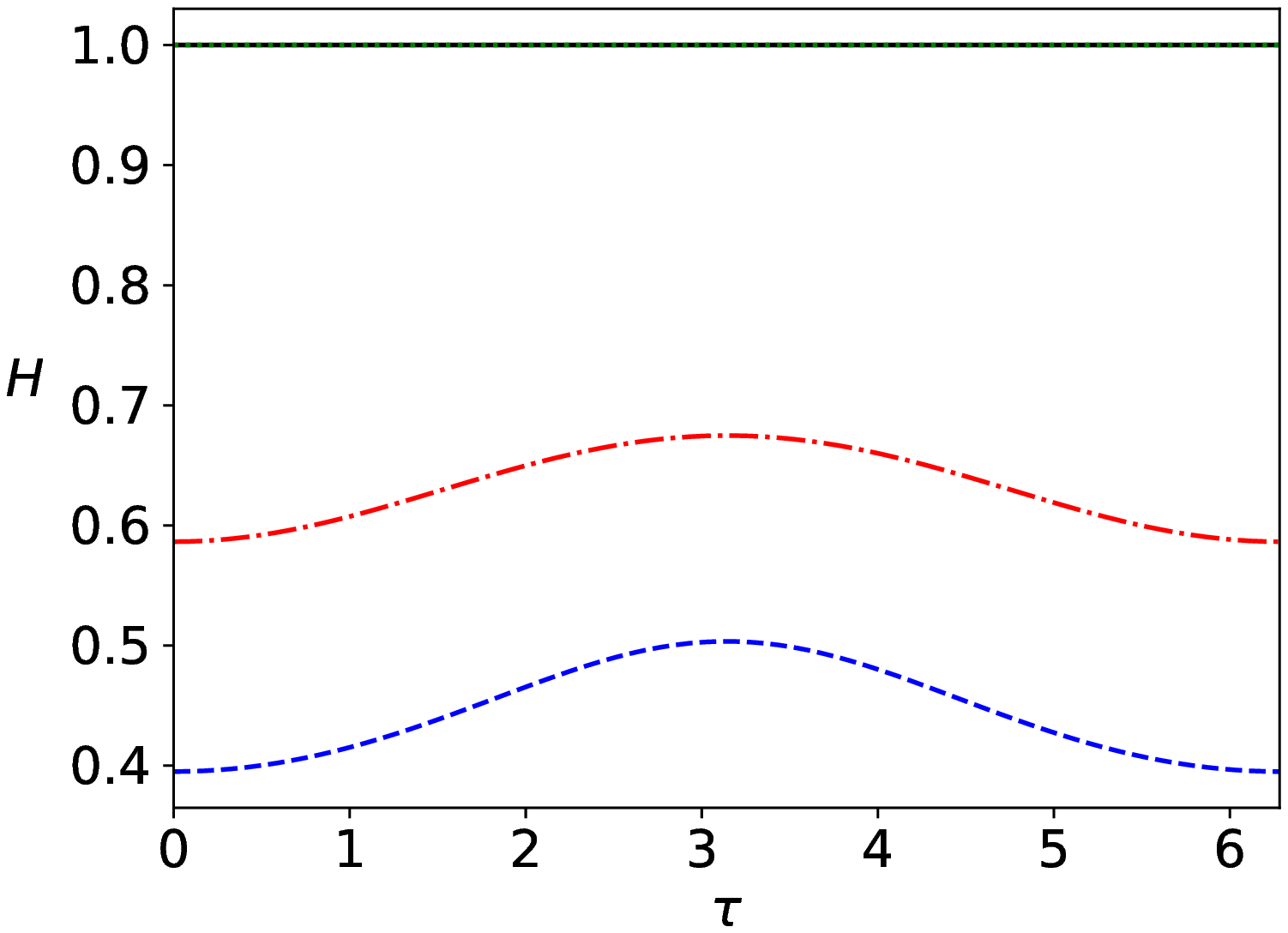}\\
  \includegraphics[width=0.49\textwidth]{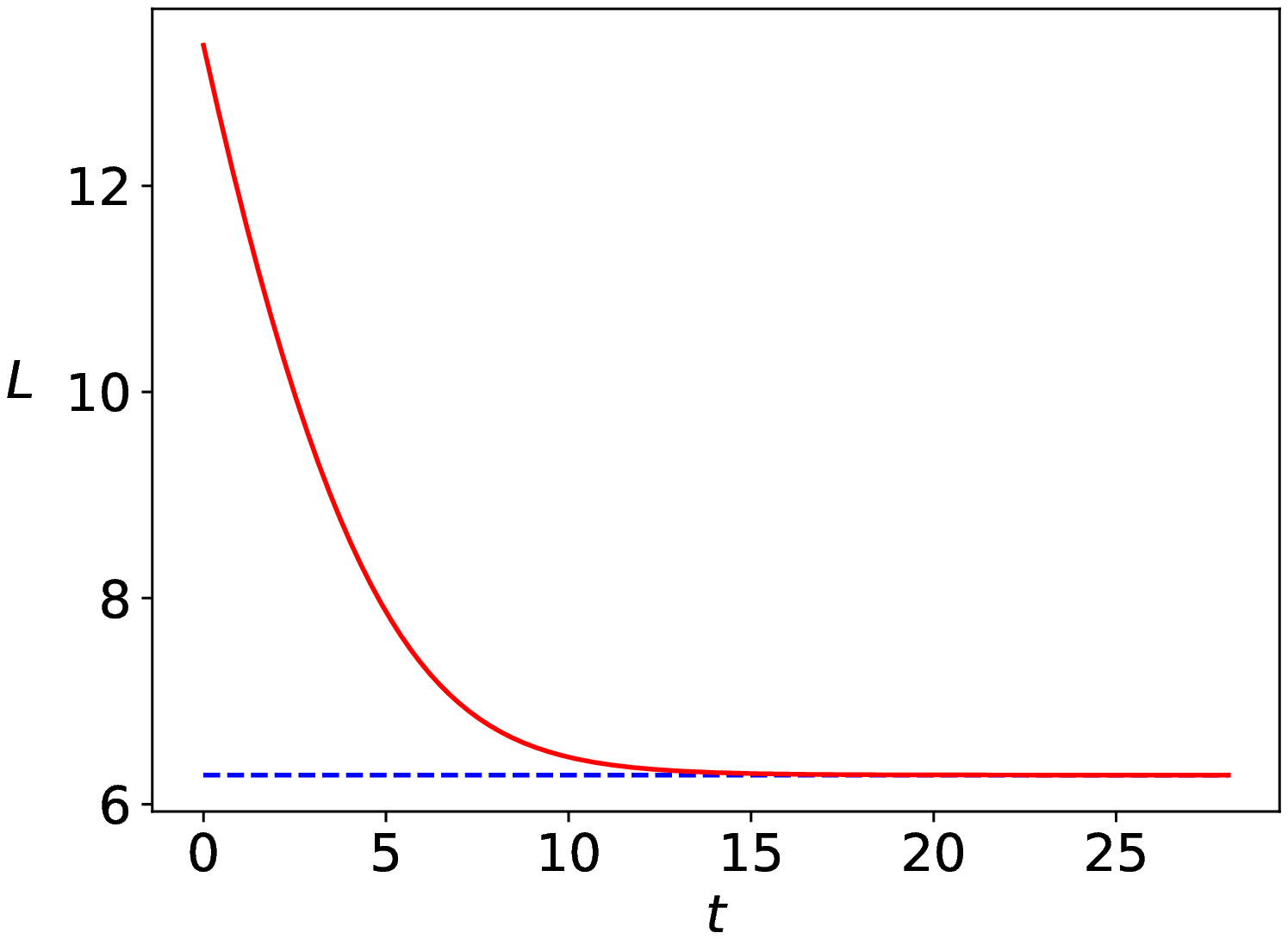}
  \includegraphics[width=0.49\textwidth]{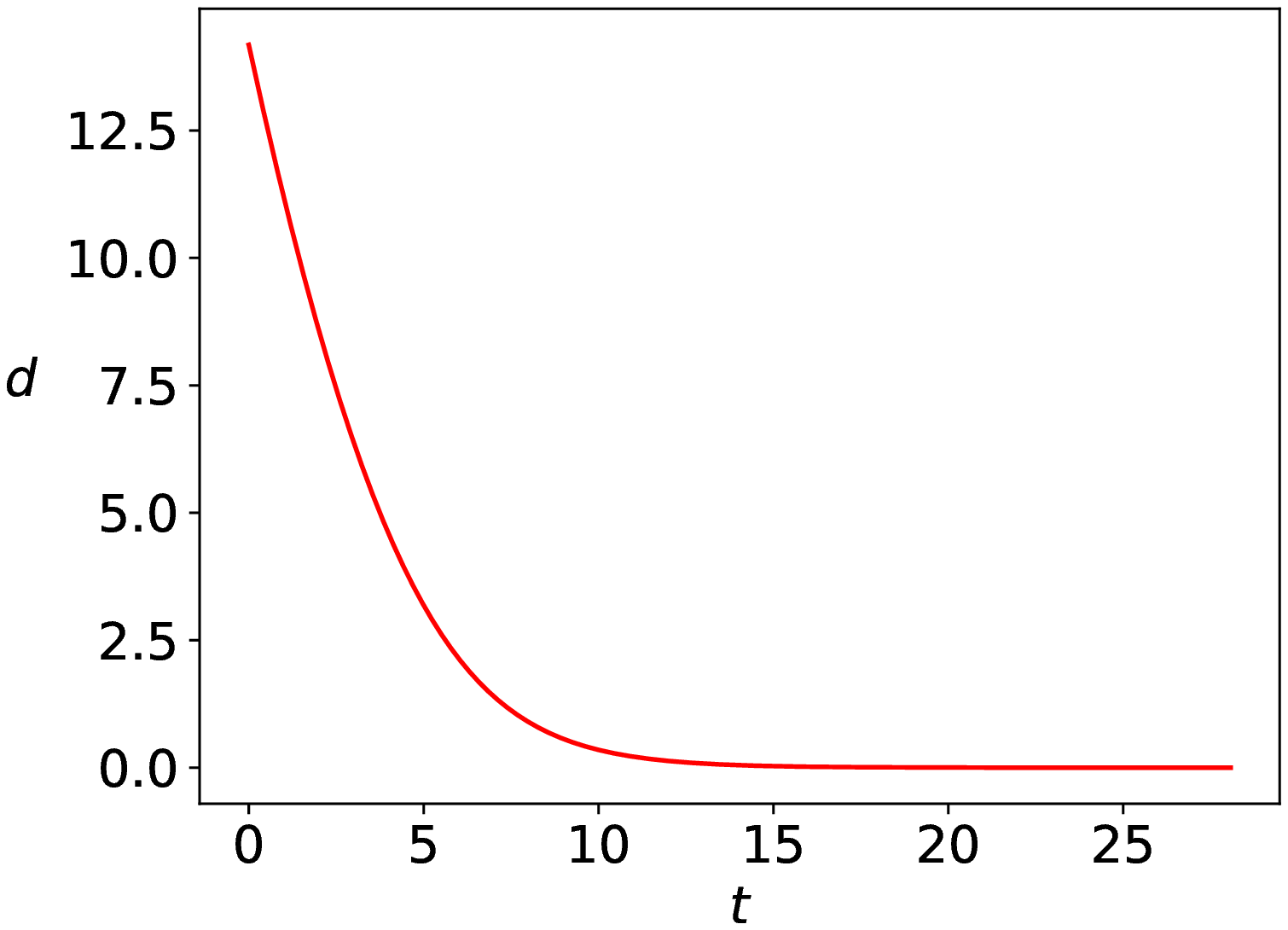}\\
  \caption{Flow with fixed Euclidean background metric.
    The initial curve is taken to be an ellipse in Weyl--Papapetrou  
    (=Euclidean cylindrical) coordinates
    (semi-major axes $\rho_0 = 4.5, z_0 = 4$) which in this case is 
    \emph{not} parametrised by arclength.
    The target curve is a circle (radius $\overline{r} = 2$).
    The same quantities as in Figure~\ref{f:static_Euclidean_circle_to_ellipse}
    are plotted.
    In the first three panels, the curves correspond to flow times $t=0$
    (dashed blue), $t=2.8$ (dash-dotted red) and $t=28.1$ (dotted
    green), with the target solution plotted in solid black.
  }
  \label{f:static_Euclidean_ellipse_to_circle_np}
\end{figure}

In Figure~\ref{f:static_Euclidean_ellipse_to_circle_np} we show a
simulation where the initial curve (an ellipse in Weyl--Papapetrou coordinates) is
\emph{not} parametrised by arclength, which causes $C$ to be non-zero
initially.
Still, the flow converges to the desired curve (a circle in this
case), and $C$ approaches zero, indicating that the final curve
\emph{is} parametrised by arc\-length.

\begin{figure}[t]
  \centering
  \includegraphics[width=0.49\textwidth]{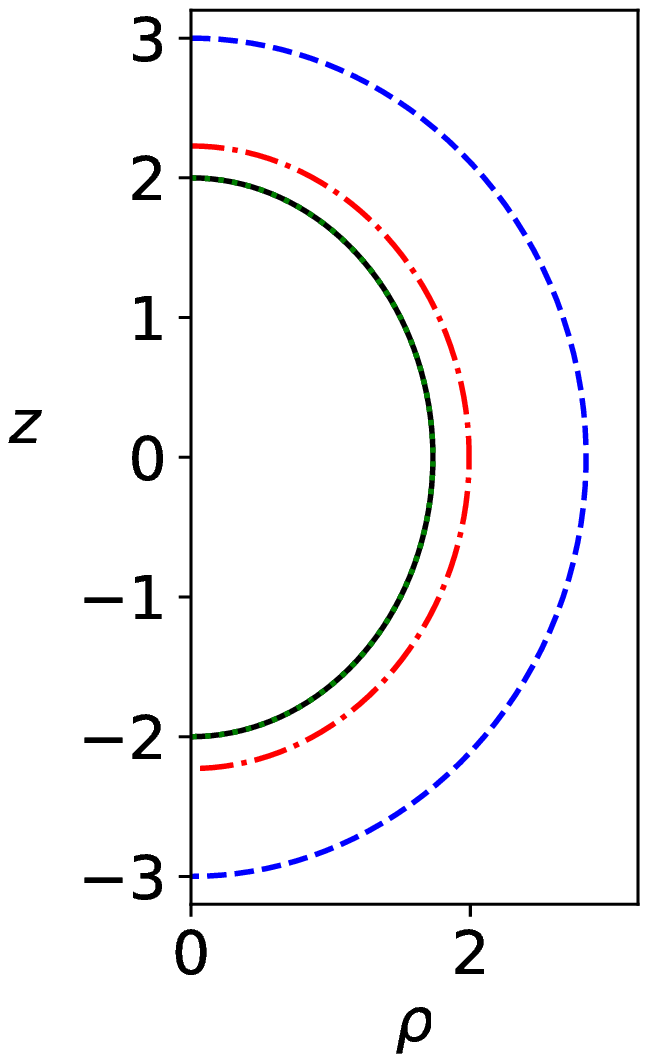}
  \includegraphics[width=0.49\textwidth]{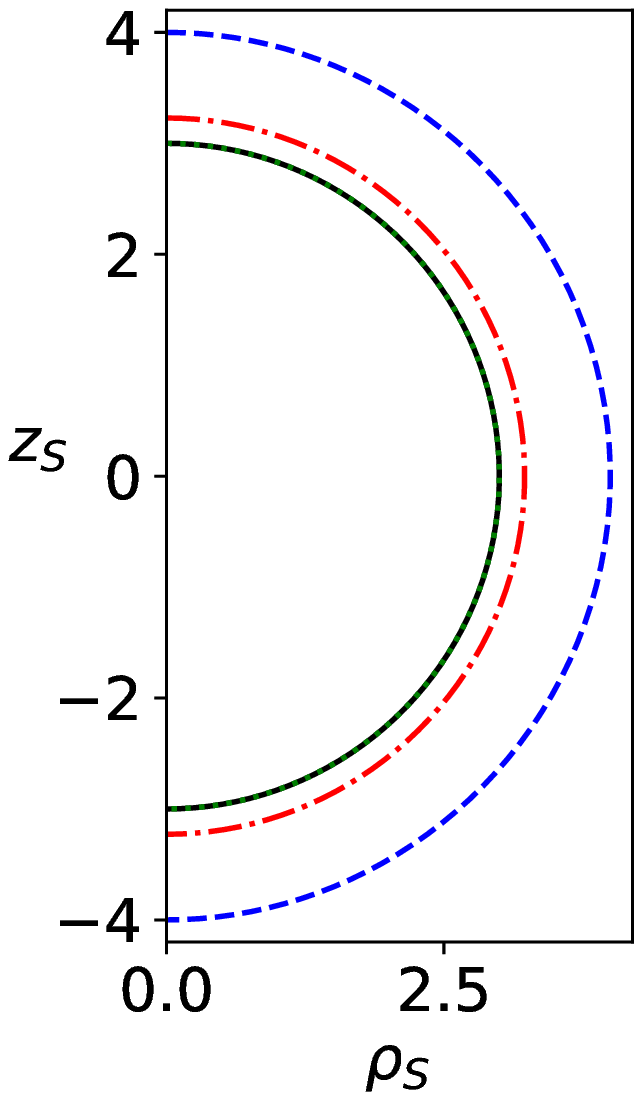}\\
  \includegraphics[width=0.49\textwidth]{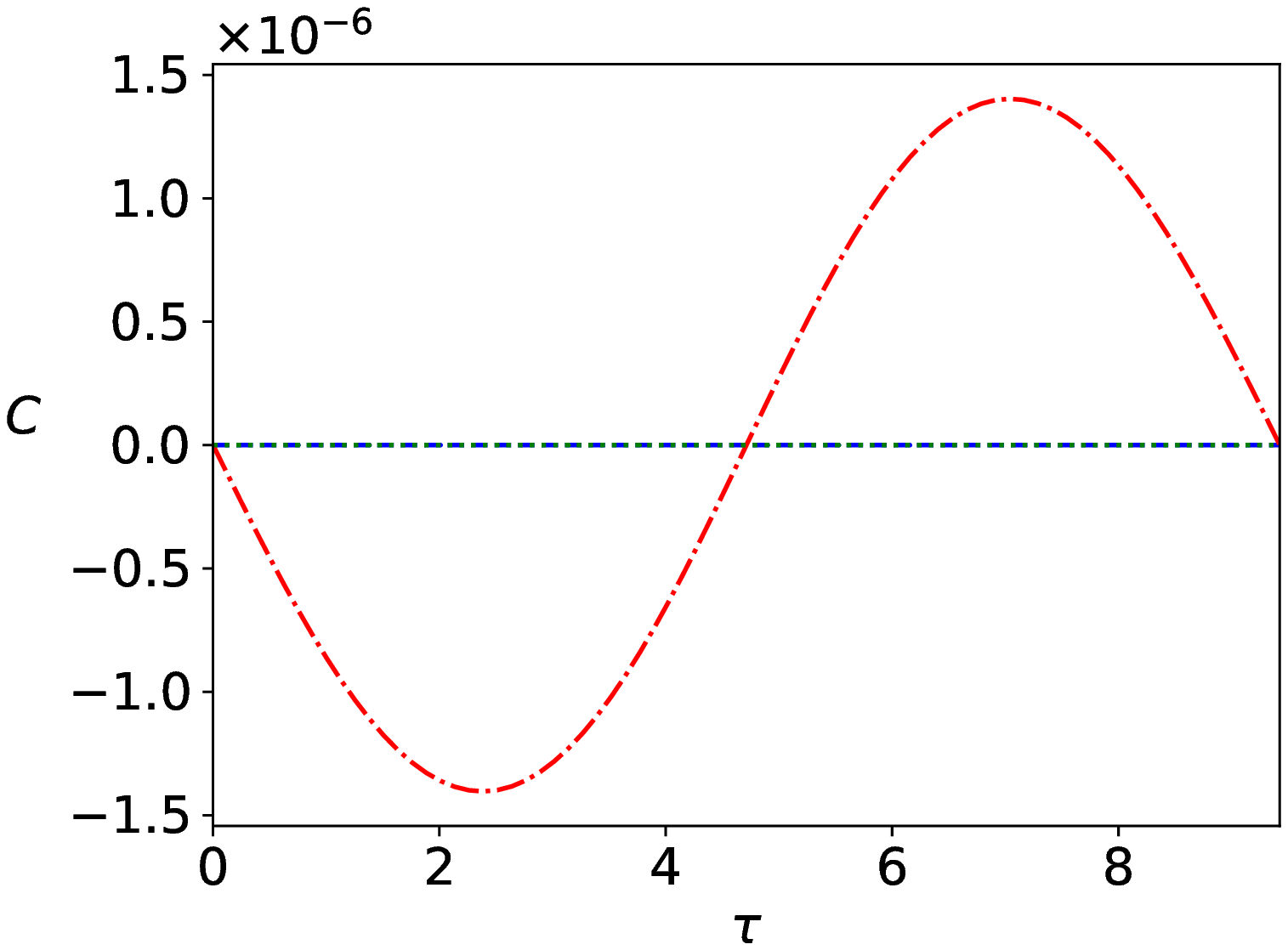}
  \includegraphics[width=0.49\textwidth]{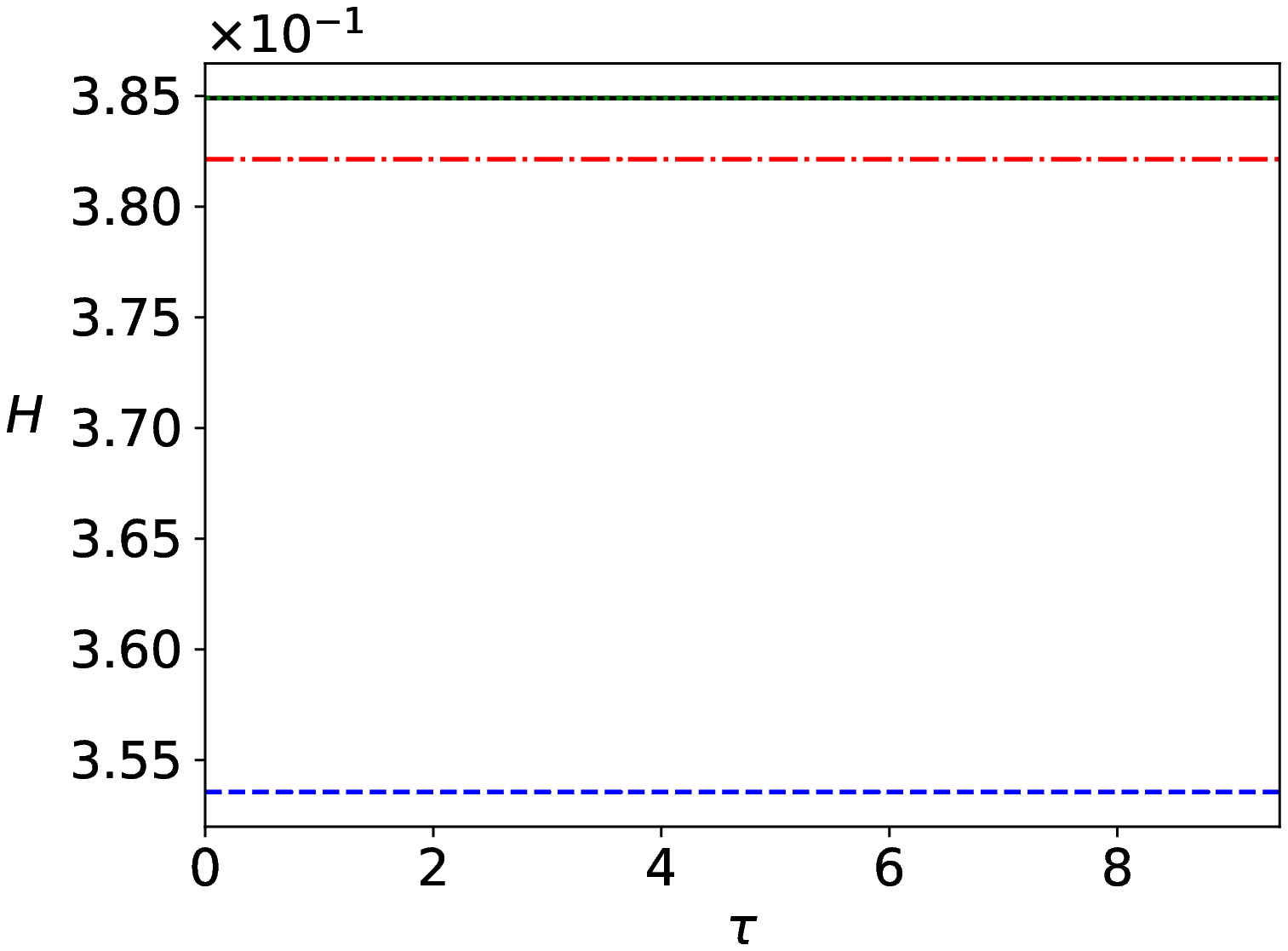}\\
  \includegraphics[width=0.49\textwidth]{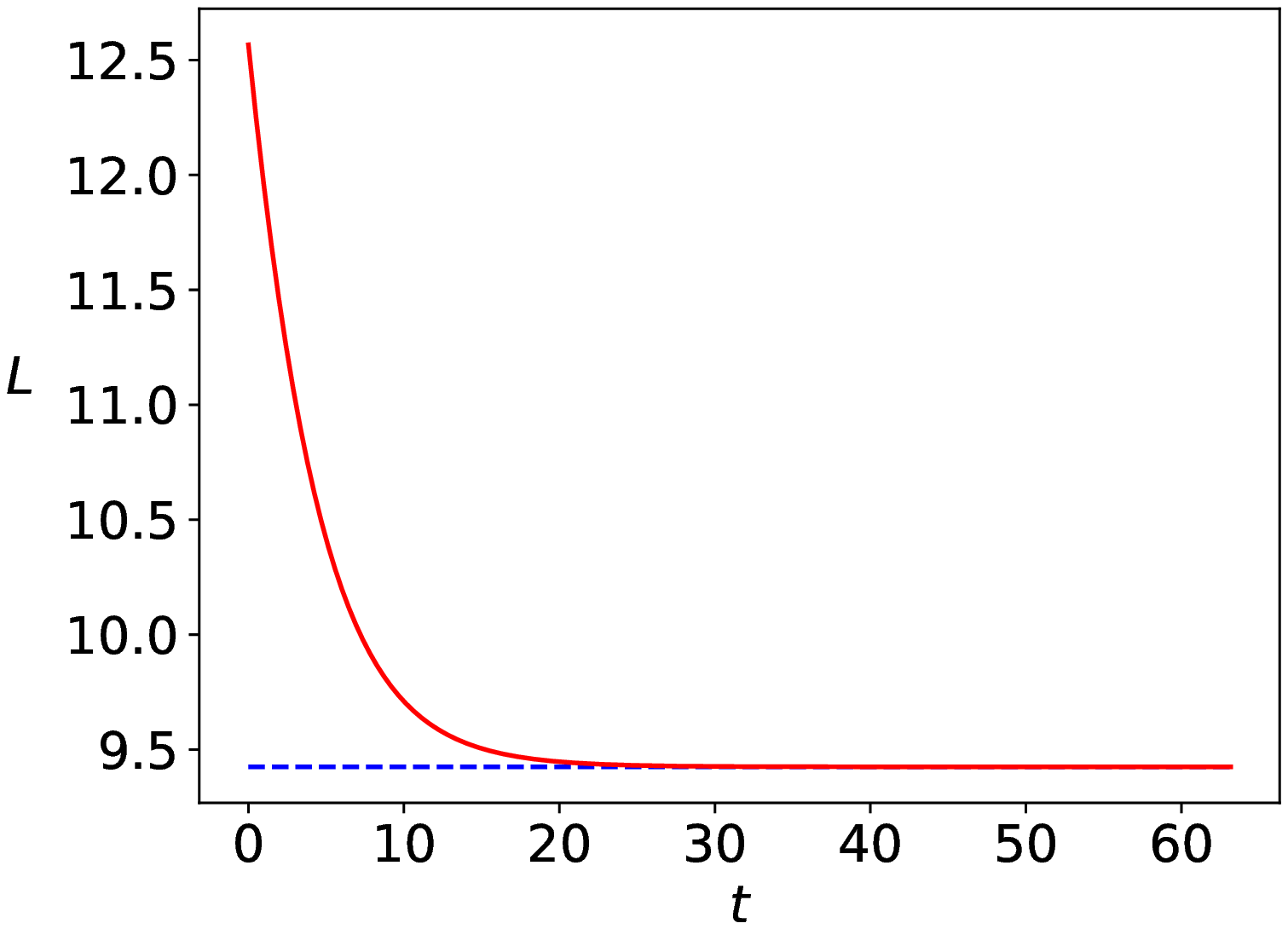}\\
  \caption{Flow with fixed background Schwarzschild metric ($M=1$).
    The initial curve is taken to be a circle in \emph{Schwarzschild} 
    coordinates (radius $r_{S,0} = 4$, parametrised by arclength) 
    and the target curve a smaller circle ($\overline{r}_S = 3$), 
    also in Schwarzschild coordinates.
    Here the top panels show the coordinate location of the curve in Weyl--Papapetrou and
    Schwarzschild coordinates, respectively.
    The remaining quantities are the same as in Figure~\ref{f:static_Euclidean_circle_to_ellipse}.
    In the first four panels, the curves correspond to flow times $t=0$
    (dashed blue), $t=6.3$ (dash-dotted red) and $t=63.2$ (dotted
    green), with the target solution plotted in solid black.
  }
  \label{f:static_ss_circle}
\end{figure}

Next, we take the metric to be a member of the Schwarzschild family.
In Figure~\ref{f:static_ss_circle} we start off with a circle in
Schwarzschild coordinates and let it flow to a circle, again in
Schwarzschild coordinates, with a different radius. 
In accordance with the analysis in Section~\ref{s:analysis.spherical},
the curve remains a coordinate circle in Schwarzschild coordinates
during the entire flow.

\begin{figure}[t]
  \centering
  \includegraphics[width=0.49\textwidth]{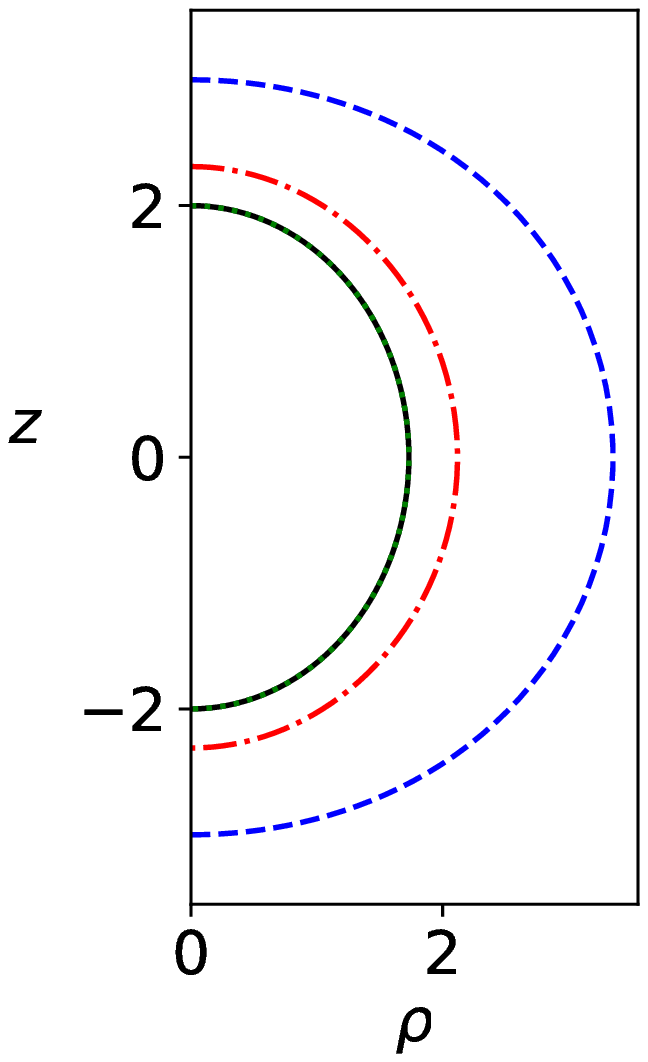}
  \includegraphics[width=0.49\textwidth]{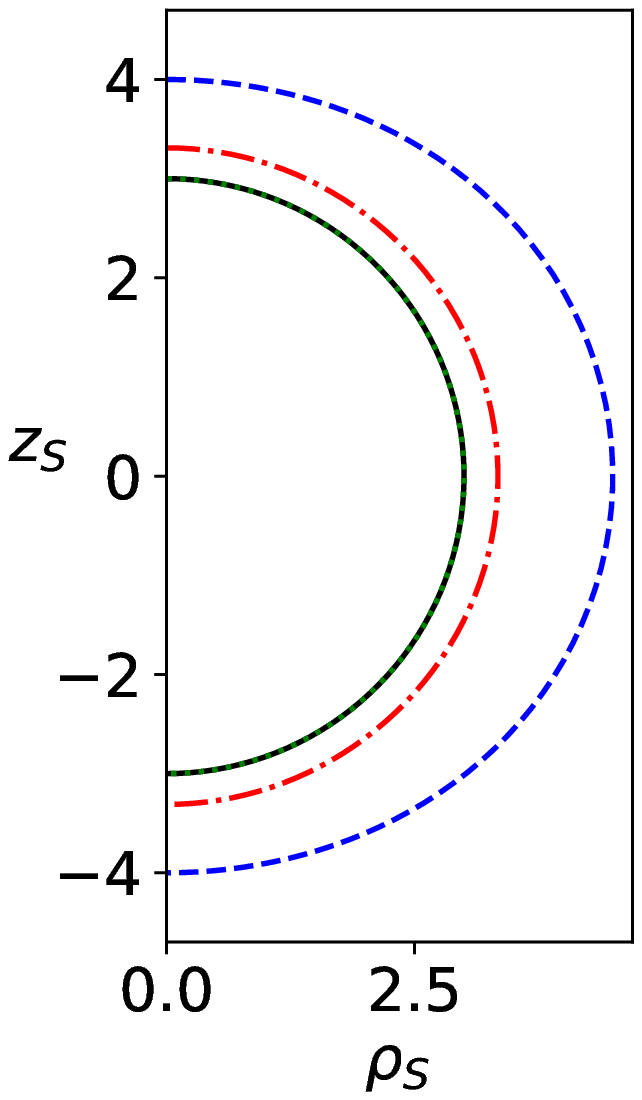}\\
  \includegraphics[width=0.49\textwidth]{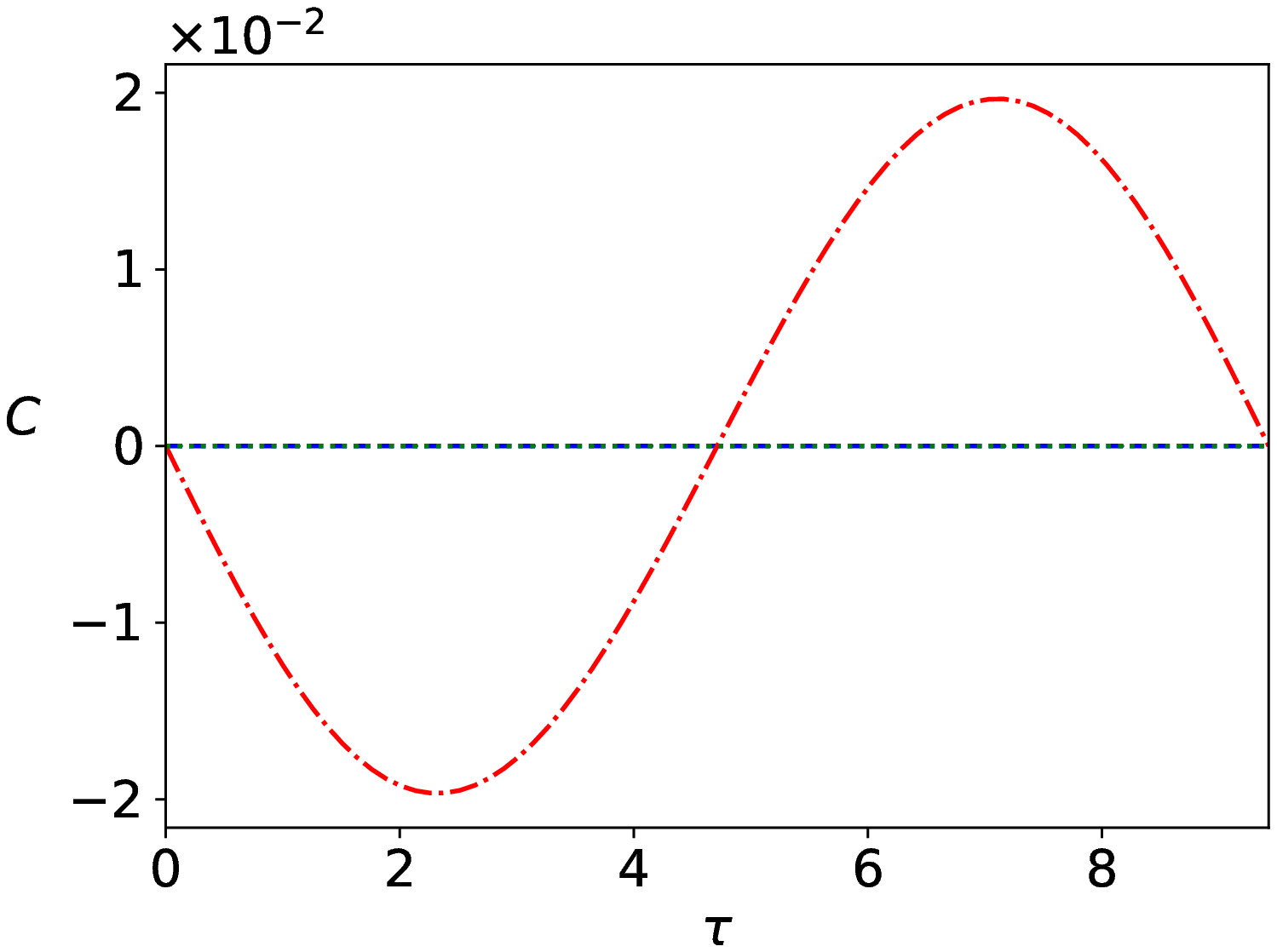}
  \includegraphics[width=0.49\textwidth]{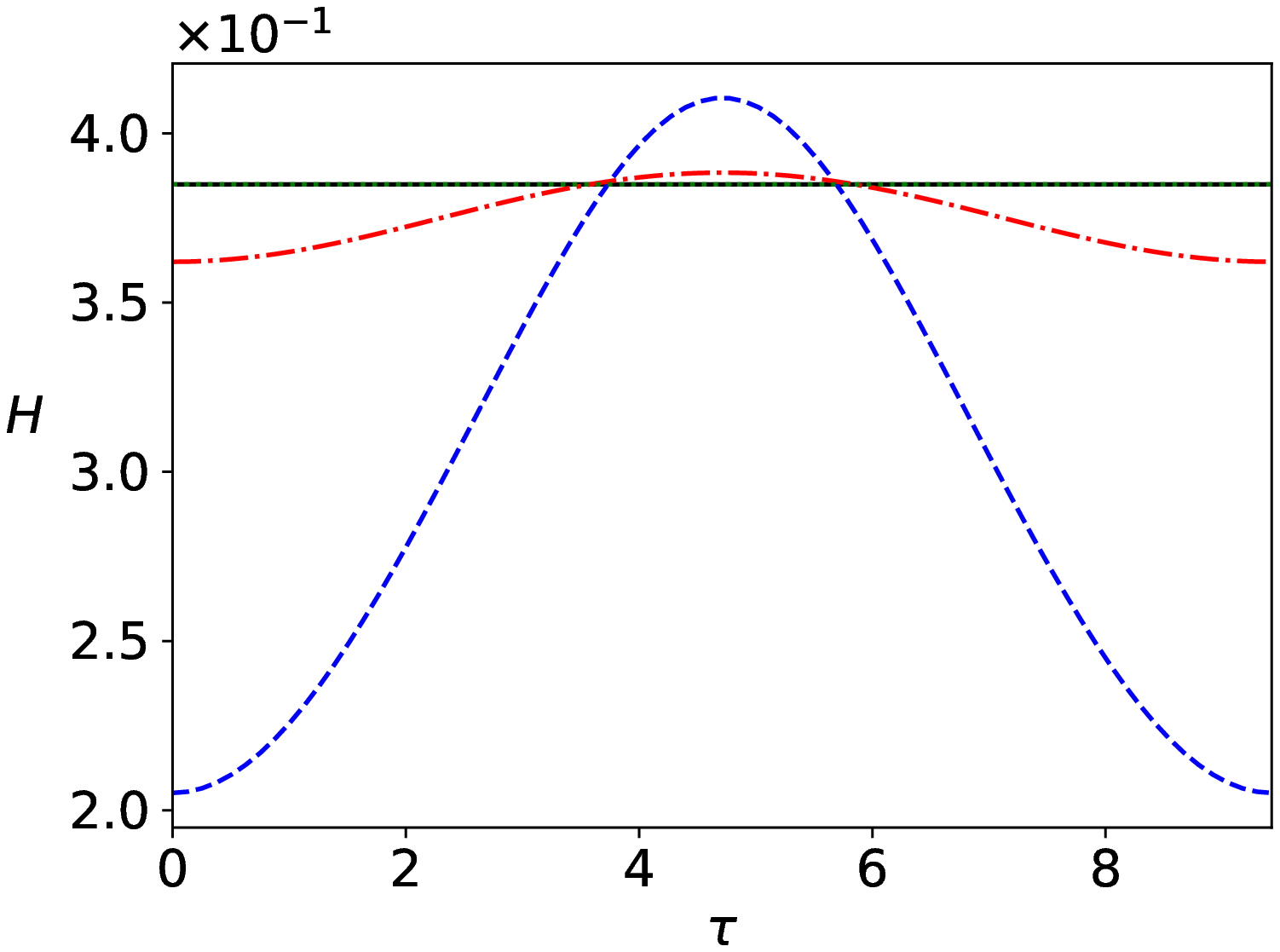}\\
  \includegraphics[width=0.49\textwidth]{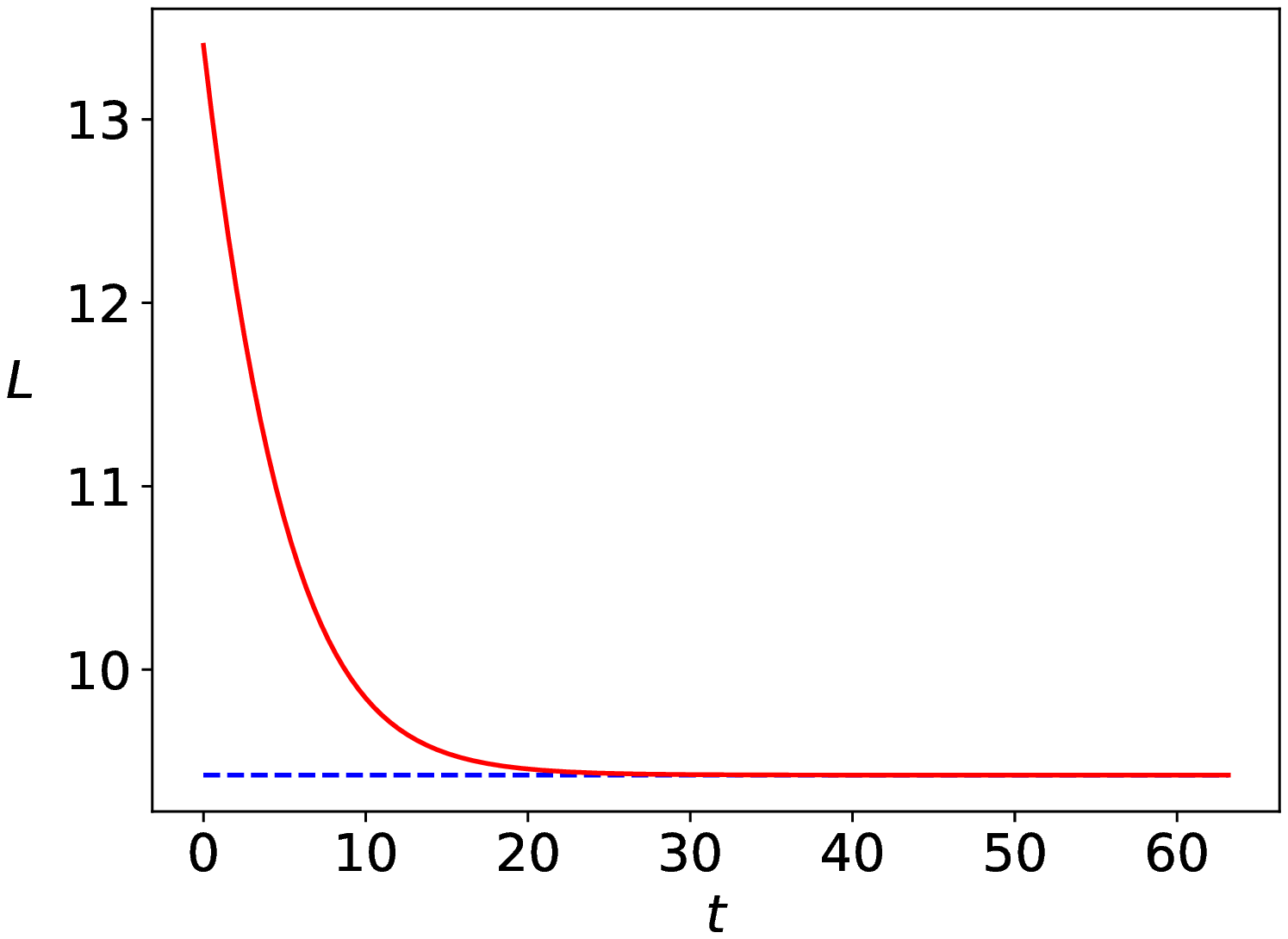}\\
  \caption{Flow with fixed background Schwarzschild metric ($M=1$).
    The initial curve is taken to be an ellipse in Schwarzschild 
    coordinates (semi-major axes $\rho_{S,0} = 4.5, z_{S,0} = 4$, 
    parametrised by arclength) and the target curve a circle ($\overline{r}_S = 3$), also in Schwarzschild coordinates.
    The same quantities as in Figure~\ref{f:static_ss_circle} are plotted.
    In the first four panels, the curves correspond to flow times $t=0$
    (dashed blue), $t=6.3$ (dash-dotted red) and $t=63.2$ (dotted
    green), with the target solution plotted in solid black.
  }
  \label{f:static_ss_ellipse}
\end{figure}

The flow also correctly finds the target curve if we start off with an
initial curve that is not a Schwarzschild coordinate circle (Figure
\ref{f:static_ss_ellipse}).

\begin{figure}[t]
  \centering
  \includegraphics[width=0.49\textwidth]{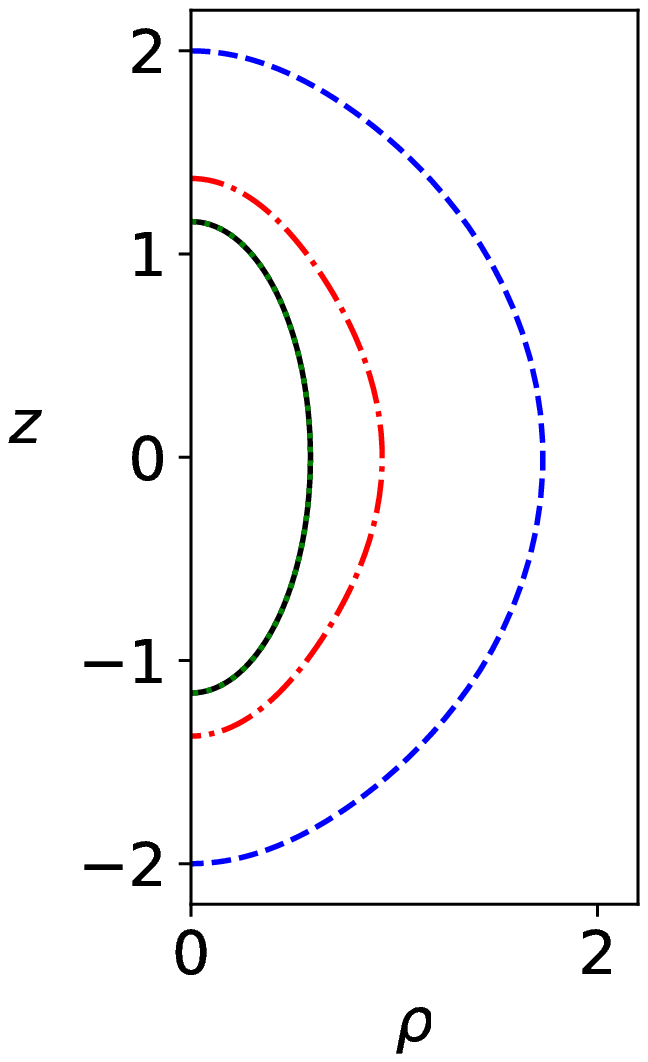}
  \includegraphics[width=0.49\textwidth]{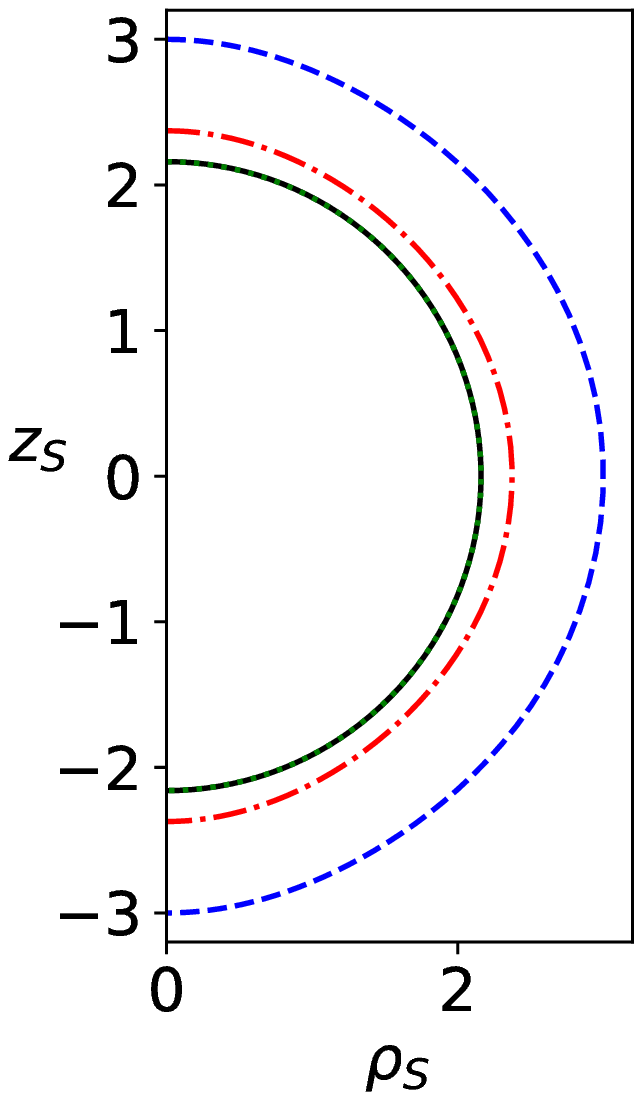}\\
  \includegraphics[width=0.49\textwidth]{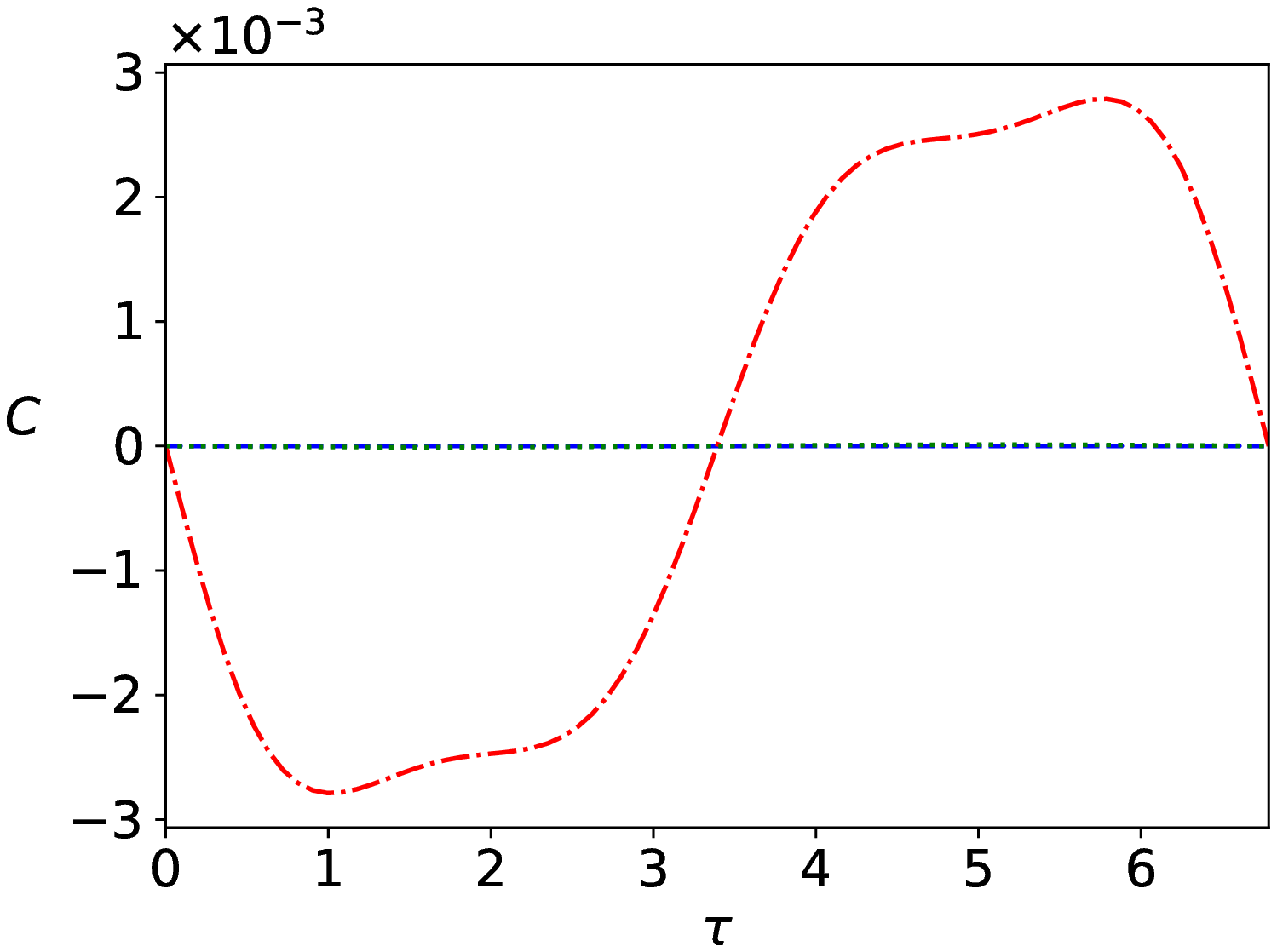}
  \includegraphics[width=0.49\textwidth]{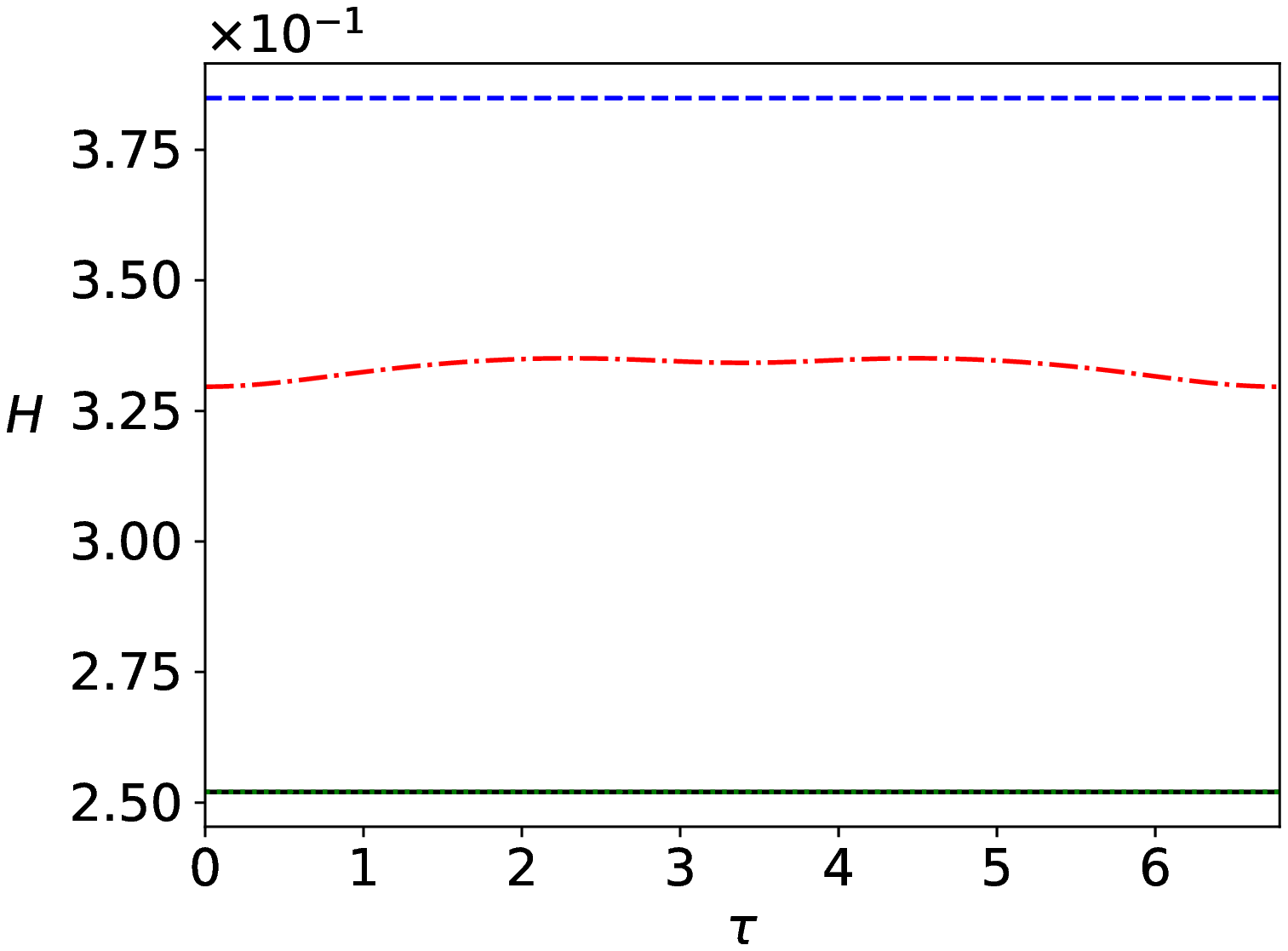}\\
  \includegraphics[width=0.49\textwidth]{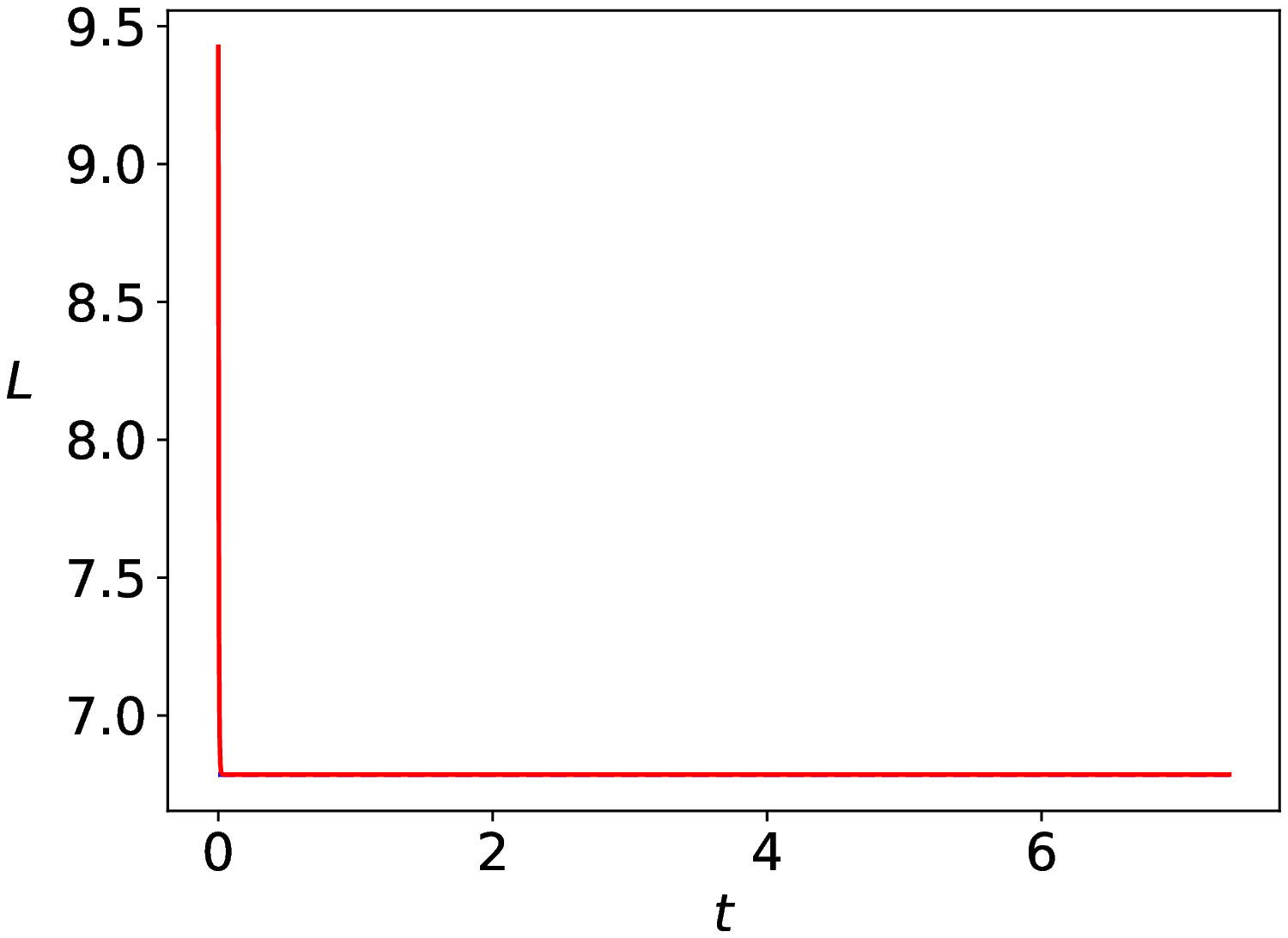}\\
  \caption{Flow with fixed background Schwarzschild metric ($M=1$).
    The initial curve is taken to be a circle in Schwarzschild 
    coordinates (radius $r_{S,0} = 3$, parametrised by arclength) 
    and the target curve a circle close to the horizon
    ($\overline{r}_S = 2.16$), also in Schwarzschild coordinates.
    The same quantities as in Figure~\ref{f:static_ss_circle} are
    plotted.
    In the first four panels, the curves correspond to flow times $t=0$
    (dashed blue), $t=0.0041$ (dash-dotted red) and $t=7.4$ (dotted
    green), with the target solution plotted in solid black.
    The total curve length $L$ decreases rapidly to its target
    value on the time scale shown.
  }
  \label{f:static_ss_circle_extremal}
\end{figure}

Our numerical method breaks down when the flowing curve gets too
close to the horizon, which as discussed above degenerates to a line
on the $\rho=0$ axis in Weyl--Papapetrou coordinates.
In Figure~\ref{f:static_ss_circle_extremal} we choose as a target
curve a circle in Schwarzschild coordinates that is as close to the
horizon as we can get ($\overline{r}_S = 2.16 M$).
In this case we had to increase the value of $\kappa$ from $4$ to $4000$
and decrease the time step from $0.1$ to $0.01$ in order to obtain a 
stable numerical evolution.

We have successfully tested the flow on the Zipoy--Voorhees and
Curzon--Chazy backgrounds as well, with similar results.

\clearpage


\subsection{Evolving metric}
\label{s:numresults.evolved}

In this section, we let the metric evolve along with the flow by
solving the Weyl--Papapetrou equations \eqref{e:Ueqn}, \eqref{e:Veqns} for $U$ and $V$ as 
described in Section~\ref{s:nummethod.einstein}.
The numerical resolution is taken to be $N=30$ collocation points throughout
this subsection.

In order to obtain some insight into how the static metric extensions
change during
the flow, we compute three masses at each time step, namely the
(total) ADM mass $m_\mathrm{ADM}$ \eqref{def:mADM}, 
the (quasi-local) Hawking mass $m_\mathrm{H}$ \eqref{def:mH}, 
and the pseudo-Newtonian mass $m_\mathrm{PN}$ \eqref{def:mPN} 
of the boundary surface corresponding to the flowing curve 
$\Gamma_{t}=x^{a}_{t}(\tau)$ in the static metric extension
corresponding to $U_{t}$ and $V_{t}$. 
We remind the reader that the relations between these masses were 
briefly discussed in Section~\ref{s:masses}.

As the ADM mass can easily be seen to be the leading order term in an expansion of $U$ in inverse
powers of $r$,
\begin{align}
  U &= -\frac{m_\mathrm{ADM}}{r} + \mathcal{O}\left(\frac{1}{r^2}\right),
\end{align}
and can thus be computed as the first expansion coefficient in the Legendre 
expansion \eqref{e:LegendreU} that we use to solve the
Weyl--Papapetrou equations:
\begin{align}
m_\mathrm{ADM} &= a_0.
\end{align}
Combining the definition of the Hawking mass~\eqref{def:mH} with our expression for the mean curvature~\eqref{e:H} and for the induced $2$-metric~\eqref{e:2dmetric}, the Hawking mass of the boundary surface corresponding to $\Gamma_{t}=x^{a}_{t}(\tau)$ is obtained as
\begin{align}
\begin{split}
  m_\mathrm{H} &= \sqrt{\frac{1}{8} \int_0^{\overline{L}} \ell \, r \,
    \sin\theta\, e^{-U\circ \gamma} \, \diff \tau}\\
  &\qquad\times\left( 1- \frac{1}{8} \int_0^{\overline{L}} H^2 \, \ell \, r \,
    \sin\theta\, e^{-U\circ \gamma} \, \diff \tau \right).
    \end{split}
\end{align}
Similarly, the pseudo-Newtonian mass~\eqref{def:mPN} of the boundary surface corresponding to $\Gamma_{t}=x^{a}_{t}(\tau)$ can be computed to be
\begin{align}
  m_\mathrm{PN} &= \frac{1}{2} \int_0^{\overline{L}} (r^2 \theta'\,
  U_{,r} \circ \gamma - r' \, U_{,\theta} \circ\gamma) \sin\theta \, \diff \tau,
\end{align}
where we used \eqref{e:lapse} and \eqref{e:2dmetric}.

In the following figures, the masses are shown as functions of flow
time $t$ in the bottom right panel: the ADM mass (solid blue), 
the Hawking mass (dashed green) and the pseudo-Newtonian mass (dash-dotted red).
In many plots the curves coincide. 
Recall from Section~\ref{s:masses} that the ADM and the
pseudo-Newtonian mass must theoretically 
be identical, 
while the Hawking mass will in general be smaller than the other two masses, except on round spheres in Euclidean space and on centred round spheres in Schwarzschild, where all three notions of mass coincide.

Our numerical analysis for the coupled flow--Weyl--Papapetrou system now
proceeds as follows. In this subsection, we construct Bartnik data by
specifying the coordinate location of the target curve in a given
background metric, as in Section~\ref{s:numresults.fixed}.
 (We will consider more general Bartnik data in Section~\ref{s:perturbed}.)

\begin{figure}[t]
  \centering
\includegraphics[width=0.49\textwidth]{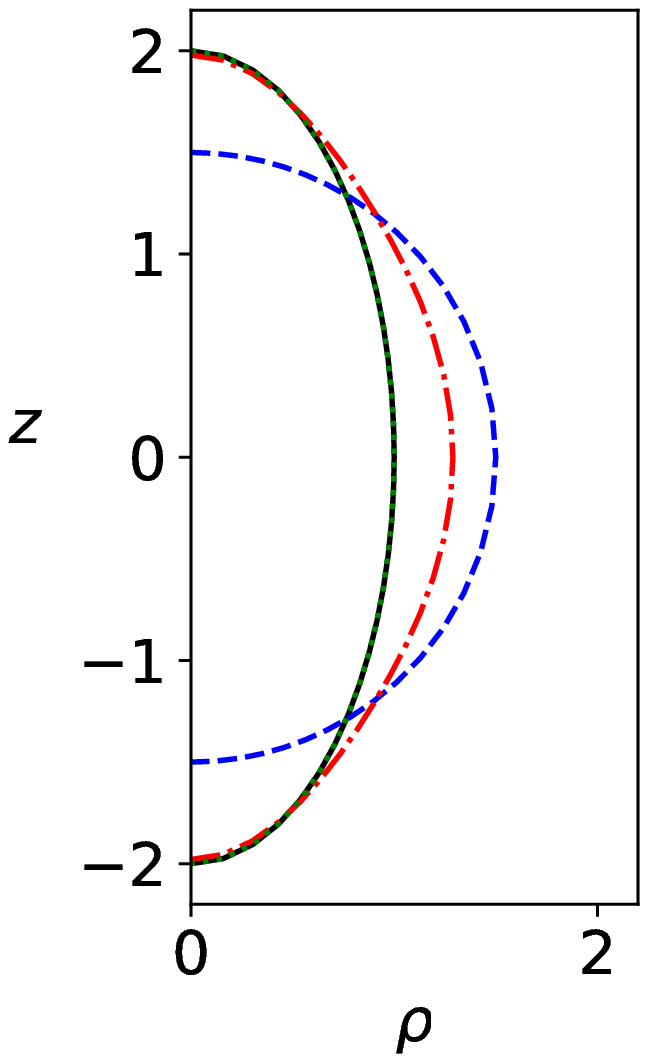}\\
  \includegraphics[width=0.49\textwidth]{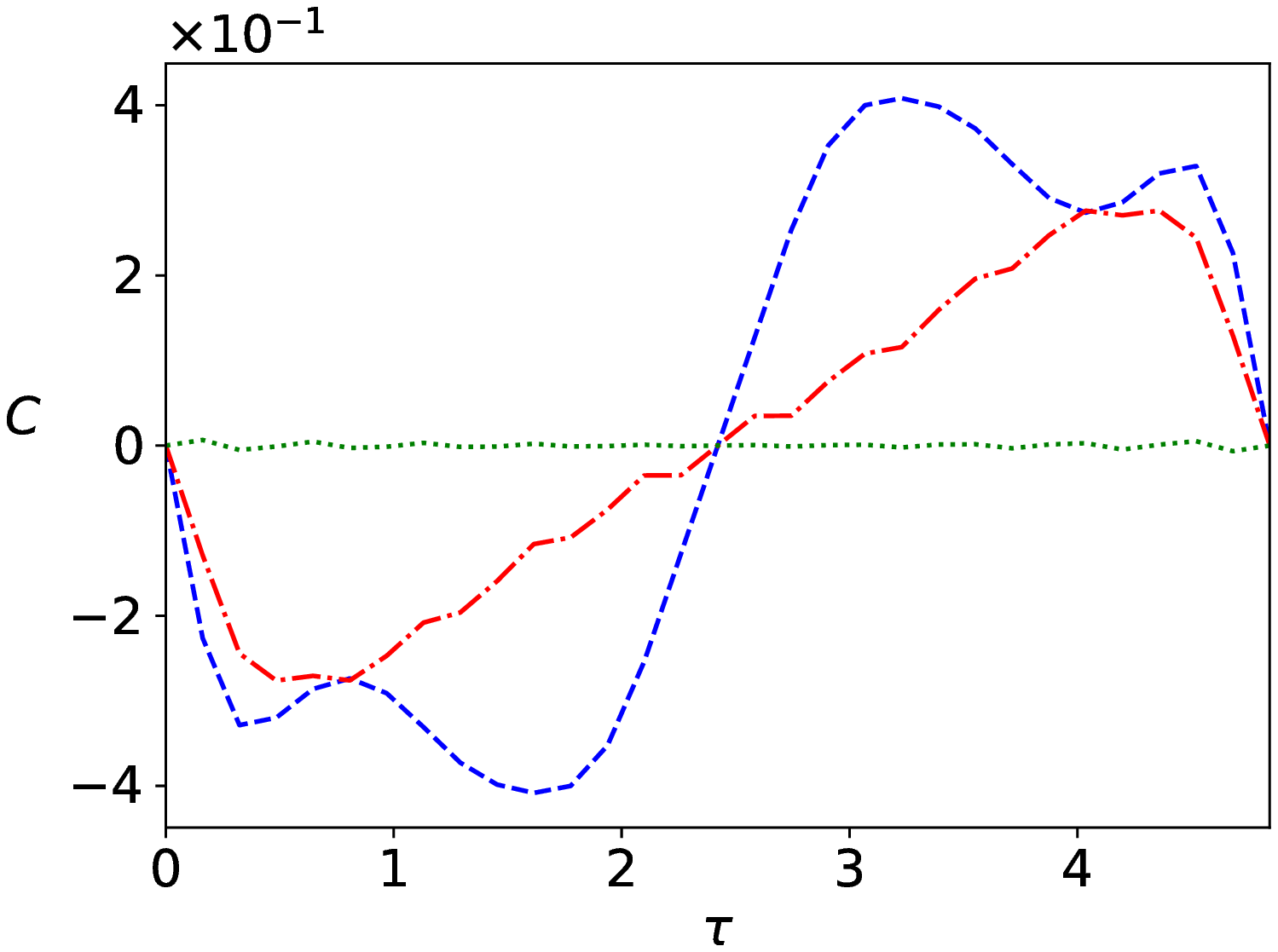}
  \includegraphics[width=0.49\textwidth]{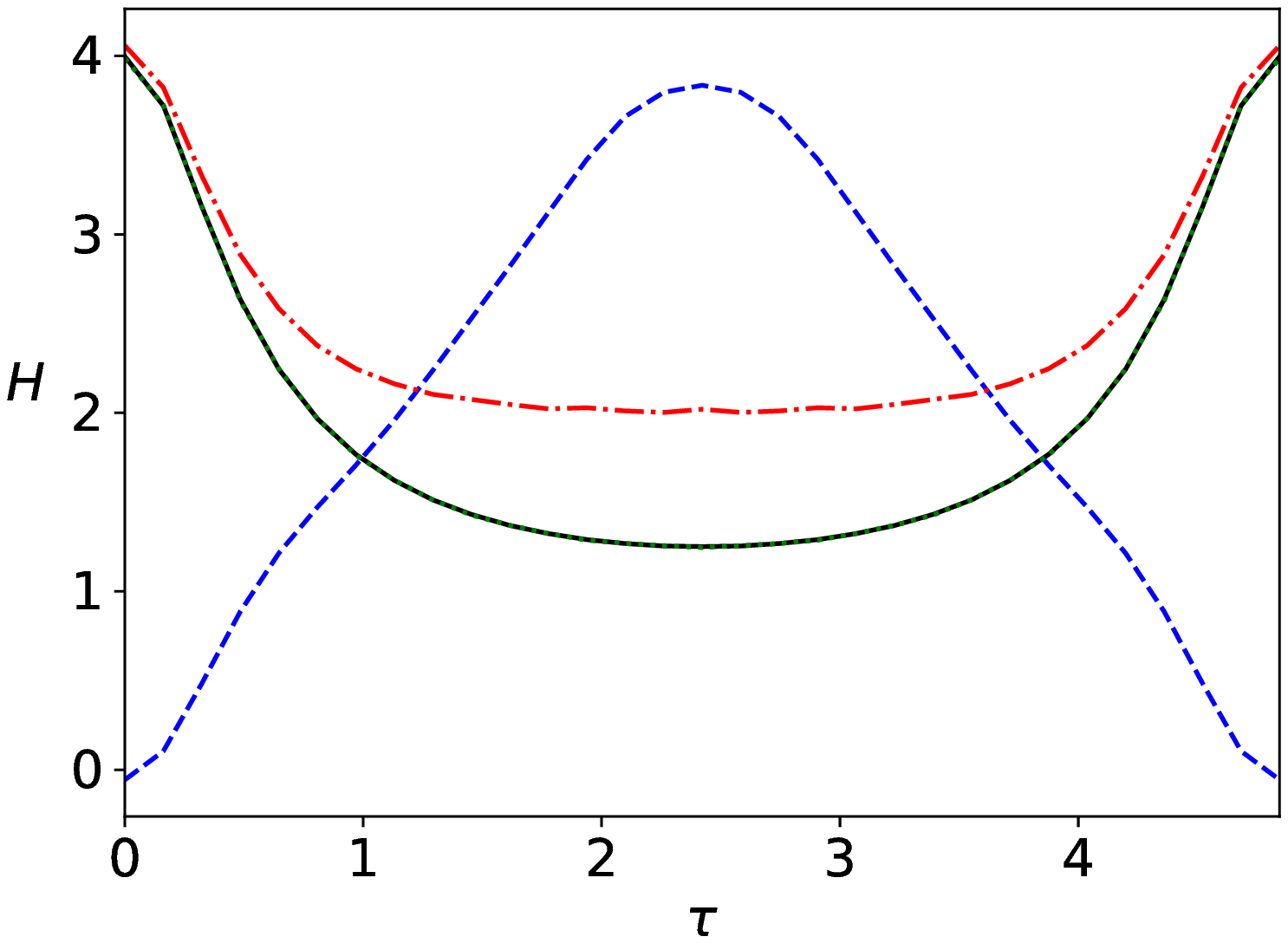}\\
  \includegraphics[width=0.49\textwidth]{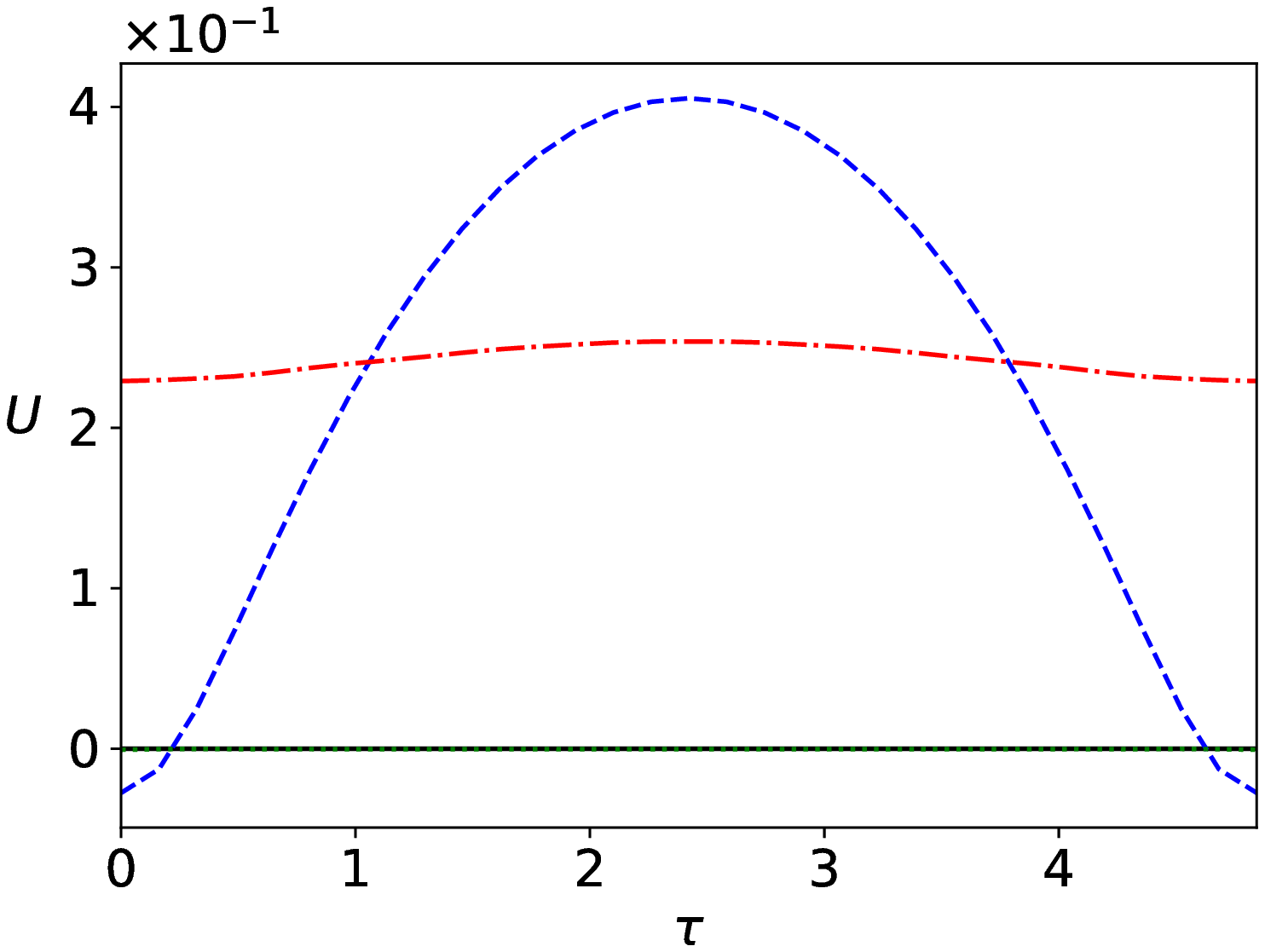}
  \includegraphics[width=0.49\textwidth]{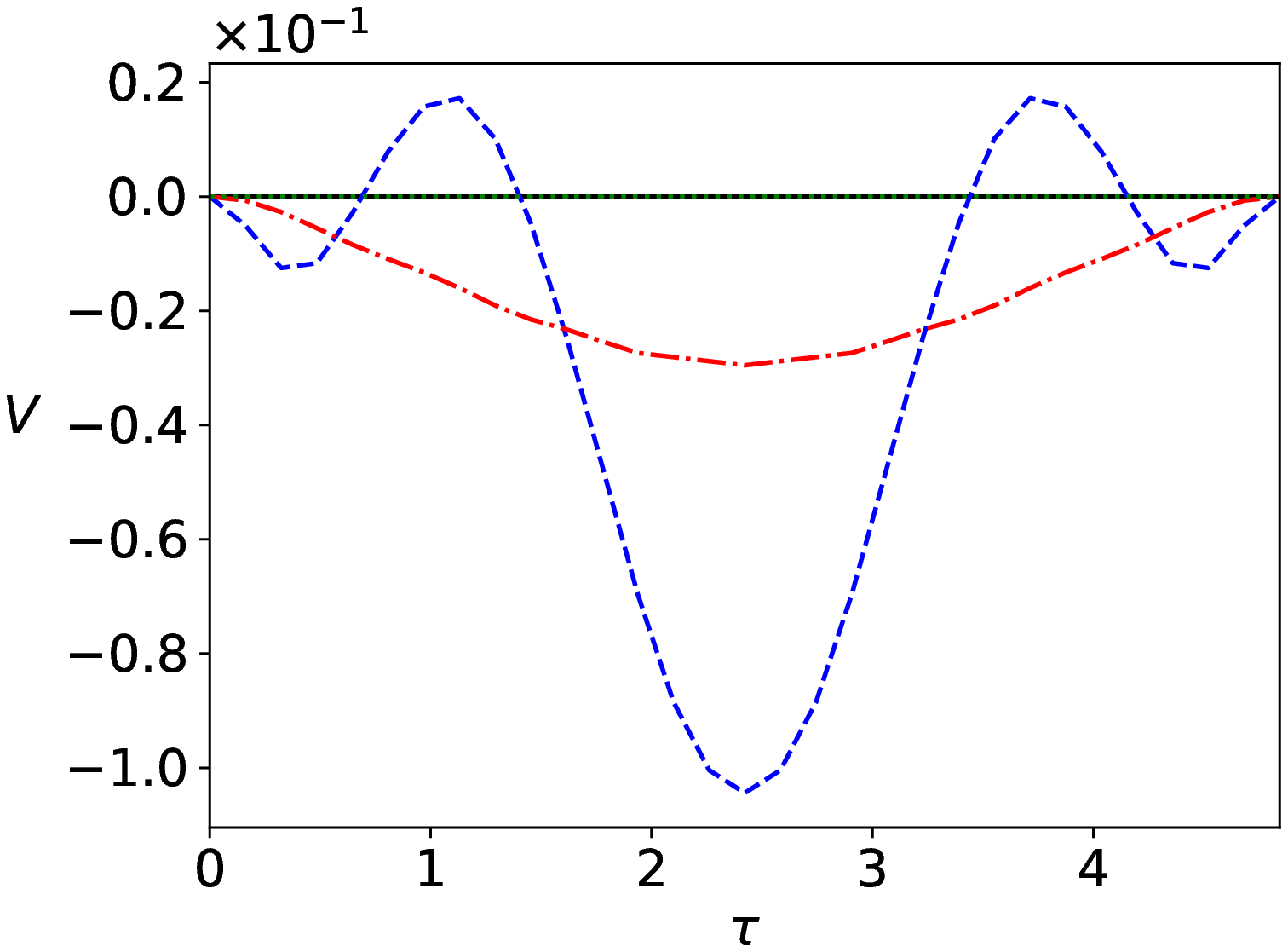}\\
  \includegraphics[width=0.49\textwidth]{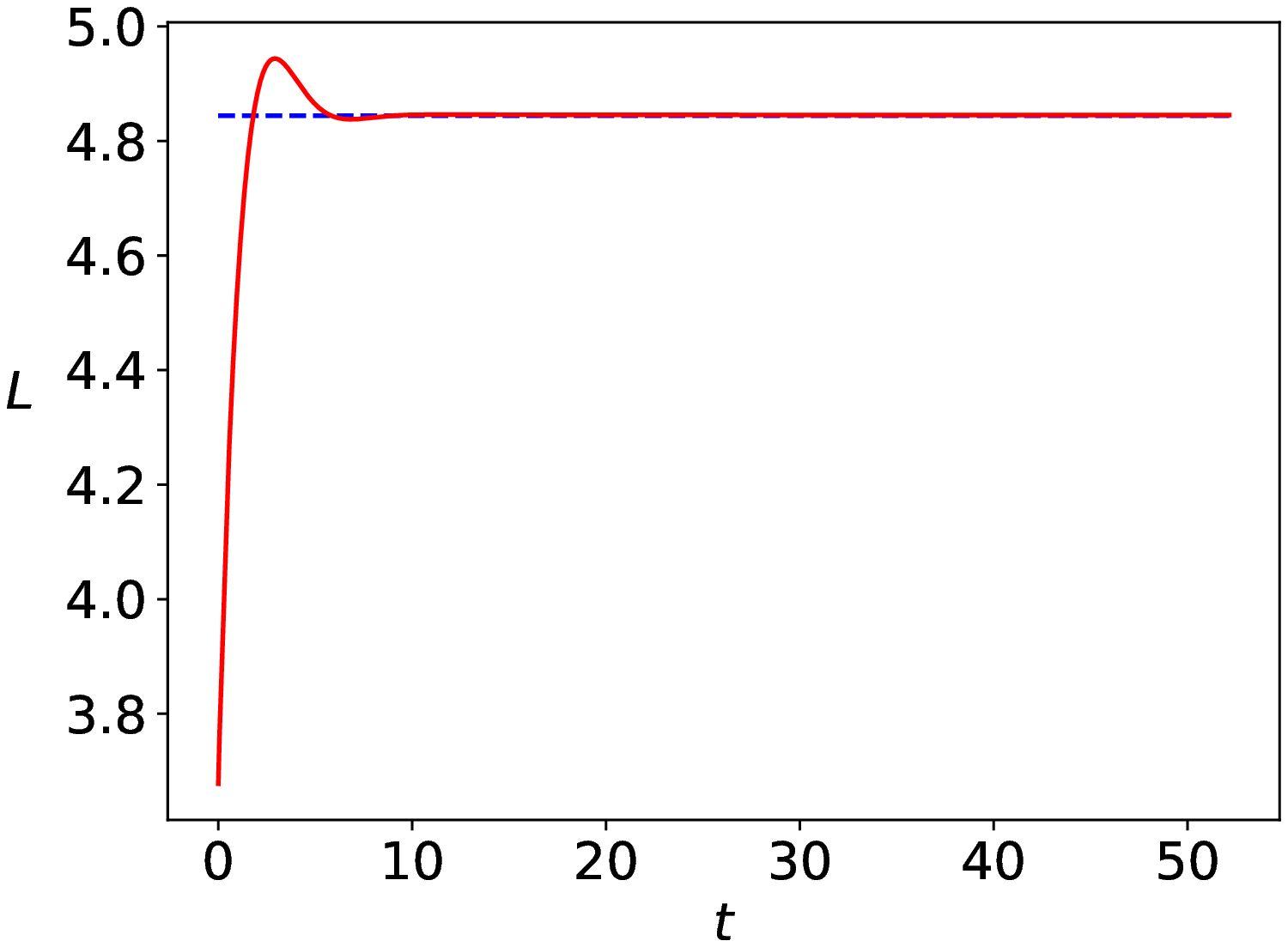}
  \includegraphics[width=0.49\textwidth]{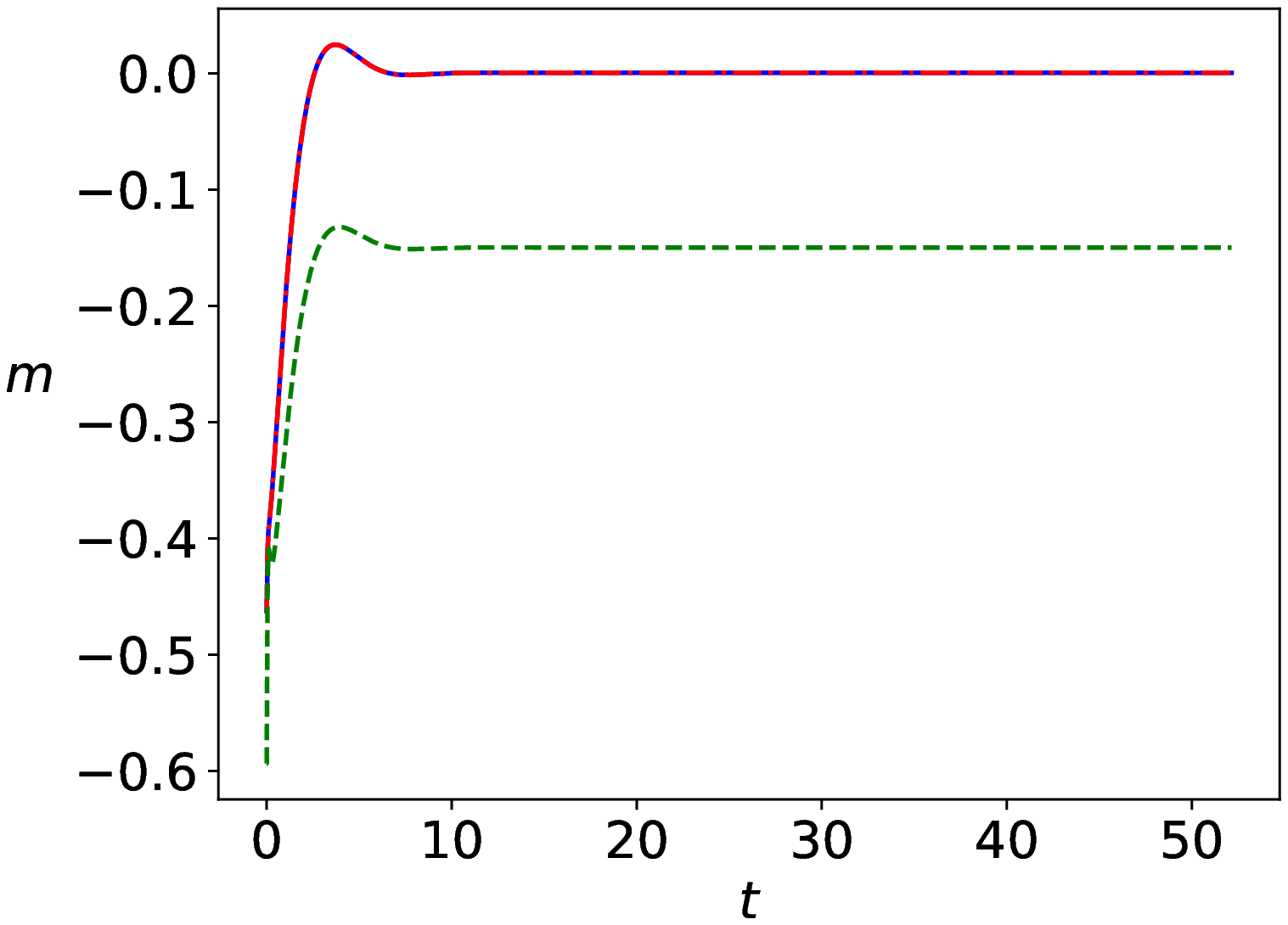}
  \caption{\small Flow with evolving metric. 
    The initial curve is taken to be a circle in Euclidean space ($r_0 = 1.5$), 
    the target curve an ellipse in Euclidean space (semi-major axes 
    $\overline{\rho} = 1, \overline{z} = 2$).
    In addition to the quantities described in Figure~\ref{f:static_Euclidean_circle_to_ellipse}, we also show the metric fields $U$ and $V$ along the flowing curve here.
    The different curves correspond to flow times $t=0$
    (dashed blue), $t=0.26$ (dash-dotted red) and $t=52.1$ (dotted
    green), with the target solution plotted in solid black.
    In the bottom right panel we plot the three masses described in
    the main text: ADM (solid blue), Hawking (dashed green) and
    pseudo-Newtonian (dash-dotted red, indistinguishable from the
    solid blue line).
  }
  \label{f:fly_Euclidean_circle_to_ellipse_exz}
\end{figure}

In Figure~\ref{f:fly_Euclidean_circle_to_ellipse_exz}, we choose as target curve
an ellipse in Euclidean space.
The initial curve is taken to be a circle in Euclidean space.
Already at the initial time, the flow departs from Euclidean space
because the initial (vanishing) $U$ and $V$ are replaced with the (non-trivial)
solution to the Weyl--Papapetrou equations with boundary data determined by
the prescribed function $\overline{\lambda}(\tau)$ 
according to \eqref{e:Ubc}, see also Section~\ref{s:Euclidean}.
This also implies that the initial curve is no longer parametrised by
arclength once $U$ and $V$ have been updated, hence $C\neq 0$ initially.
The flow does converge to the target curve and the Euclidean metric as desired.
This example shows clearly how the Hawking mass can differ from the
other masses; in this case it is negative, whereas the other masses are zero.
It should be noted that the expressions we use for the masses are only physically meaningful
as $t\to\infty$ as the curve becomes parametrised by arclength.

\begin{figure}[t]
  \centering
  \includegraphics[width=0.49\textwidth]{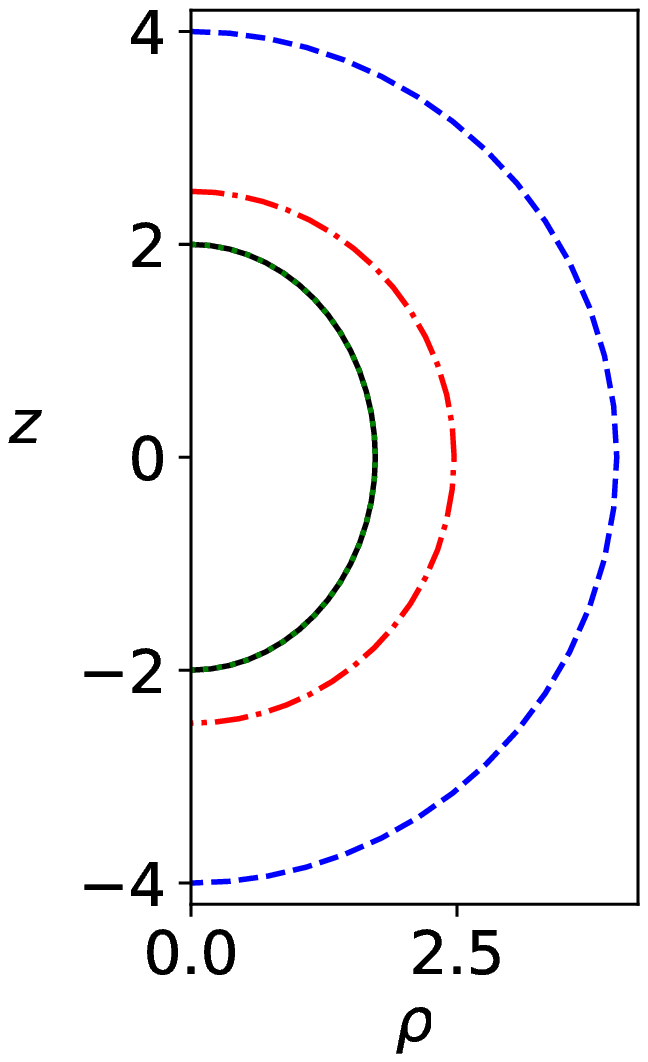}\\
  \includegraphics[width=0.49\textwidth]{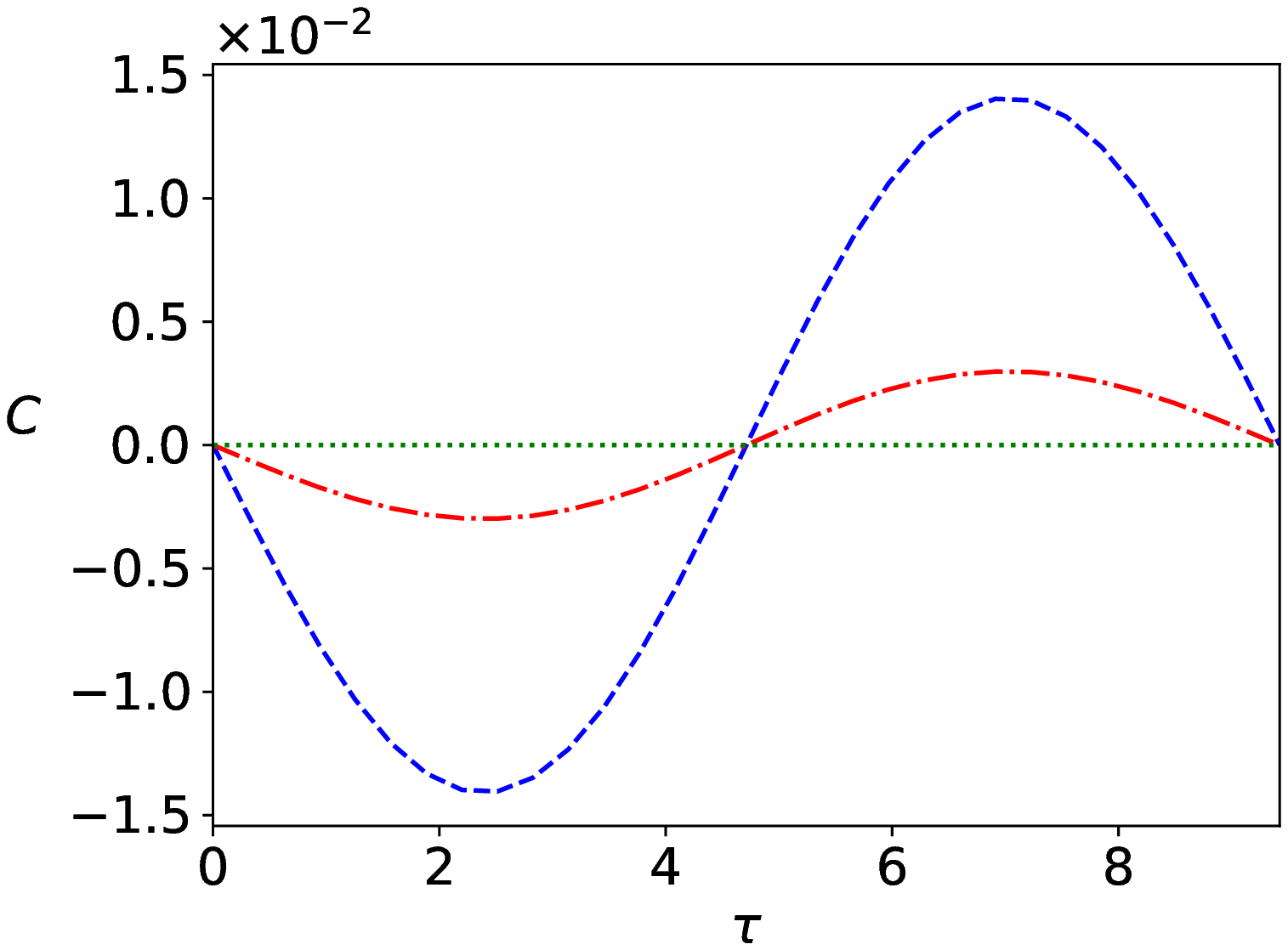}
  \includegraphics[width=0.49\textwidth]{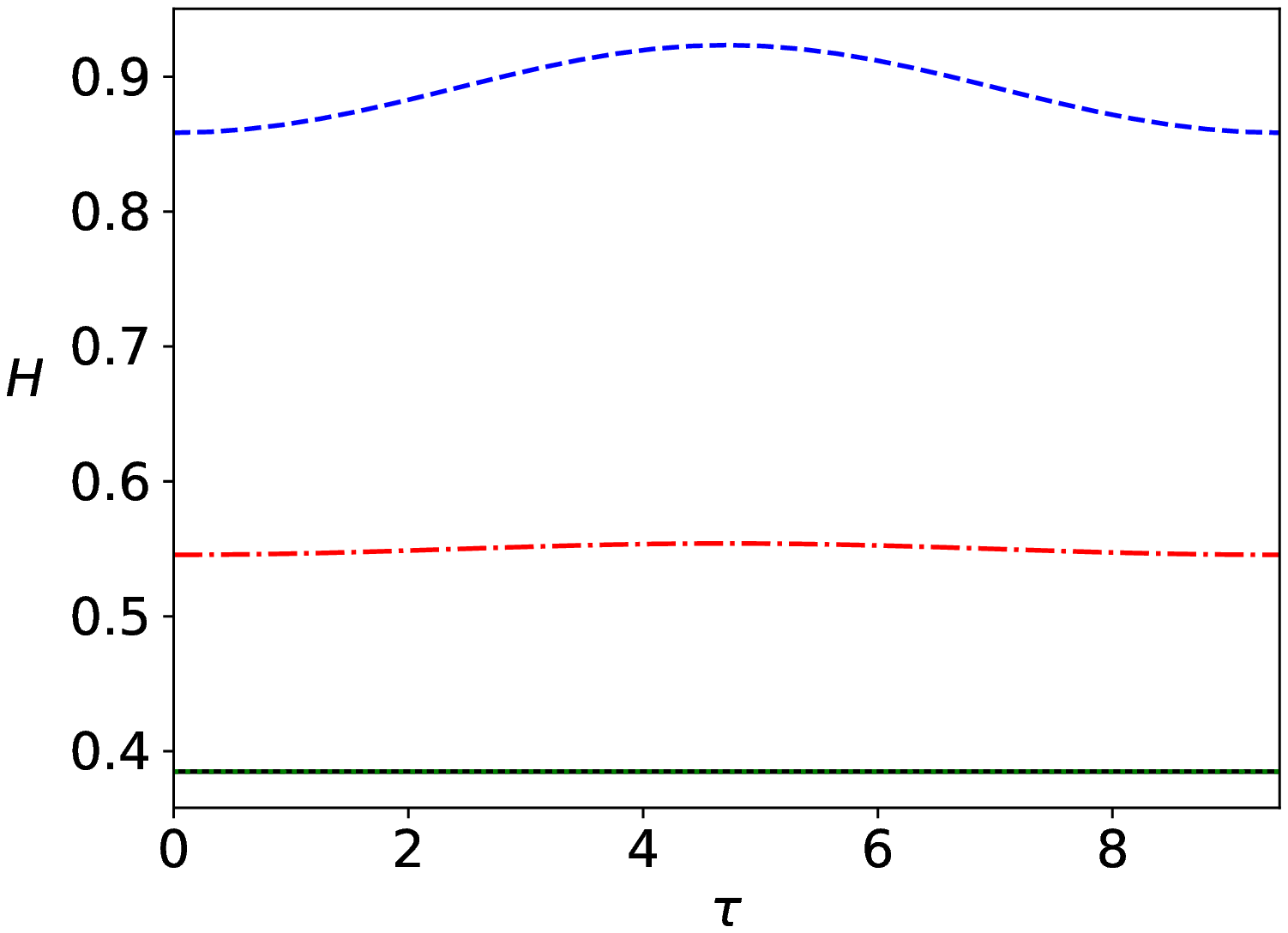}\\
  \includegraphics[width=0.49\textwidth]{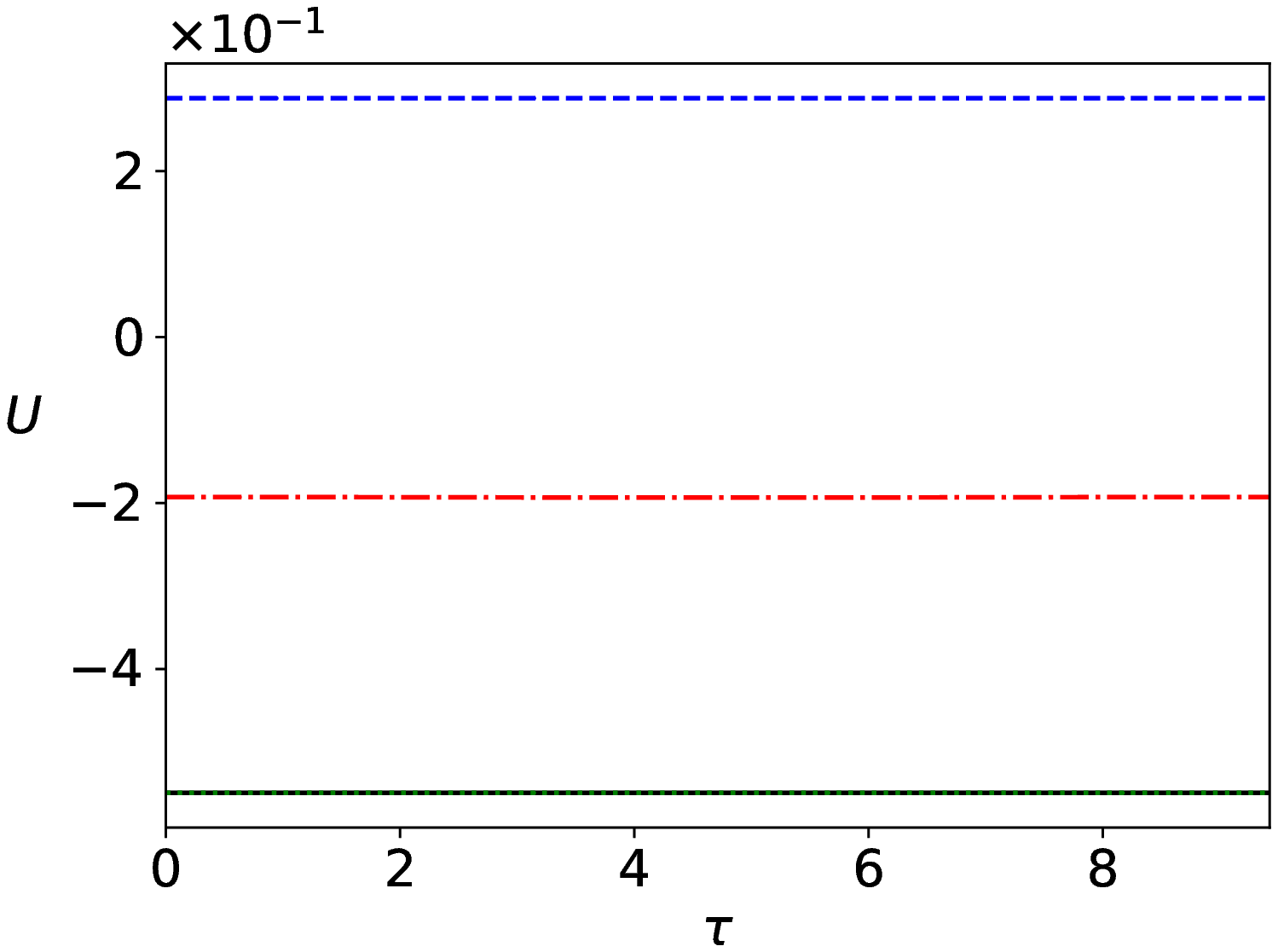}
  \includegraphics[width=0.49\textwidth]{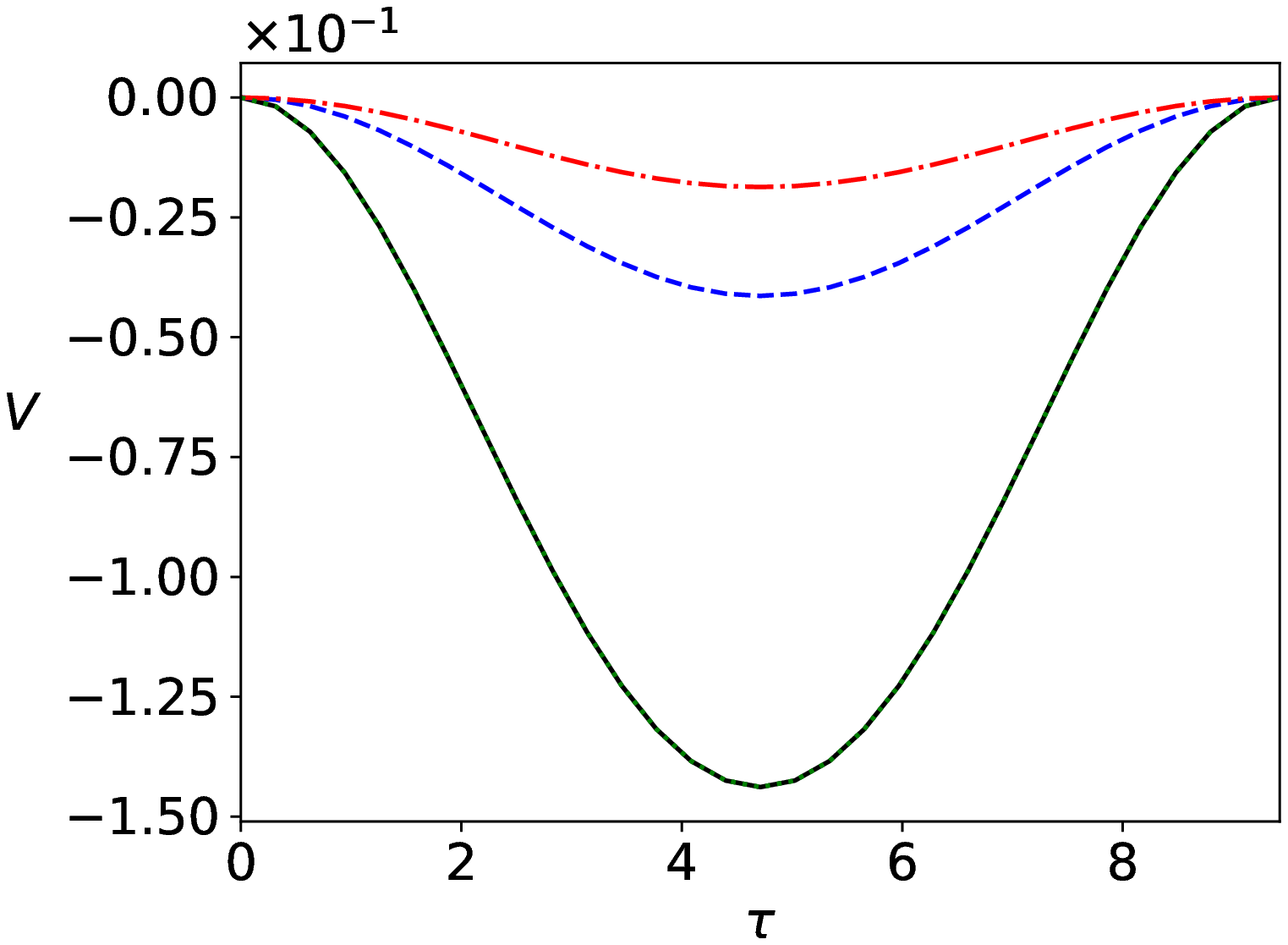}\\
  \includegraphics[width=0.49\textwidth]{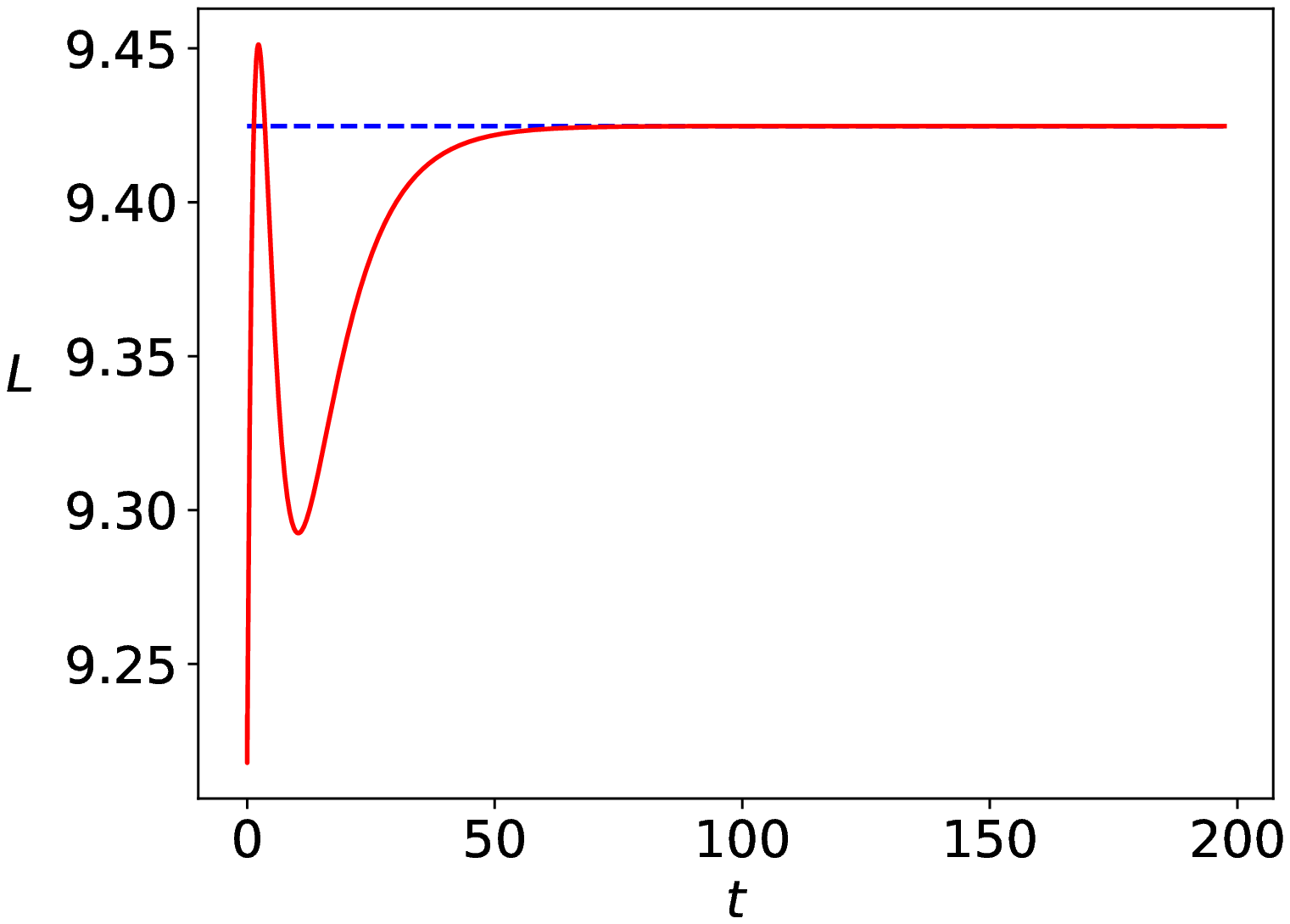}
  \includegraphics[width=0.49\textwidth]{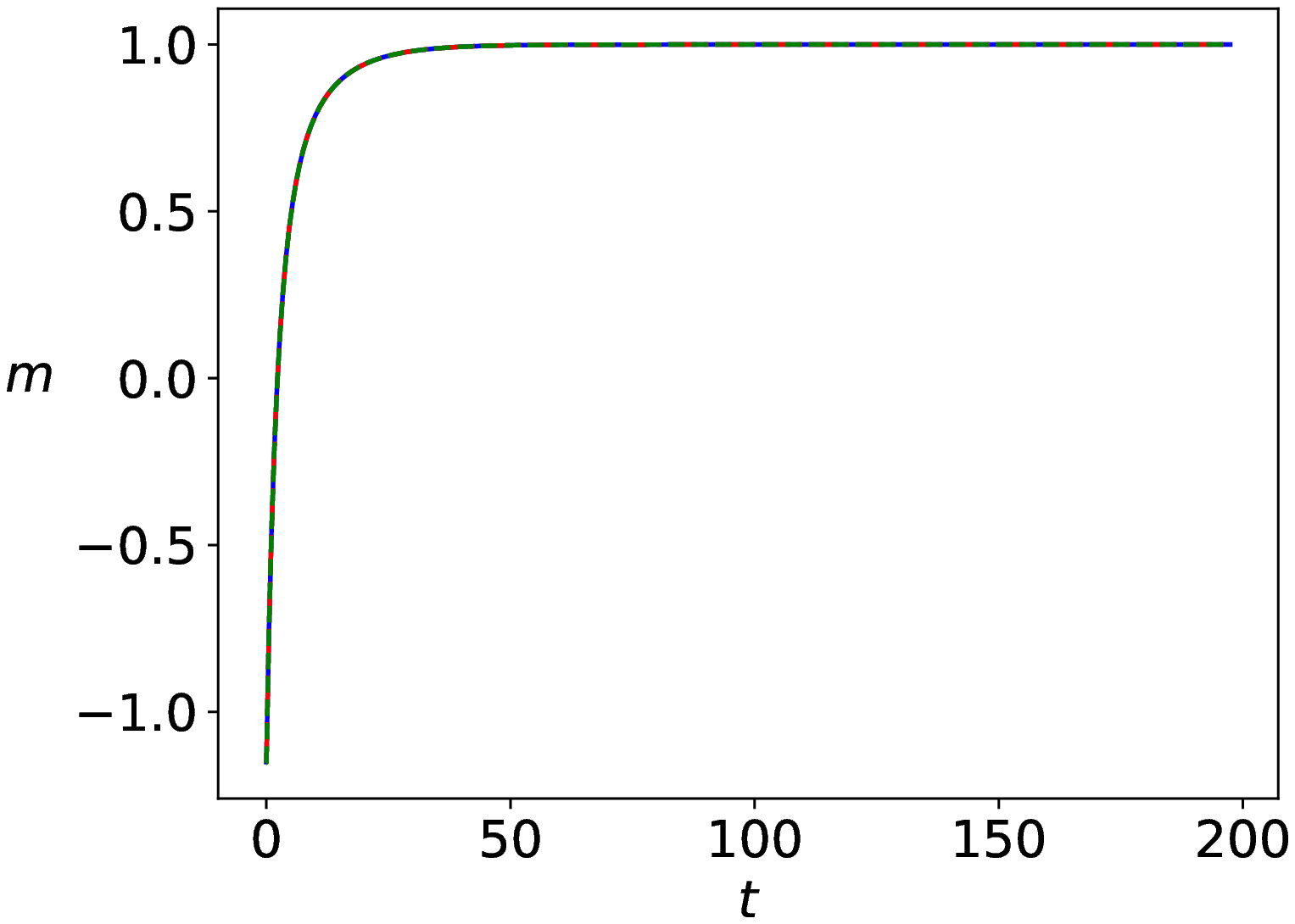}
  \caption{\small Flow with evolving metric. 
    The initial curve is taken to be a circle in Euclidean space ($r_0 = 4$), 
    the target curve a circle in Schwarzschild coordinates
    ($\overline{r}_{S} = 3$) in Schwarzschild space of mass $M=1$.
    The same quantities as in Figure~\ref{f:fly_Euclidean_circle_to_ellipse_exz} are
    plotted.
    In the first five panels, the different curves correspond to 
    flow times $t=0$ (dashed blue), $t=4.9$ (dash-dotted red) and 
    $t=197.4$ (dotted green), with the target solution plotted in
    solid black. }
  \label{f:fly_ss_0_to_1}
\end{figure}

\begin{figure}[t]
  \centering
  \includegraphics[width=0.49\textwidth]{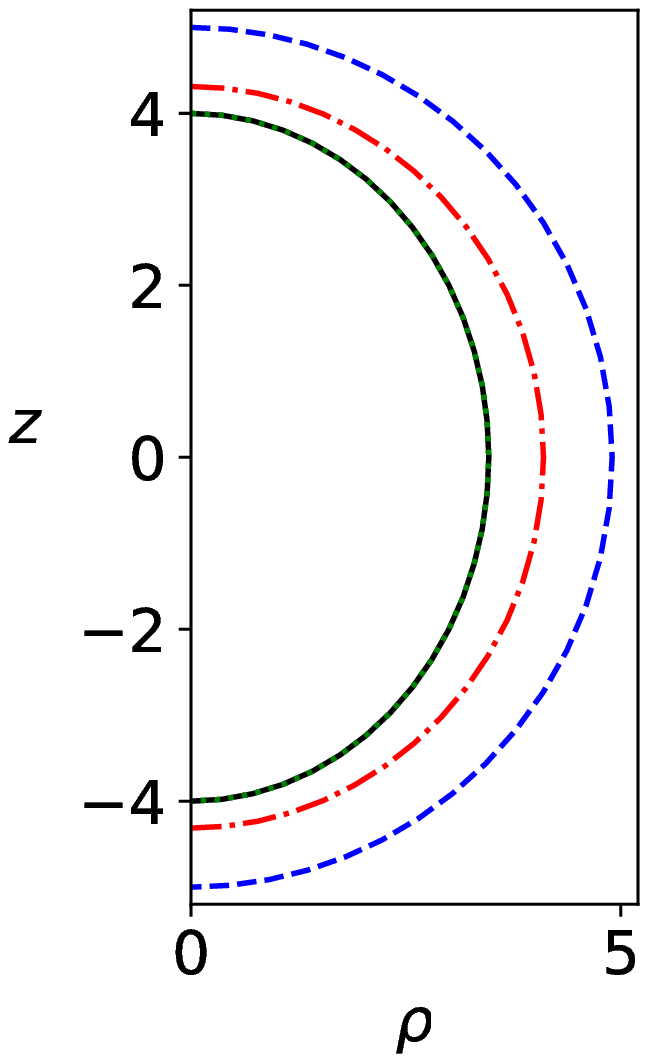}\\
  \includegraphics[width=0.49\textwidth]{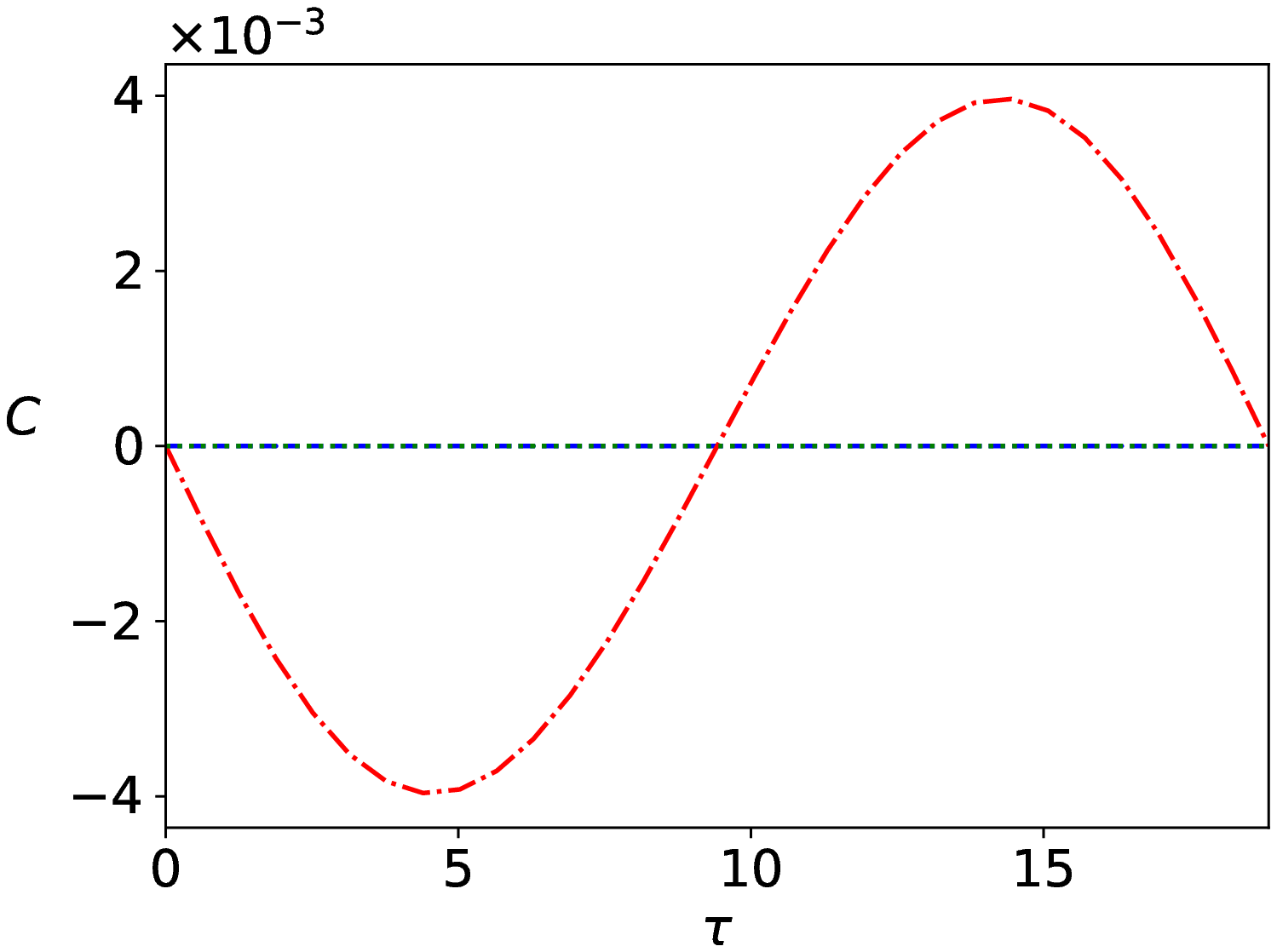}
  \includegraphics[width=0.49\textwidth]{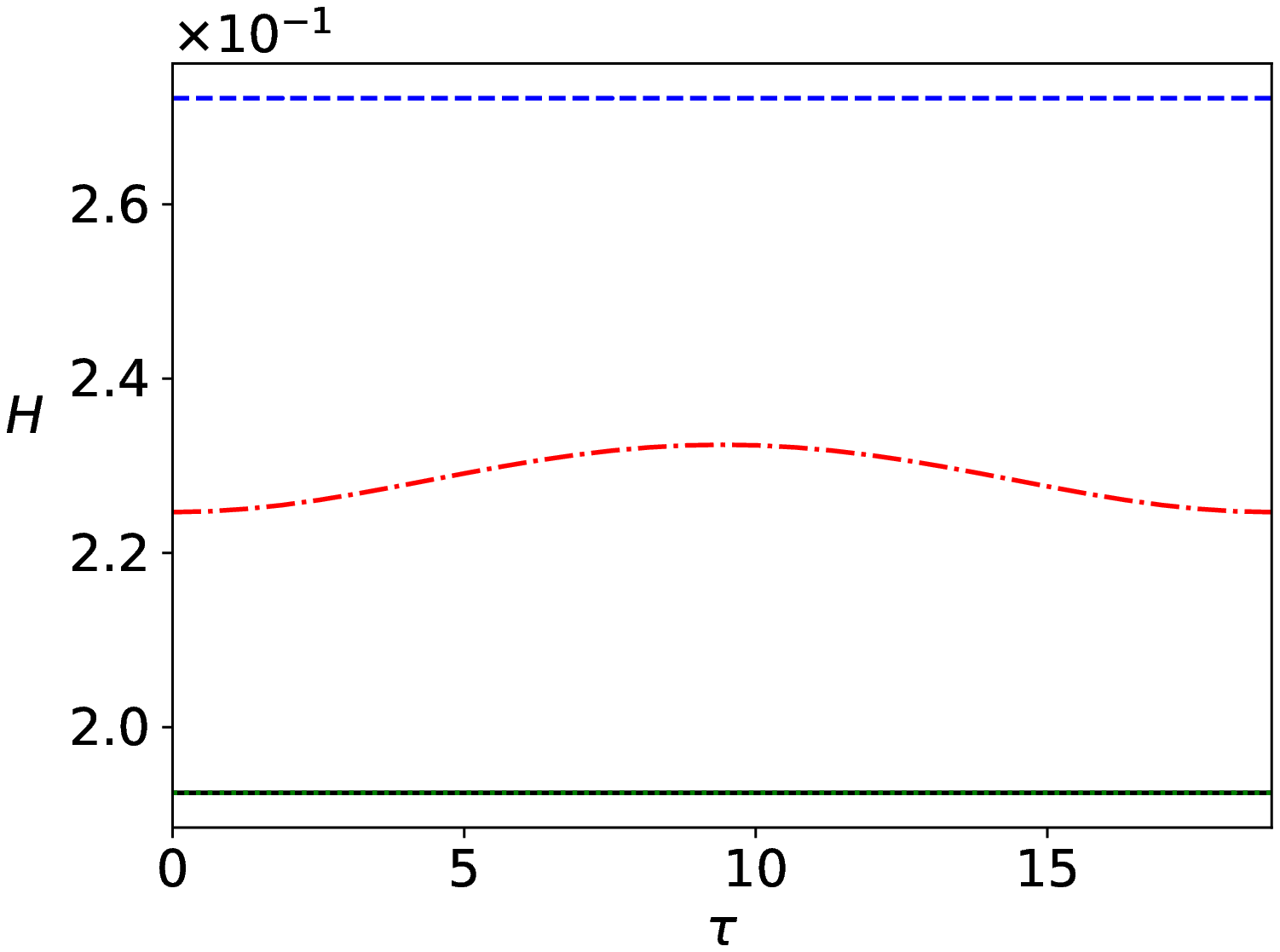}\\
  \includegraphics[width=0.49\textwidth]{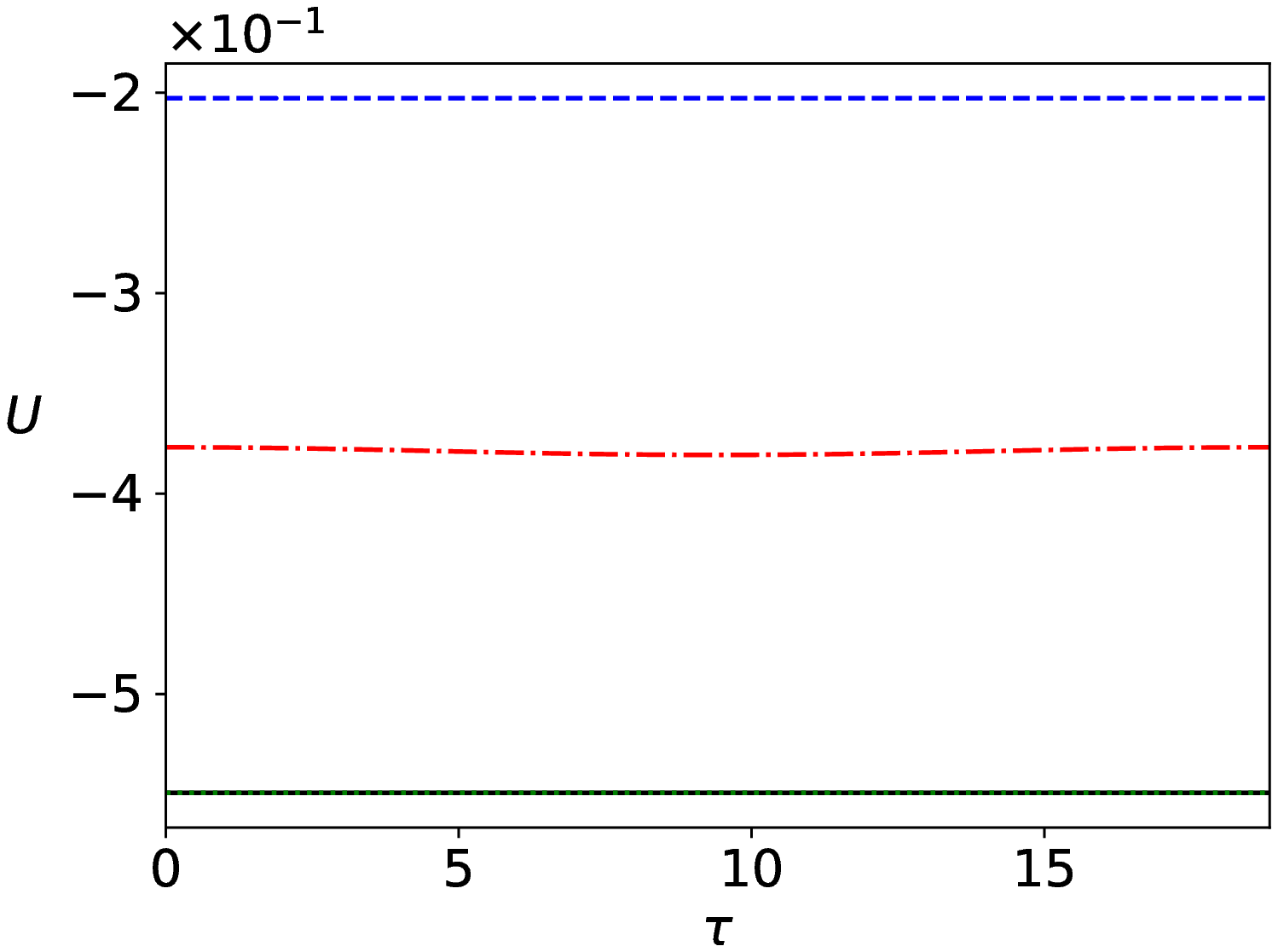}
  \includegraphics[width=0.49\textwidth]{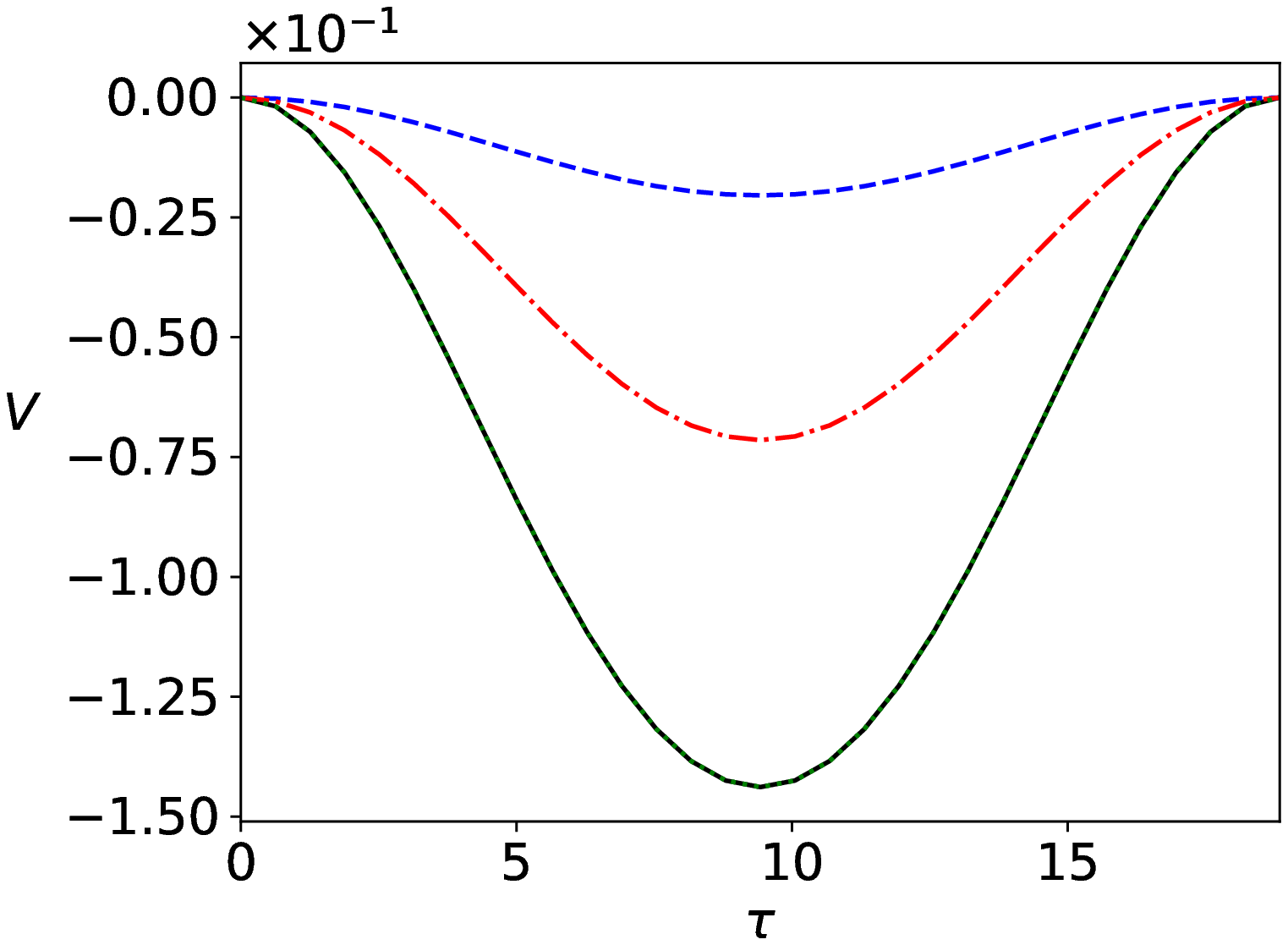}\\
  \includegraphics[width=0.49\textwidth]{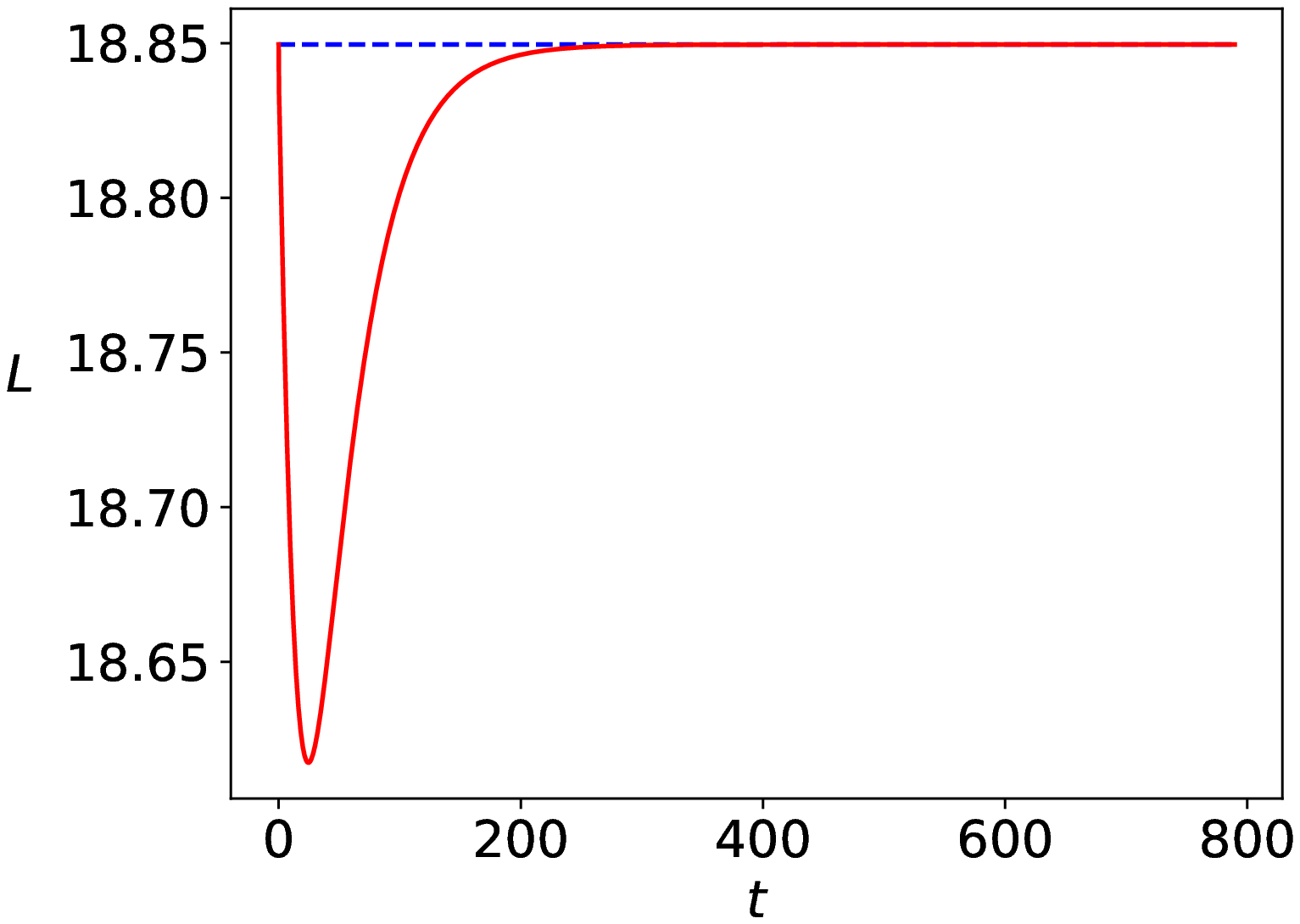}
  \includegraphics[width=0.49\textwidth]{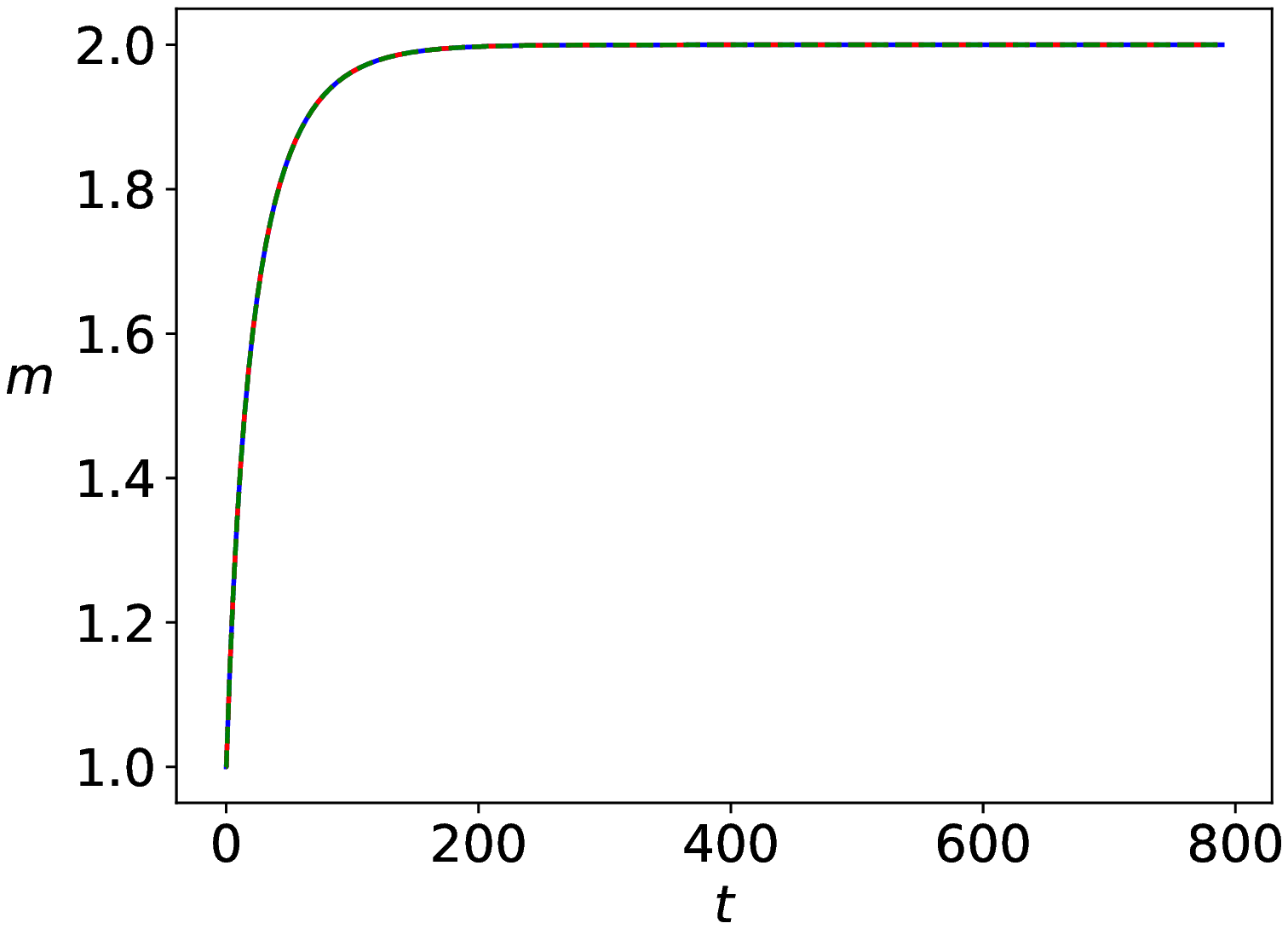}
  \caption{\small Flow with evolving metric.
    The initial curve is taken to be a circle ($r_{S,0} = 6$) 
    in $M=1$ Schwarzschild space, 
    the target curve a circle ($\overline{r}_{S} = 6$) in $M=2$
    Schwarzschild space.
    The same quantities as in Figure~\ref{f:fly_Euclidean_circle_to_ellipse_exz} are
    plotted.
    In the first five panels, the different curves correspond to 
    flow times $t=0$ (dashed blue), $t=19.7$ (dash-dotted red) and 
    $t=789.6$ (dotted green), with the target solution plotted in solid black.}
  \label{f:fly_ss_1_to_2}
\end{figure}

\begin{figure}[t]
  \centering
  \includegraphics[width=0.49\textwidth]{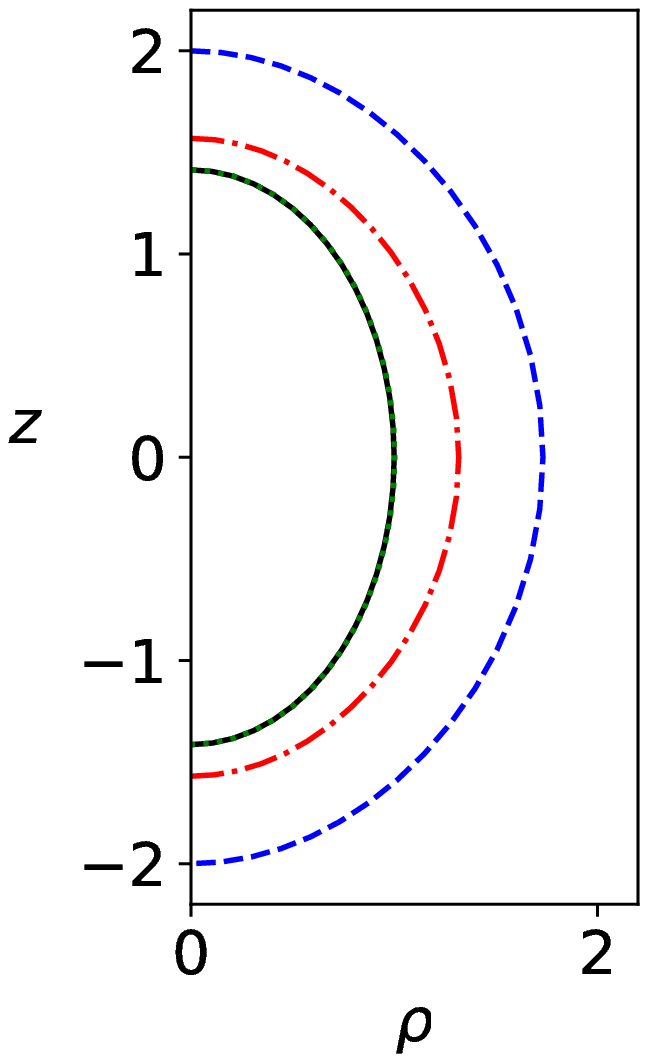}\\
  \includegraphics[width=0.49\textwidth]{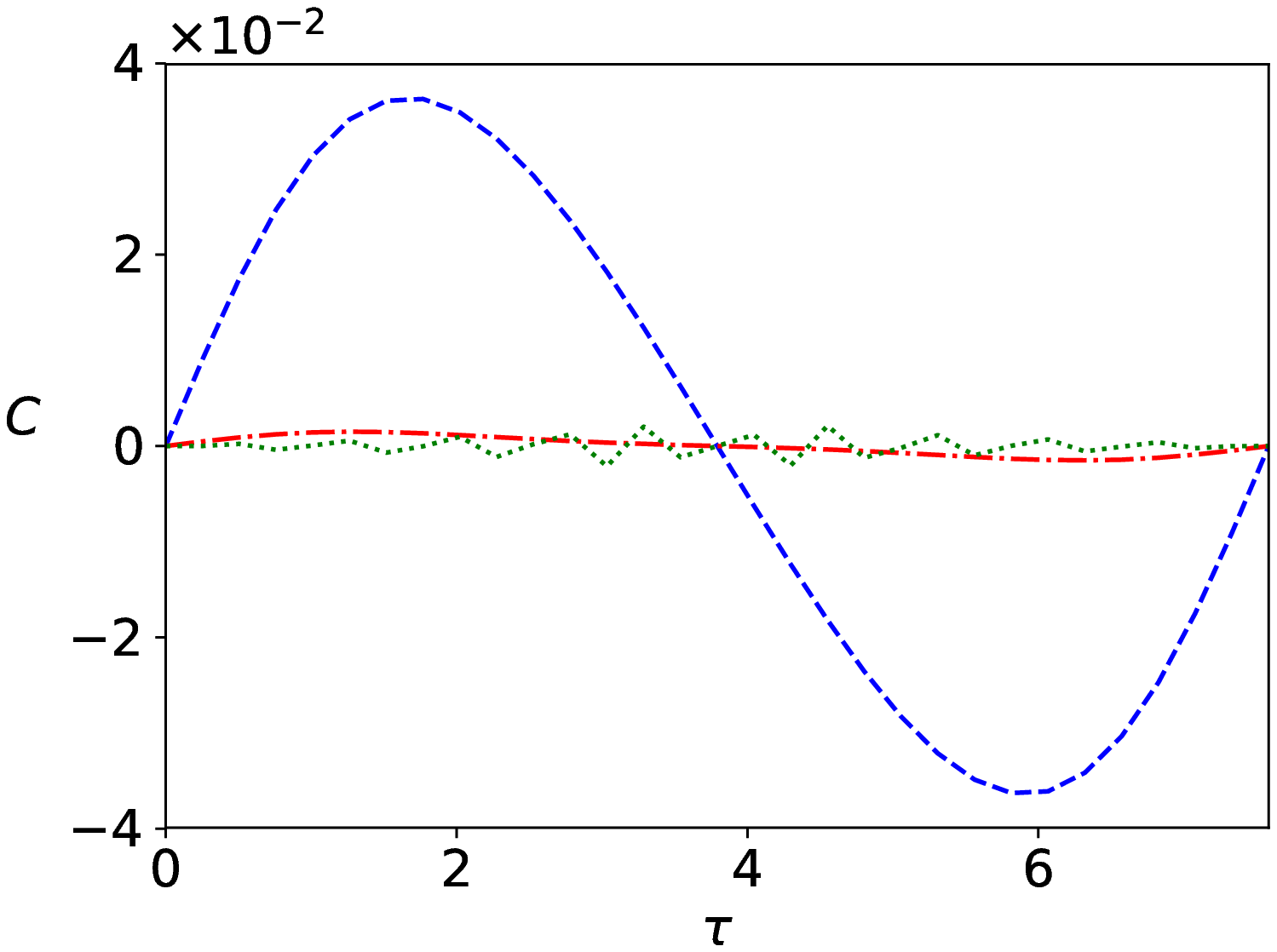}
  \includegraphics[width=0.49\textwidth]{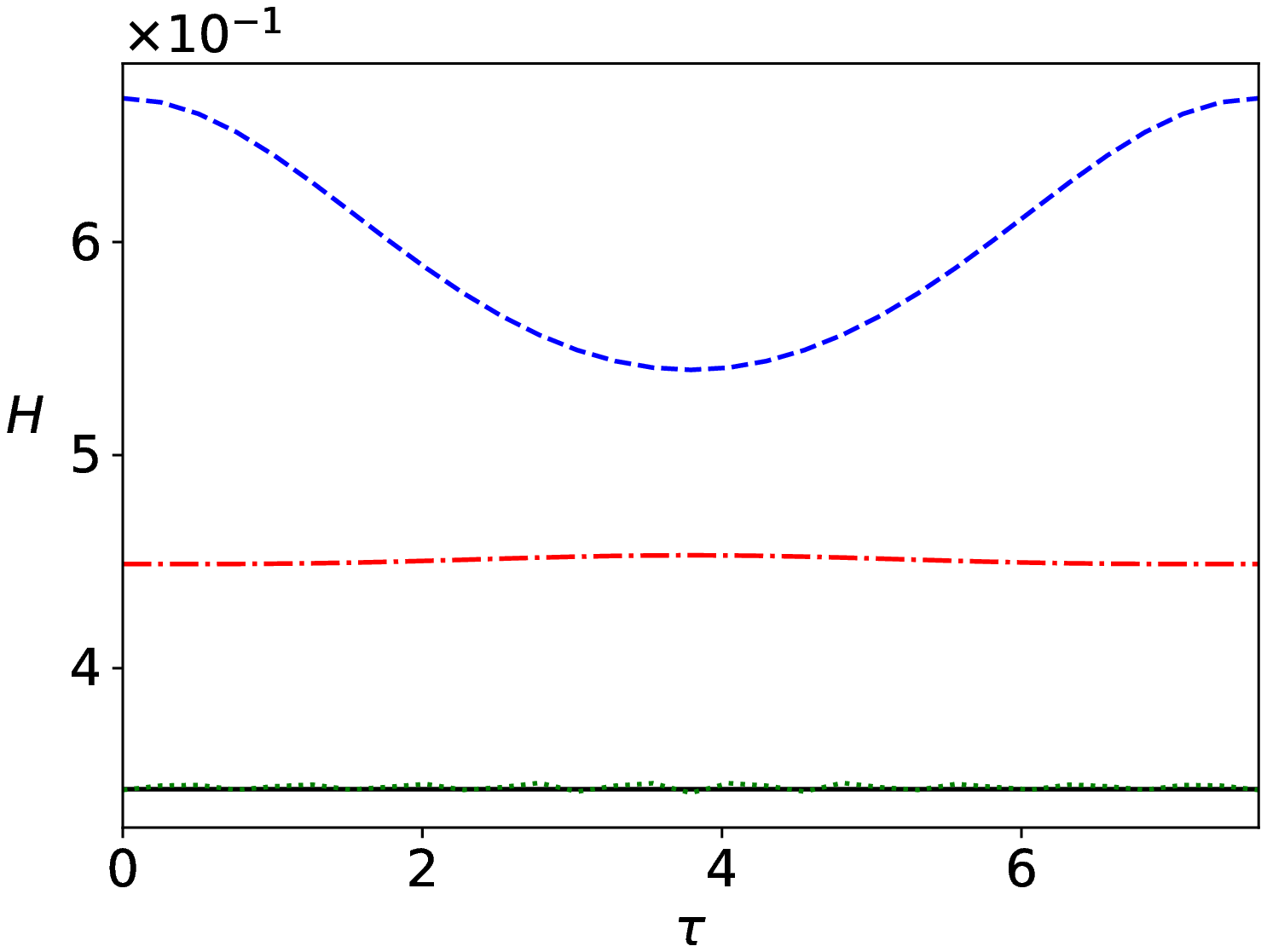}\\
  \includegraphics[width=0.49\textwidth]{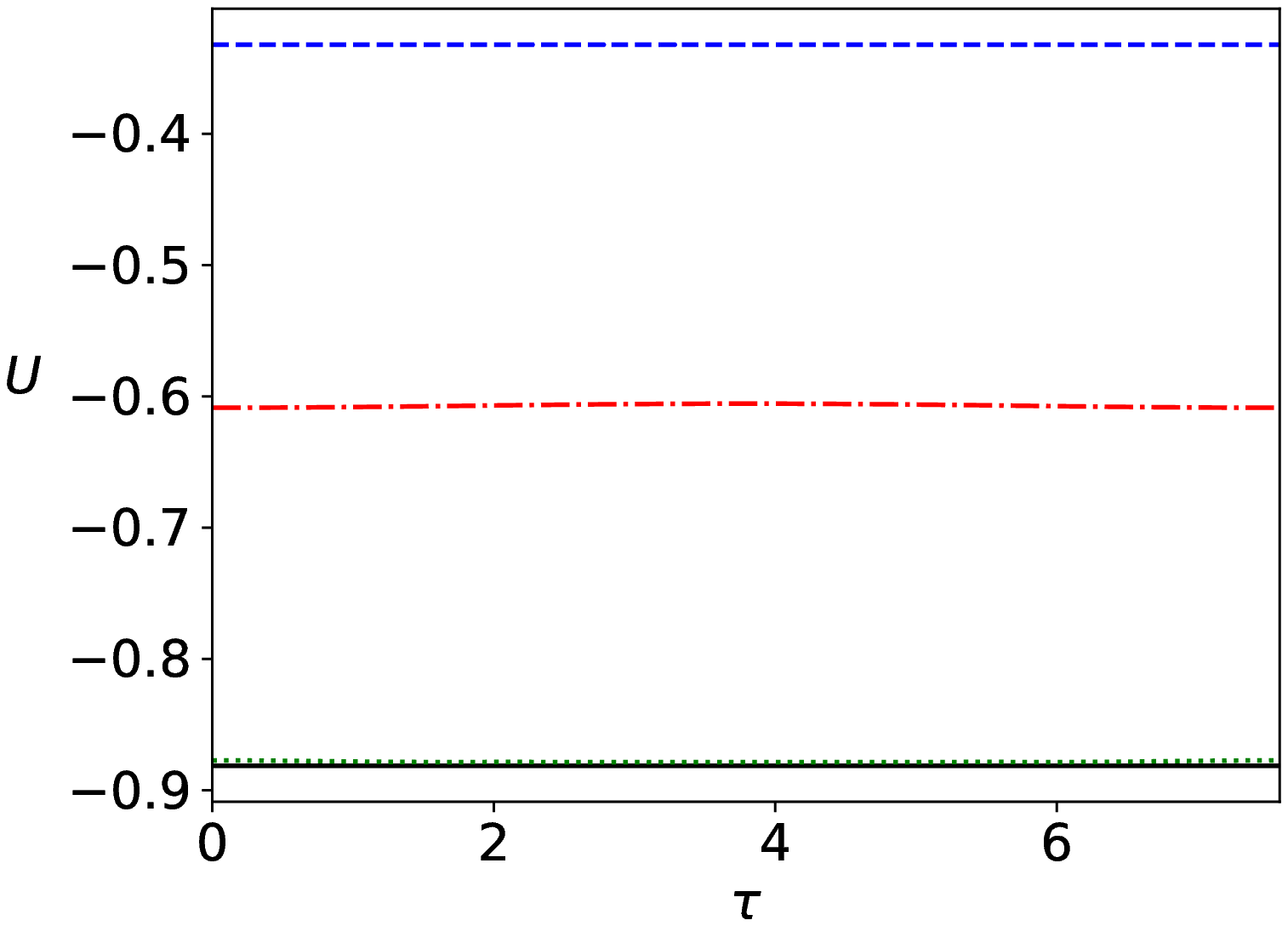}
  \includegraphics[width=0.49\textwidth]{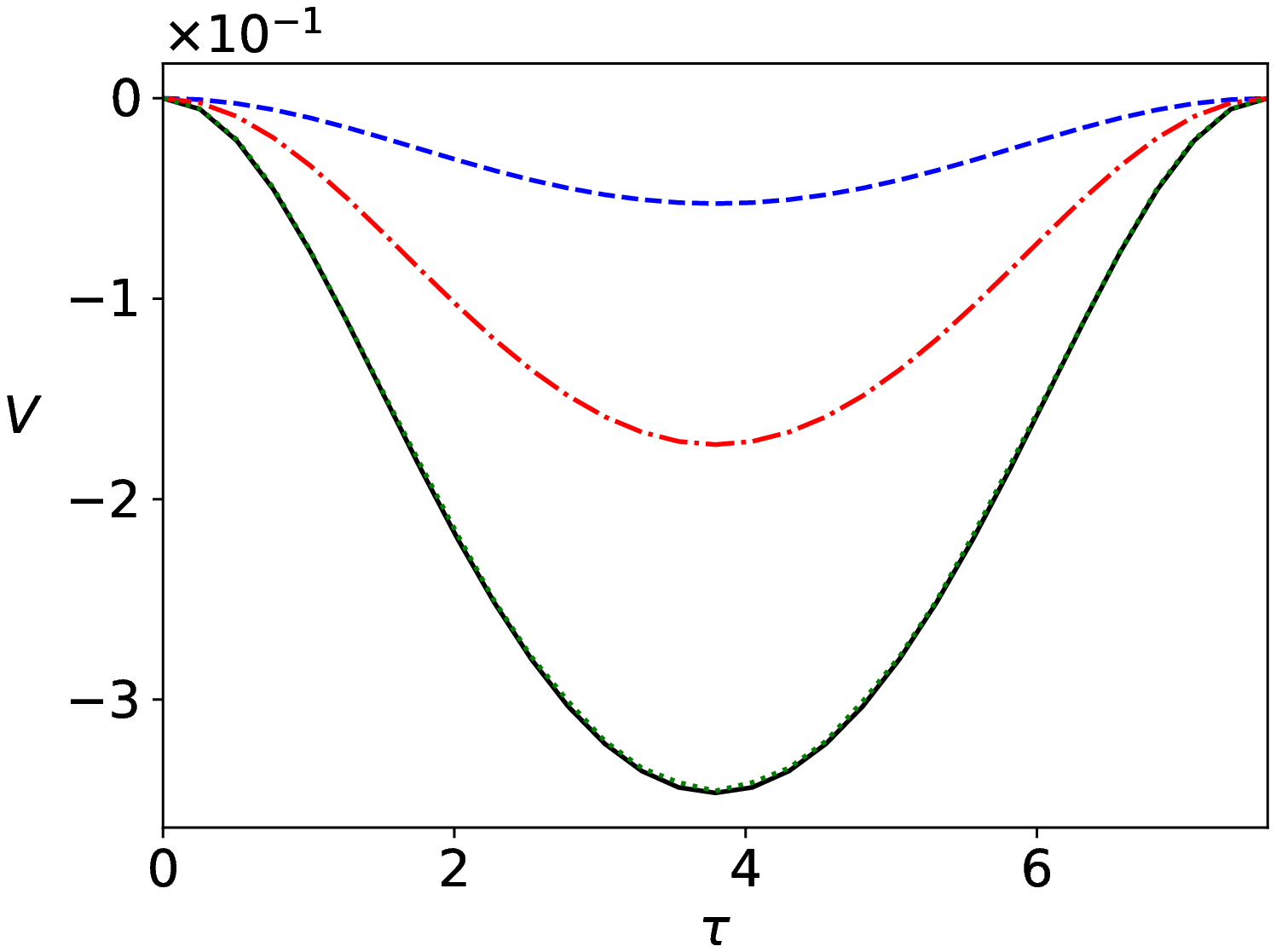}\\
  \includegraphics[width=0.49\textwidth]{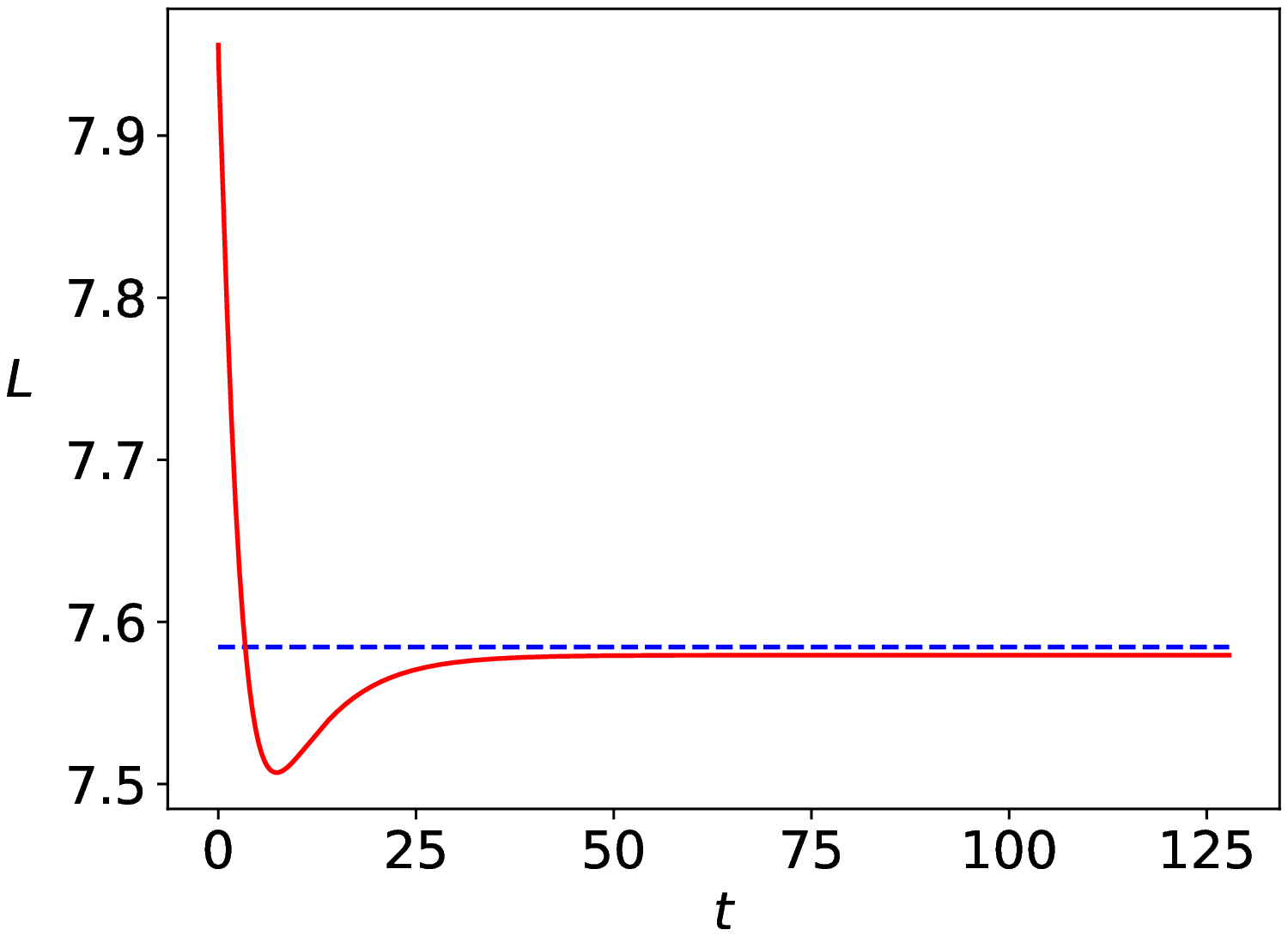}
  \includegraphics[width=0.49\textwidth]{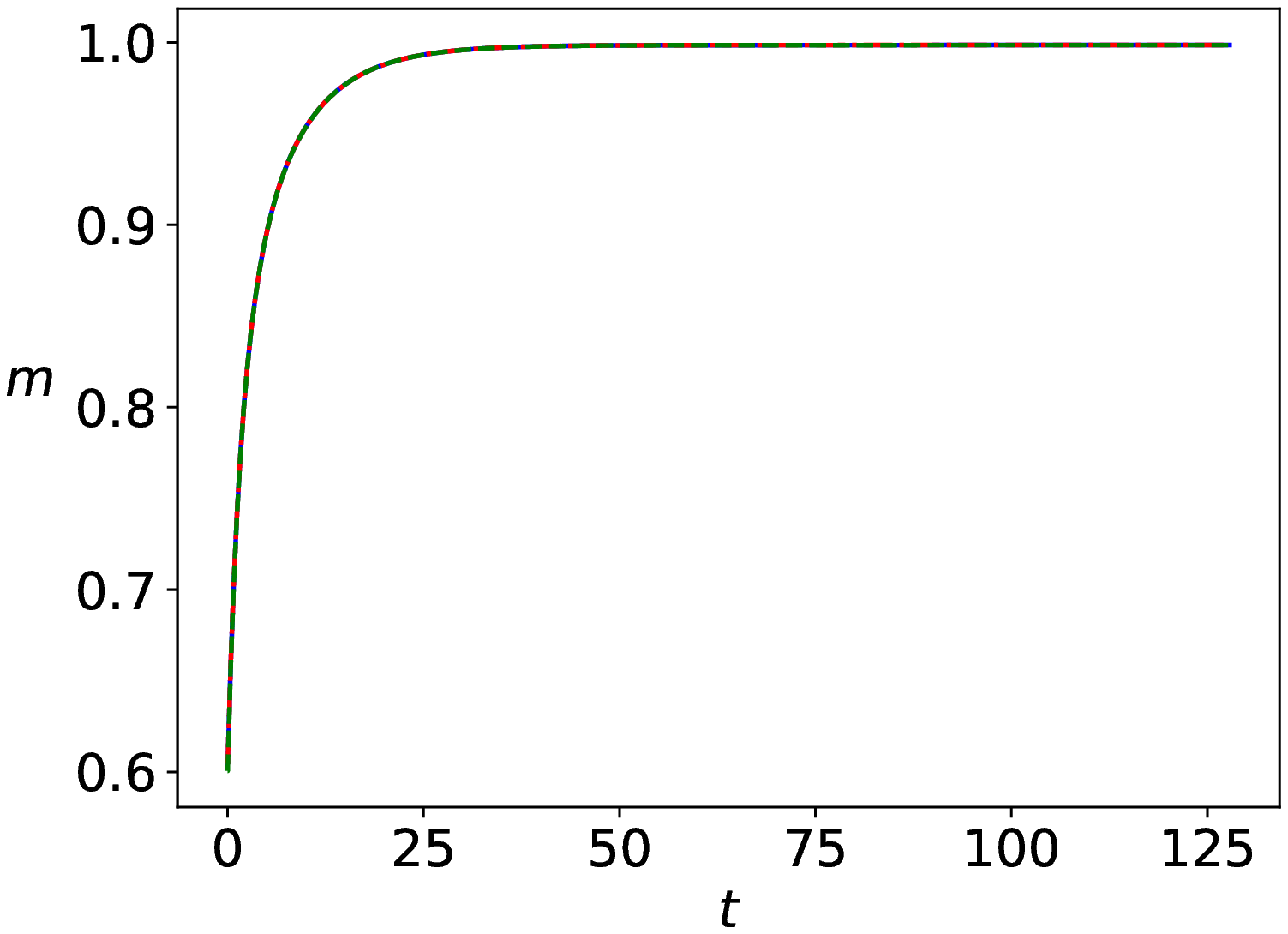}
  \caption{\small Flow with evolving metric.
    The initial curve is taken to be a (Schwarzschild coordinate) circle ($r_{S,0} = 3$) 
    in $M=1$ Schwarzschild space, 
    the target curve a (Schwarzschild coordinate) circle close to the horizon 
    ($\overline{r}_{S} = 2.41$) in Schwarzschild space of the same mass $M=1$.
    The same quantities as in Figure~\ref{f:fly_Euclidean_circle_to_ellipse_exz} are
    plotted.
    In the first five panels, the different curves correspond to 
    flow times $t=0$ (dashed blue), $t=3.2$ (dash-dotted red) and 
    $t=127.8$ (dotted green), with the target solution plotted in solid black.
  }
  \label{f:fly_ss_extremal}
\end{figure}

In Figure~\ref{f:fly_ss_0_to_1}, we present a flow from a circle in Euclidean
space to a (Schwarzschild coordinate) circle in Schwarzschild space ($M=1$).
Figure~\ref{f:fly_ss_1_to_2} shows a similar evolution of a
(Schwarzschild coordinate) circle in Schwarzschild space of mass $M=1$ to a
different Schwarzschild space of mass $M=2$.
In Figure~\ref{f:fly_ss_extremal}, we let a (Schwarzschild coordinate) circle in $M=1$
Schwarzschild space flow to a circle close to the horizon ($\overline{r}_S =
2.41$) in the same Schwarzschild system (but again notice that the coupled flow
departs from this metric at intermediate times).
This is the closest we could get to the horizon at $r_S = 2$ due to
the numerical problems associated with our method of solving the
Laplace equation for $U$ described in Section~\ref{s:nummethod.einstein}.

\begin{figure}[t]
  \centering
  \includegraphics[width=0.49\textwidth]{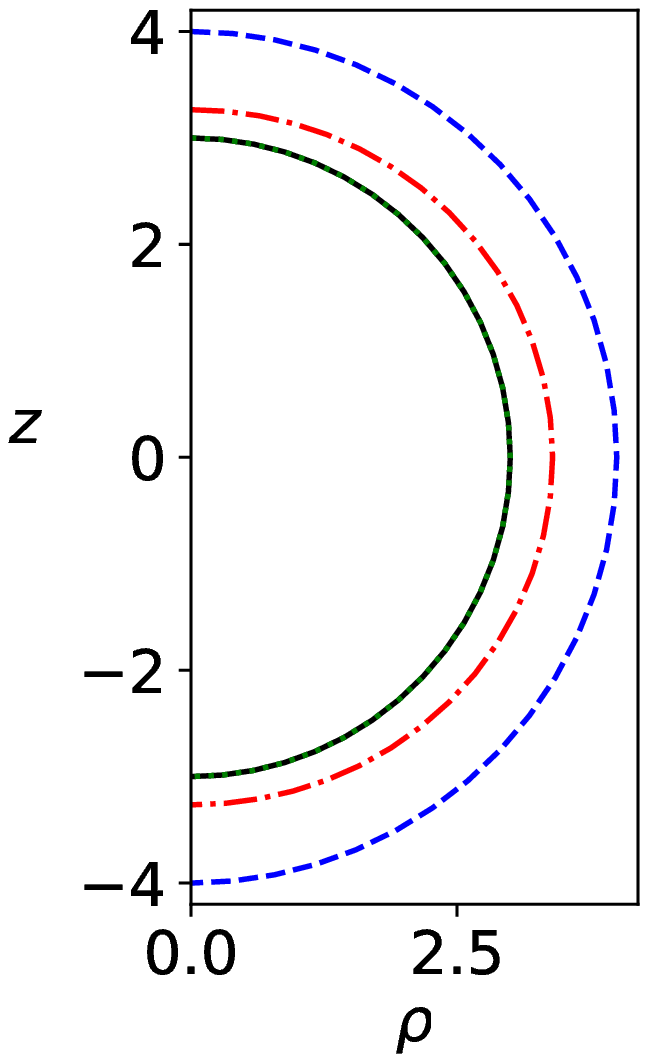}\\
  \includegraphics[width=0.49\textwidth]{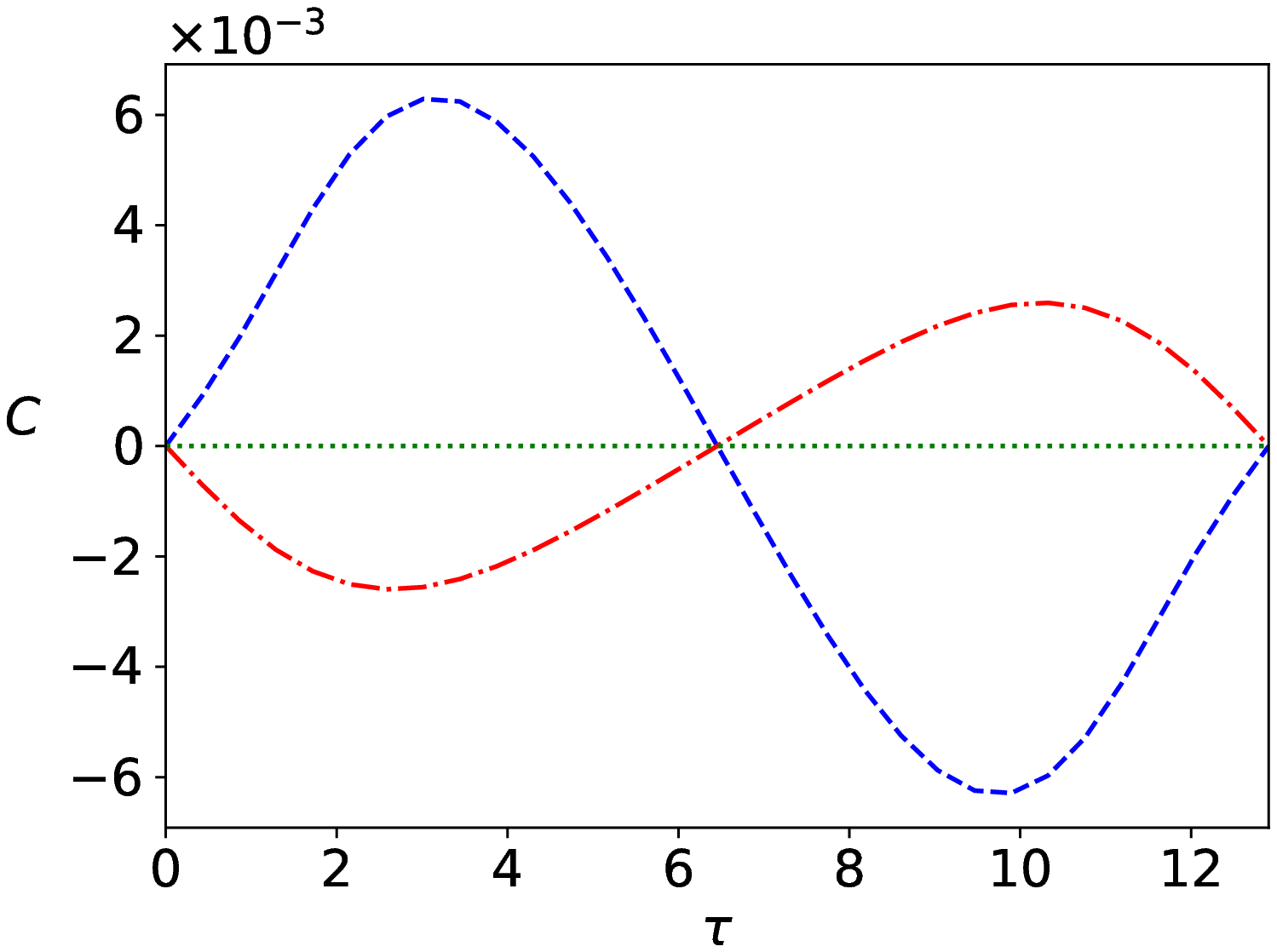}
  \includegraphics[width=0.49\textwidth]{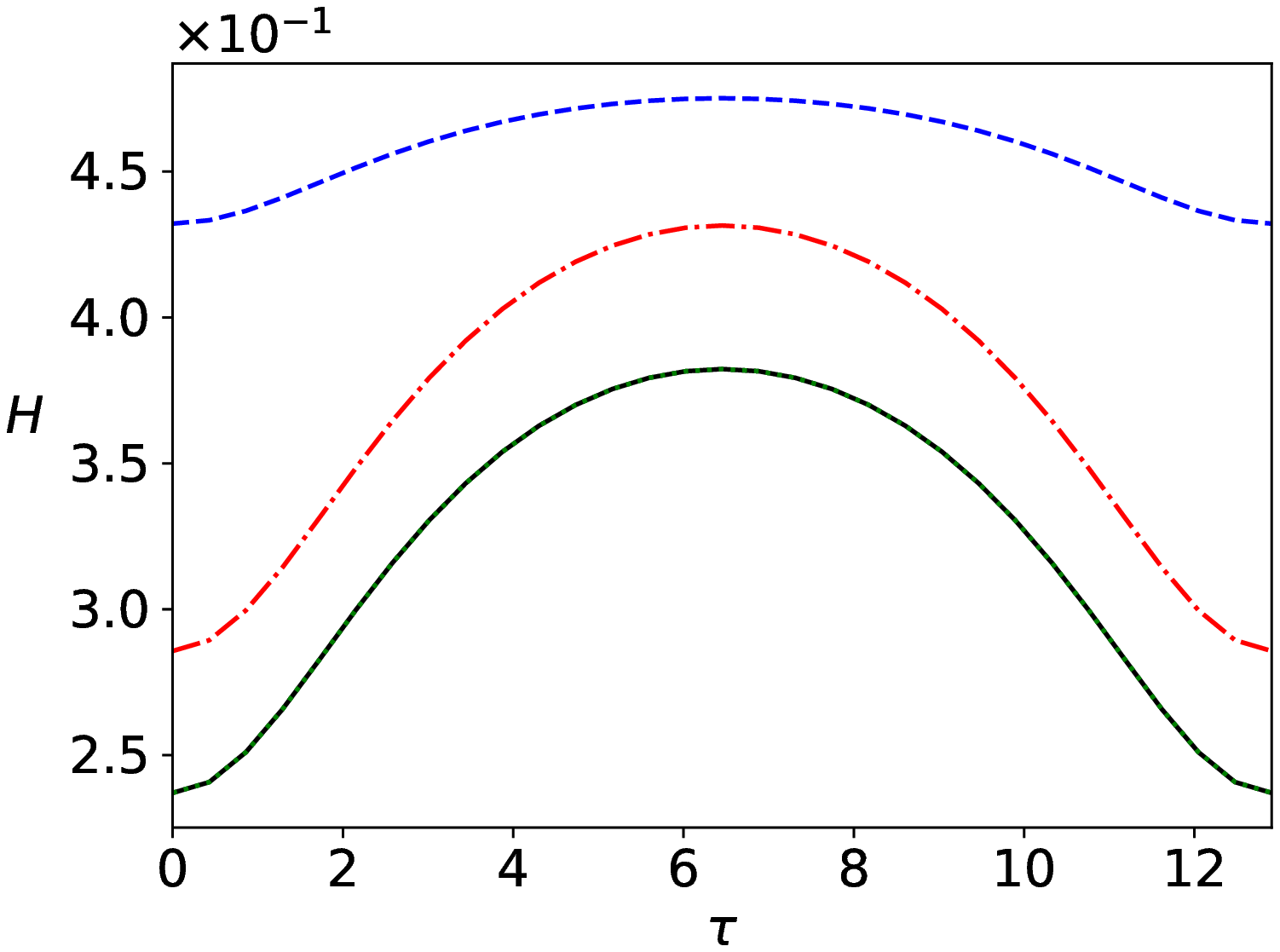}\\
  \includegraphics[width=0.49\textwidth]{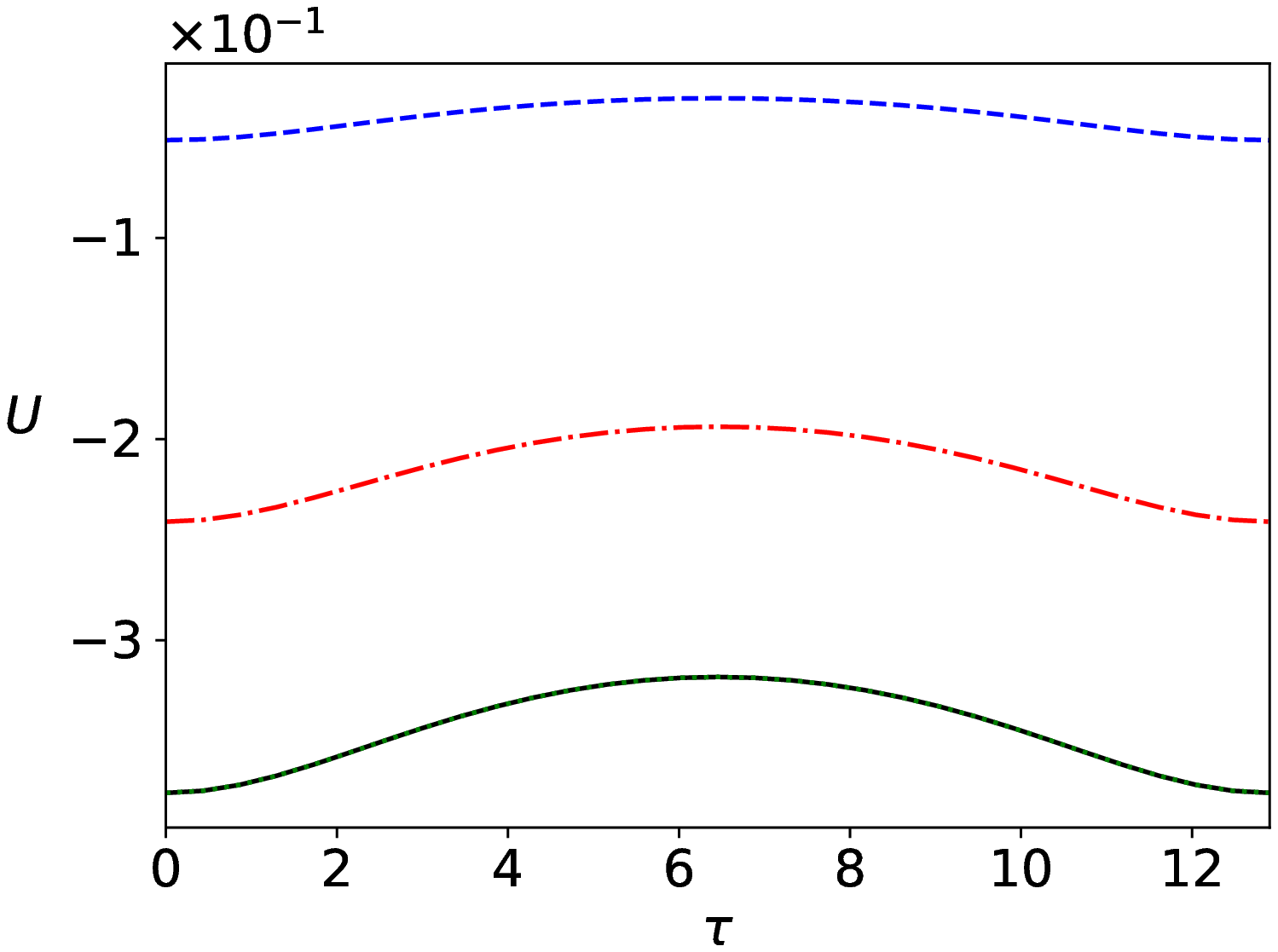}
  \includegraphics[width=0.49\textwidth]{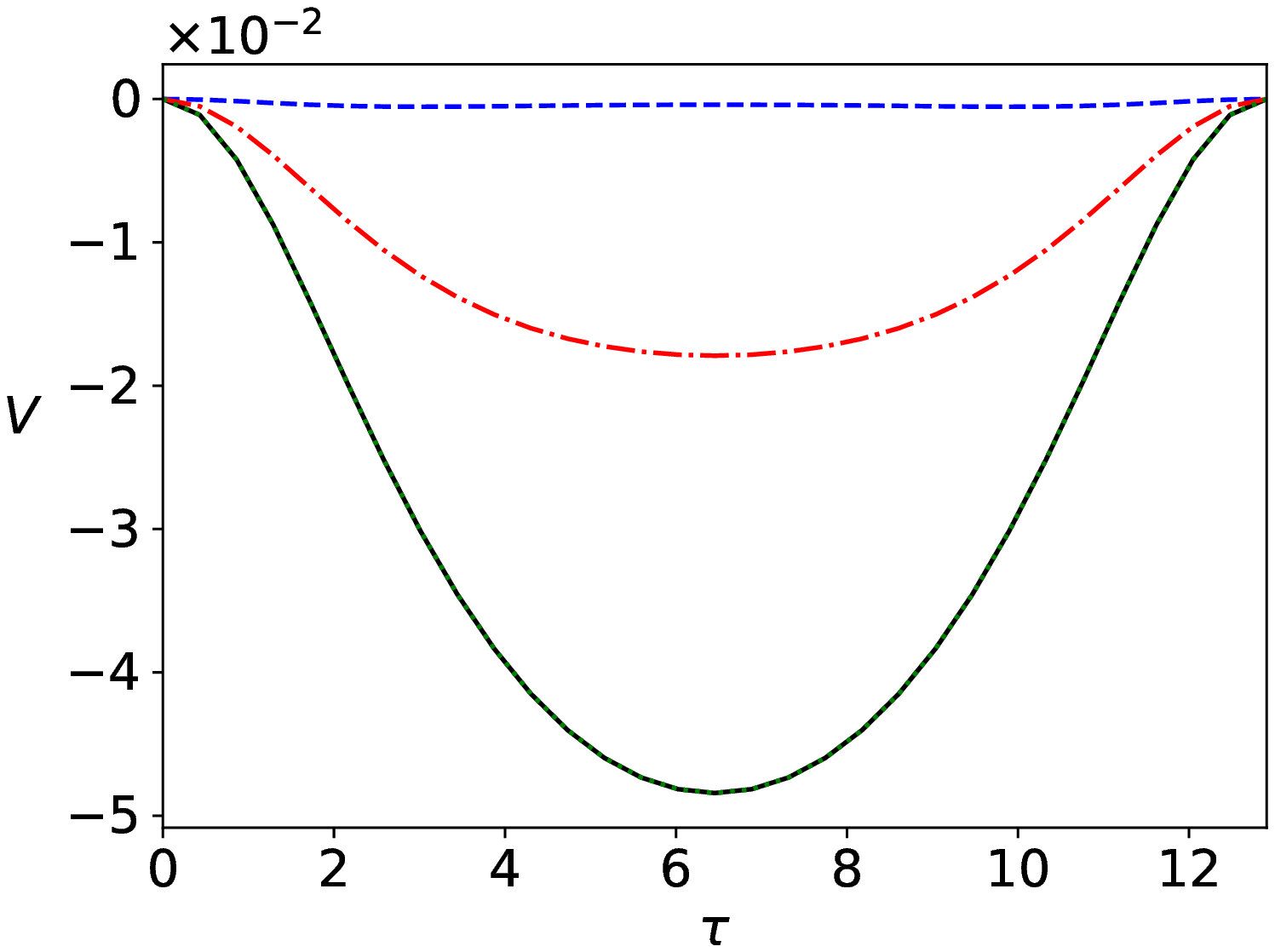}\\
  \includegraphics[width=0.49\textwidth]{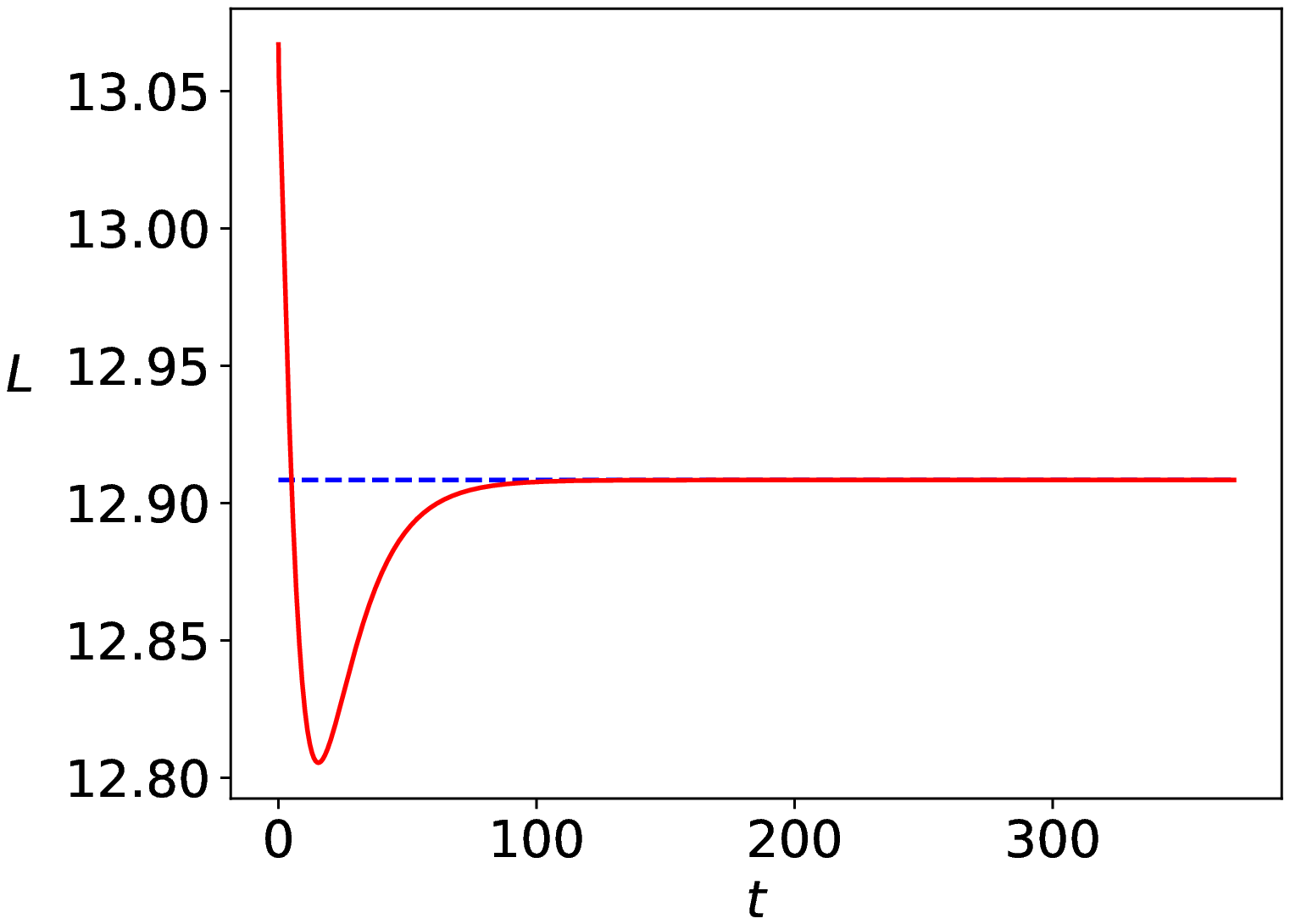}
  \includegraphics[width=0.49\textwidth]{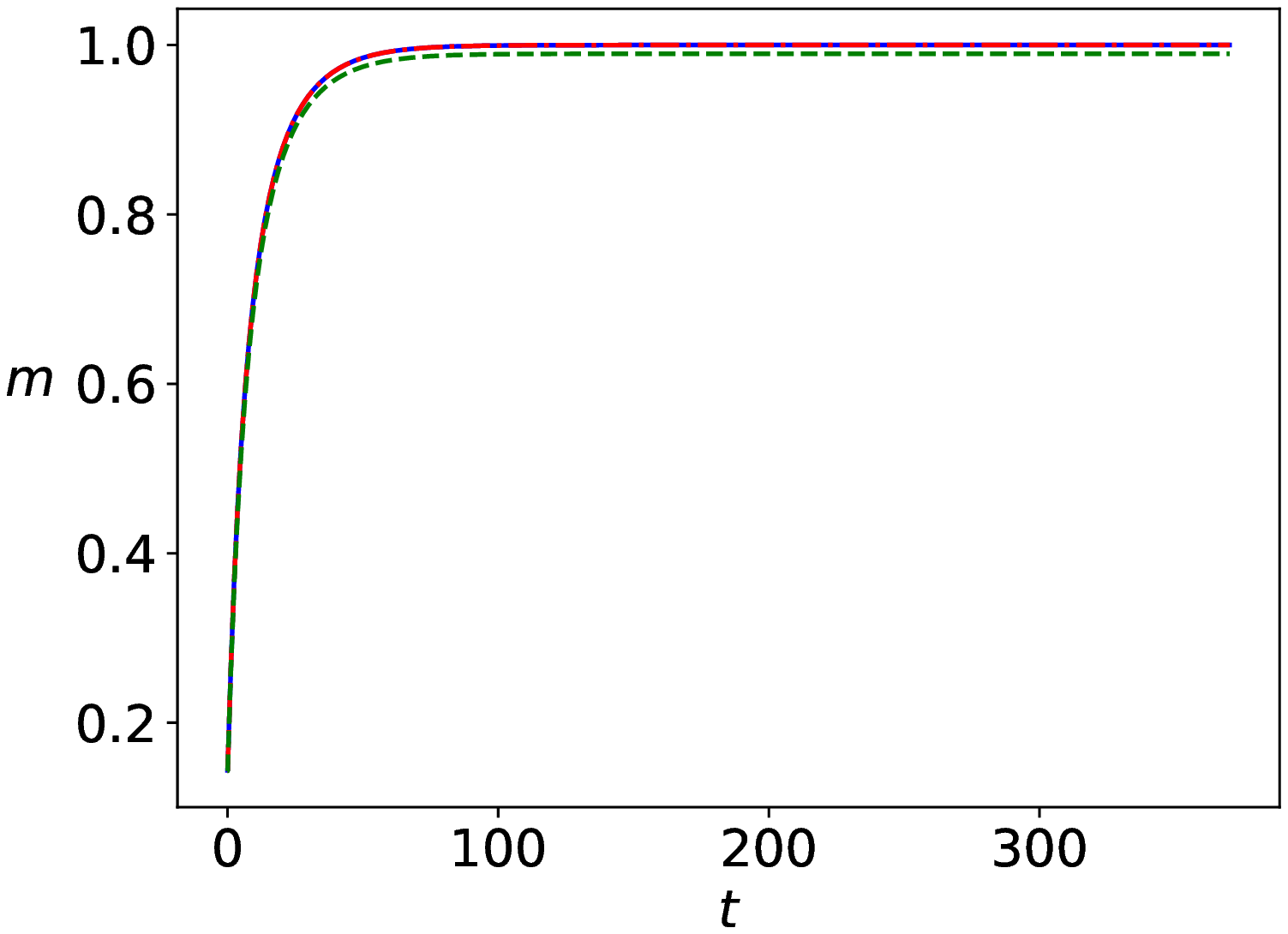}
  \caption{\small Flow with evolving metric.
    The initial curve is taken to be a Weyl--Papapetrou coordinate circle ($r_0 = 4$) 
    in the Zipoy--Voorhees space with
    parameters $M=1, \delta=0.7$,
    the target curve a Weyl--Papapetrou coordinate circle ($\overline{r} = 3$)
    in the Zipoy--Voorhees space with parameters $M=1, \delta=0.6$.
    The same quantities as in Figure~\ref{f:fly_Euclidean_circle_to_ellipse_exz} are
    plotted.
    In the first five panels, the different curves correspond to 
    flow times $t=0$ (dashed blue), $t=9.3$ (dash-dotted red) and 
    $t=370.3$ (dotted green), with the target solution plotted in solid black.
  }
  \label{f:fly_zv_d_07_to_06}
\end{figure}

\begin{figure}[t]
  \centering
  \includegraphics[width=0.49\textwidth]{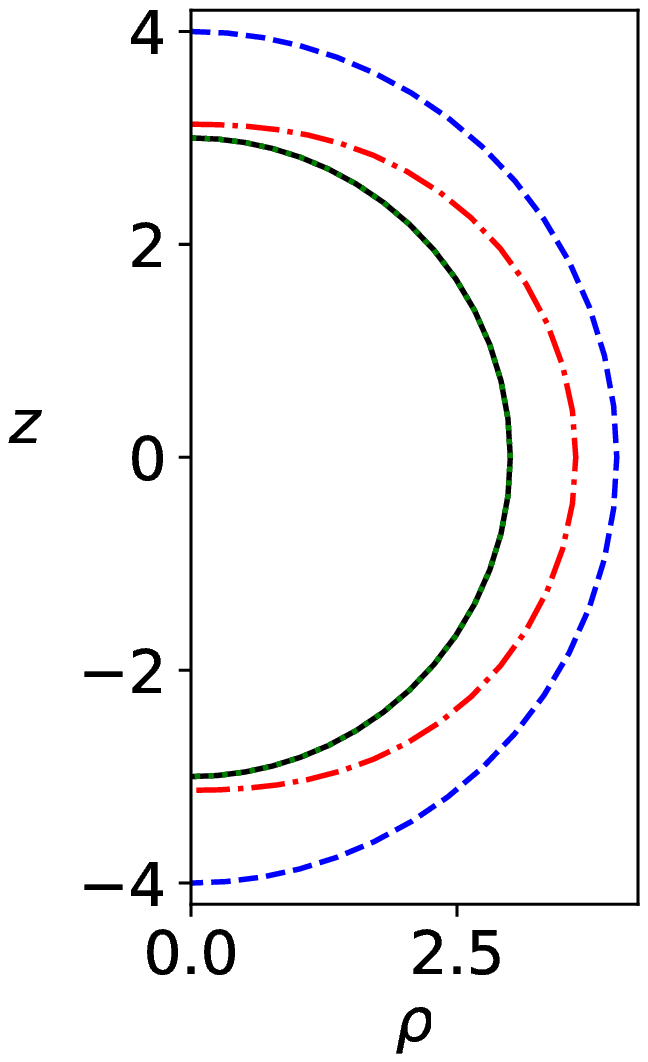}\\
  \includegraphics[width=0.49\textwidth]{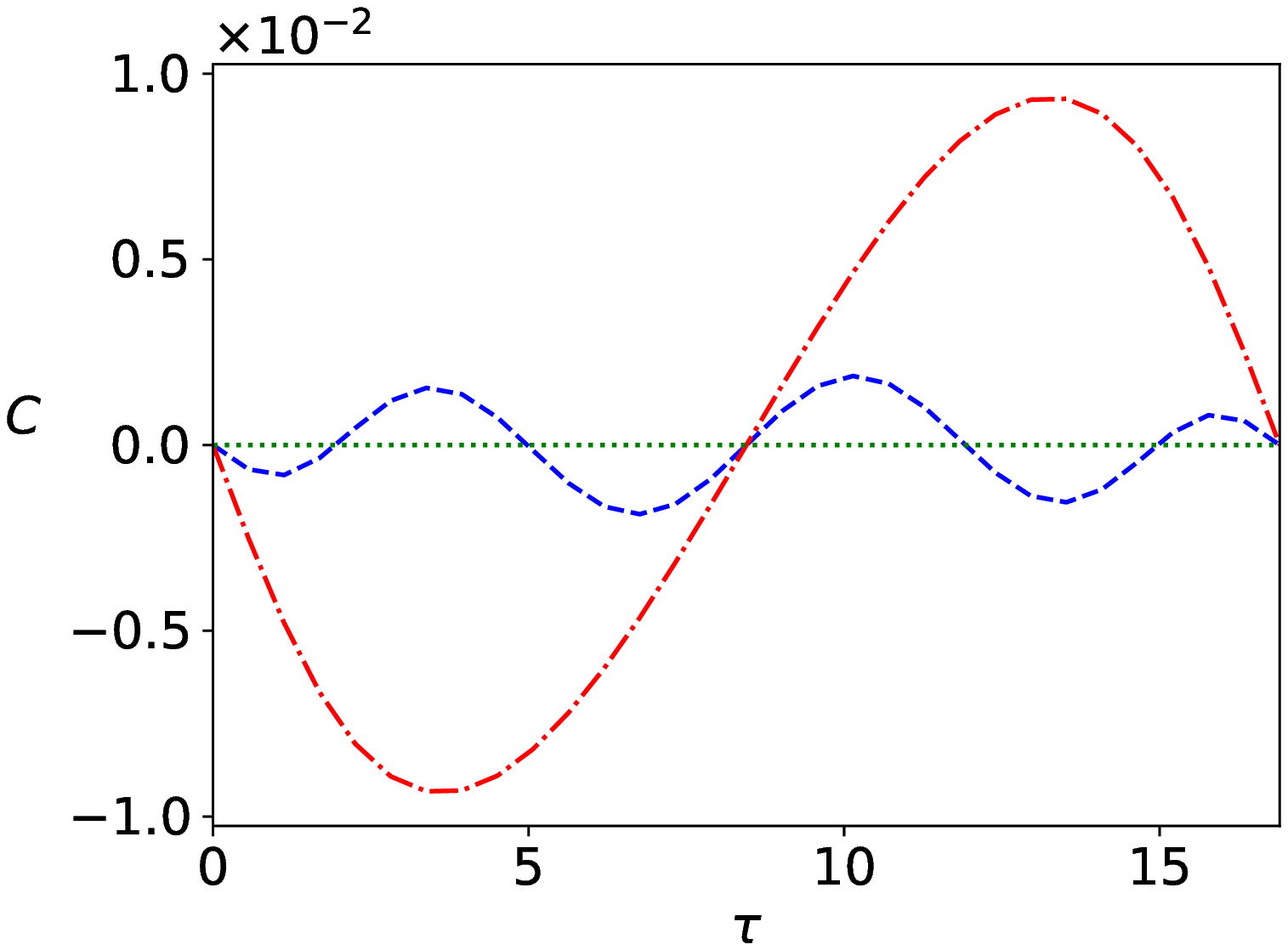}
  \includegraphics[width=0.49\textwidth]{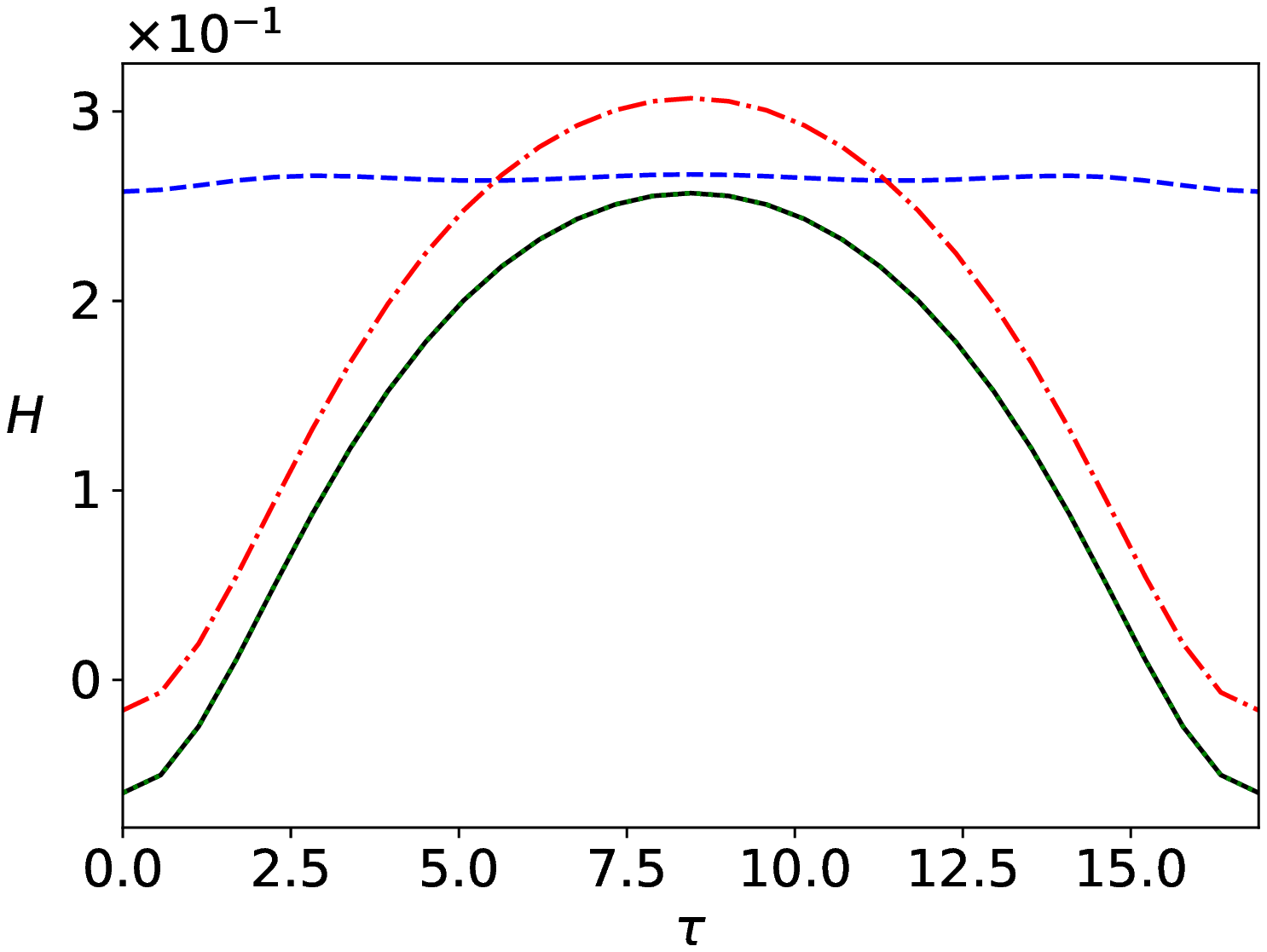}\\
  \includegraphics[width=0.49\textwidth]{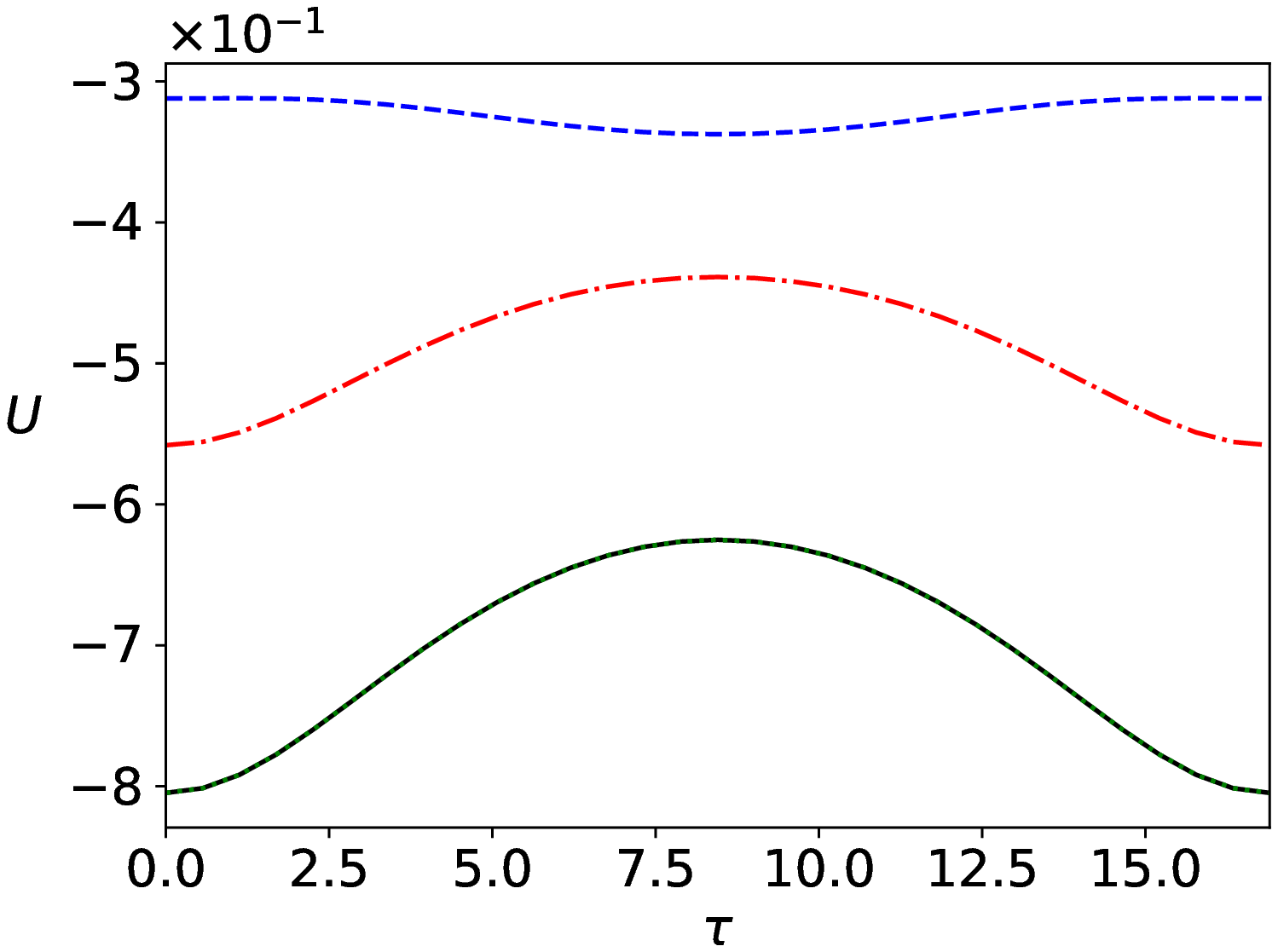}
  \includegraphics[width=0.49\textwidth]{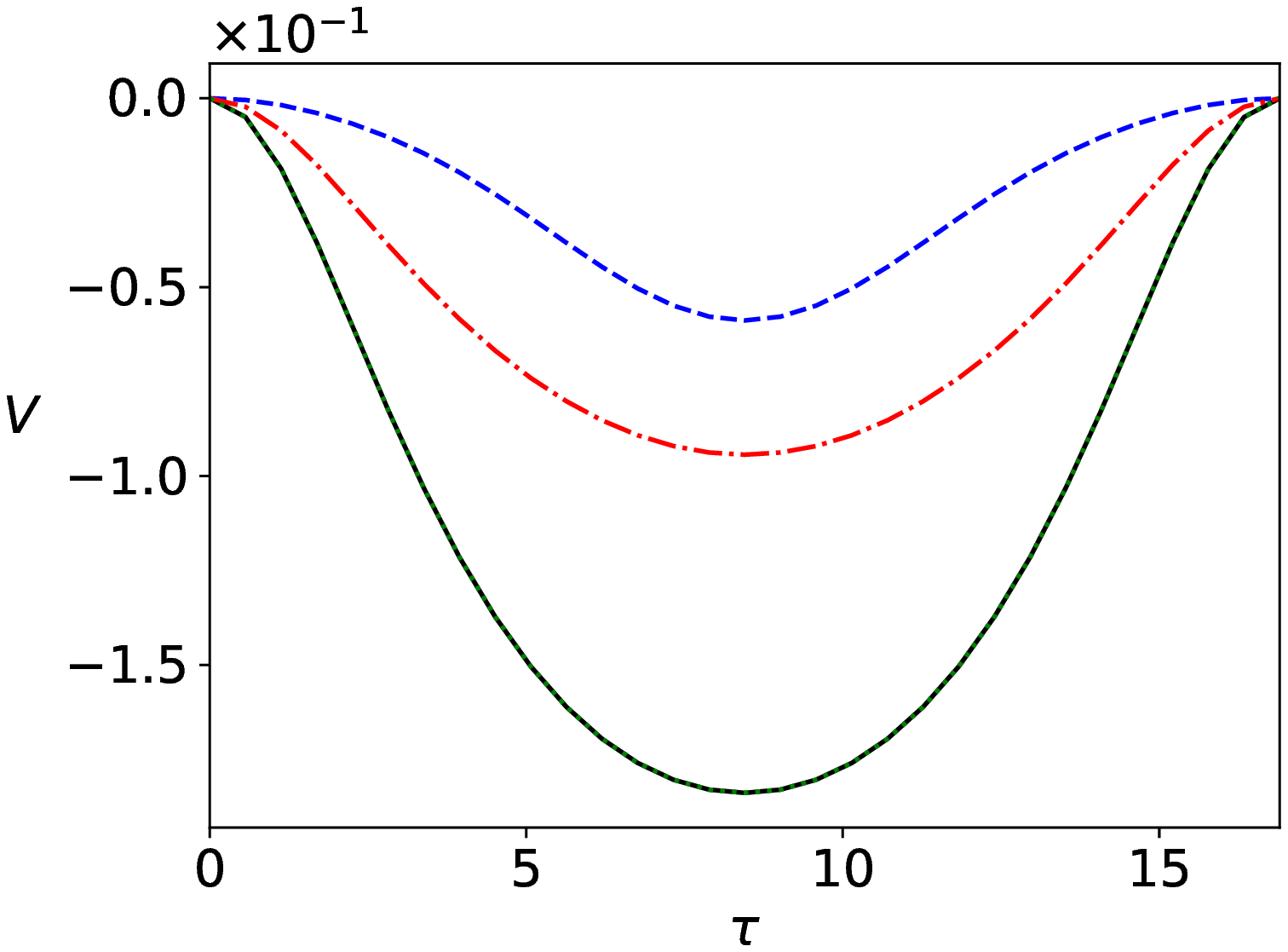}\\
  \includegraphics[width=0.49\textwidth]{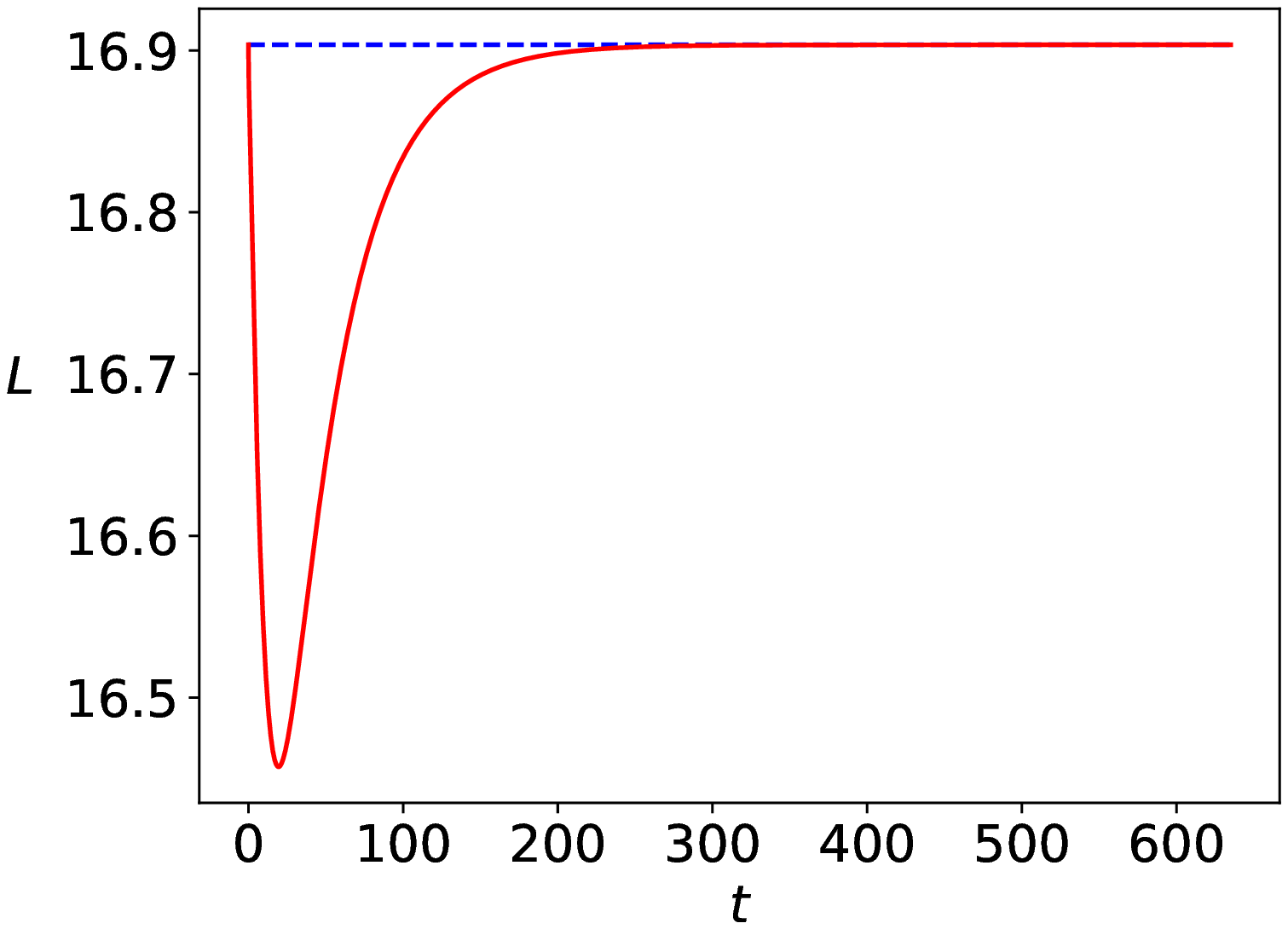}
  \includegraphics[width=0.49\textwidth]{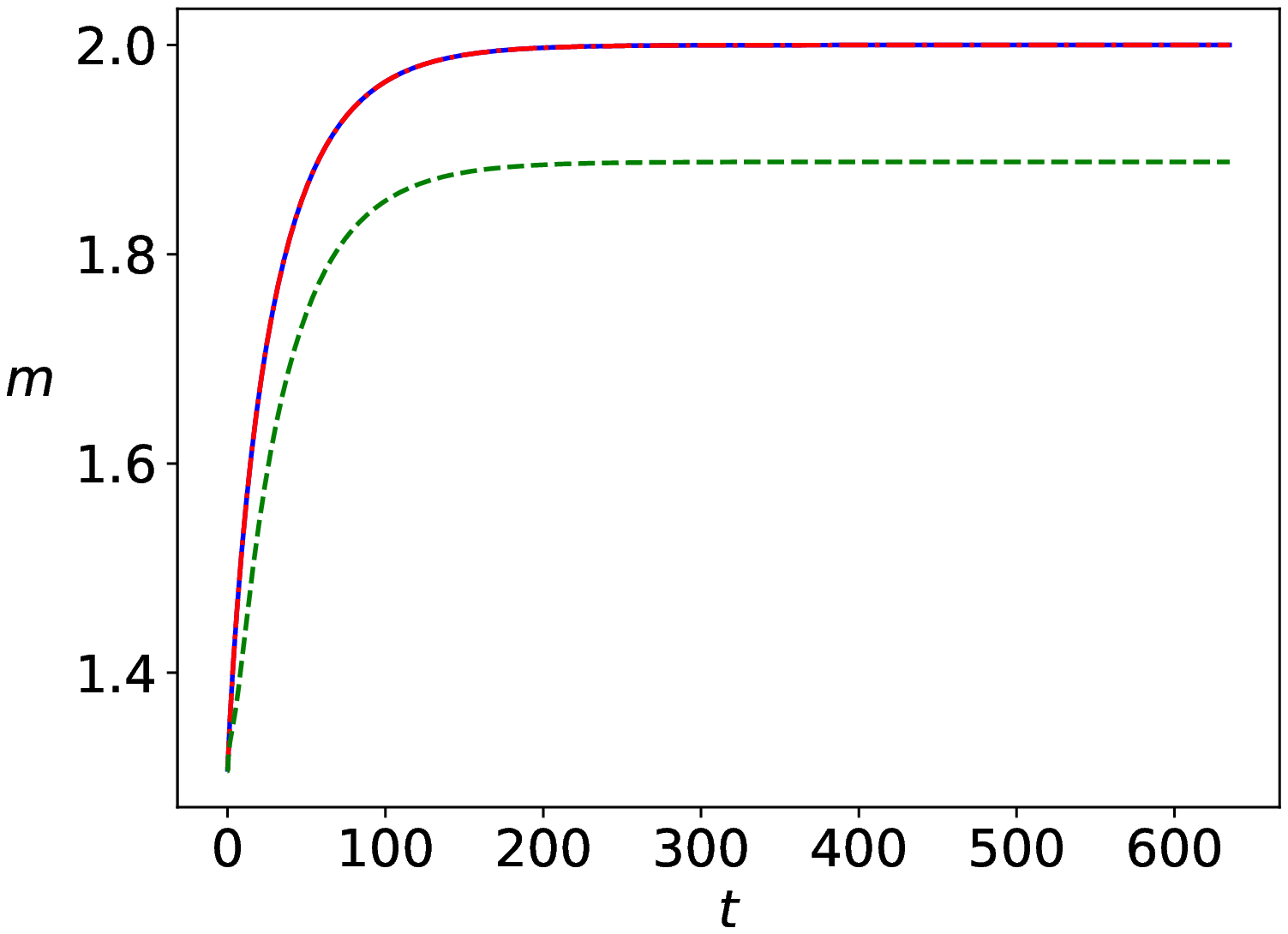}
  \caption{\small Flow with evolving metric.
    The initial curve is taken to be a Weyl--Papapetrou coordinate circle ($r_0 = 4$) 
    in the Curzon--Chazy space of mass $M=1$,
    the target curve a Weyl--Papapetrou coordinate circle ($\overline{r} = 3$)
    in the Curzon--Chazy space of mass $M=2$.
    The same quantities as in Figure~\ref{f:fly_Euclidean_circle_to_ellipse_exz} are
    plotted.
    In the first five panels, the different curves correspond to 
    flow times $t=0$ (dashed blue), $t=15.9$ (dash-dotted red) and 
    $t=634.9$ (dotted green), with the target solution plotted in solid black.
  }
  \label{f:fly_cc_1_to_2}
\end{figure}

A flow between different members of the Zipoy--Voorhees family is shown
in Figure~\ref{f:fly_zv_d_07_to_06} and between different members 
of the Curzon--Chazy family in Figure~\ref{f:fly_cc_1_to_2}.
In this case, the Hawking mass can be seen to differ noticeably
from the other mass functions, as was to be expected.

\clearpage
\subsection{Perturbed Bartnik data}\label{s:perturbed}
In the simulations presented so far, we constructed the Bartnik data
from a given Weyl--Papapetrou system, and the flow correctly recovered the
corresponding metric extension.
To conclude this section, we now study a case where we do not know
\textit{a priori} which exterior metric the Bartnik data give rise to.
To construct such data, we start with a Schwarzschild coordinate
circle at the photon sphere $r_{S\circ}=3M$. 
(We use the photon sphere here since it is a geometrically
distinguished round sphere where the mean curvature $H_{\circ}$ is maximal.)
We compute the corresponding Killing vector norm $\overline{\lambda}_\circ (\tau)$
and then perturb this function:
\begin{equation}
  \overline{\lambda}(\tau) = \overline{\lambda}_\circ(\tau) (1 + f(\tau)),
\end{equation}
where we choose a Gaussian profile
\begin{equation}
  \label{e:pert_exp}
  f(\tau) = A \exp \left( - \left( \frac{\tau-\tau_0}{\sigma} \right) \right).
\end{equation}
We compute a new Schwarzschild coordinate radius $r_S$ from
\begin{equation}
  r_S \definedas \sqrt{\frac{\vert\Sigma\vert}{4\pi}}, \qquad 
  \vert\Sigma\vert \definedas 2\pi \int_0^{\overline{L}} \overline{\lambda}(\tau) \, \diff \tau,
\end{equation}
where $\overline{L}$ is taken from the \emph{unperturbed} target curve
(note that we do not know the coordinate location of the target curve
corresponding to the perturbed Bartnik data in this case).
 We choose the mean curvature to have the constant value
  \begin{equation}
    \overline{H}(\tau) = \frac{2}{\sqrt{3} \, r_S},
  \end{equation}
  the same functional dependence as between $H_{\circ}$ and
  $r_{S\circ}$.
  
\begin{figure}[t]
  \centering
  \includegraphics[width=0.49\textwidth]{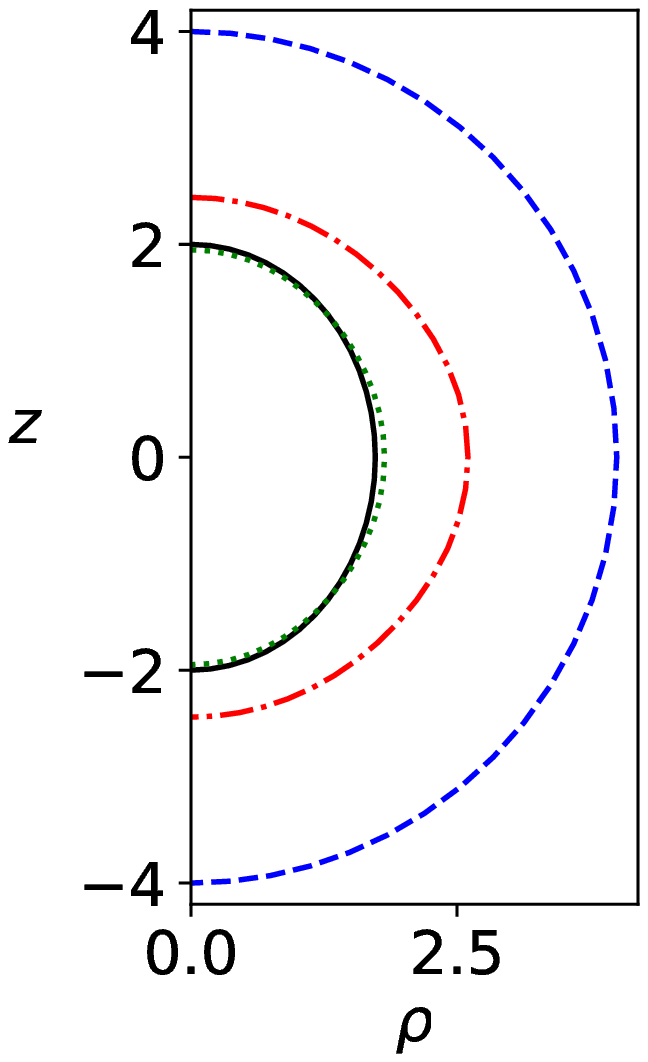}\\
  \includegraphics[width=0.49\textwidth]{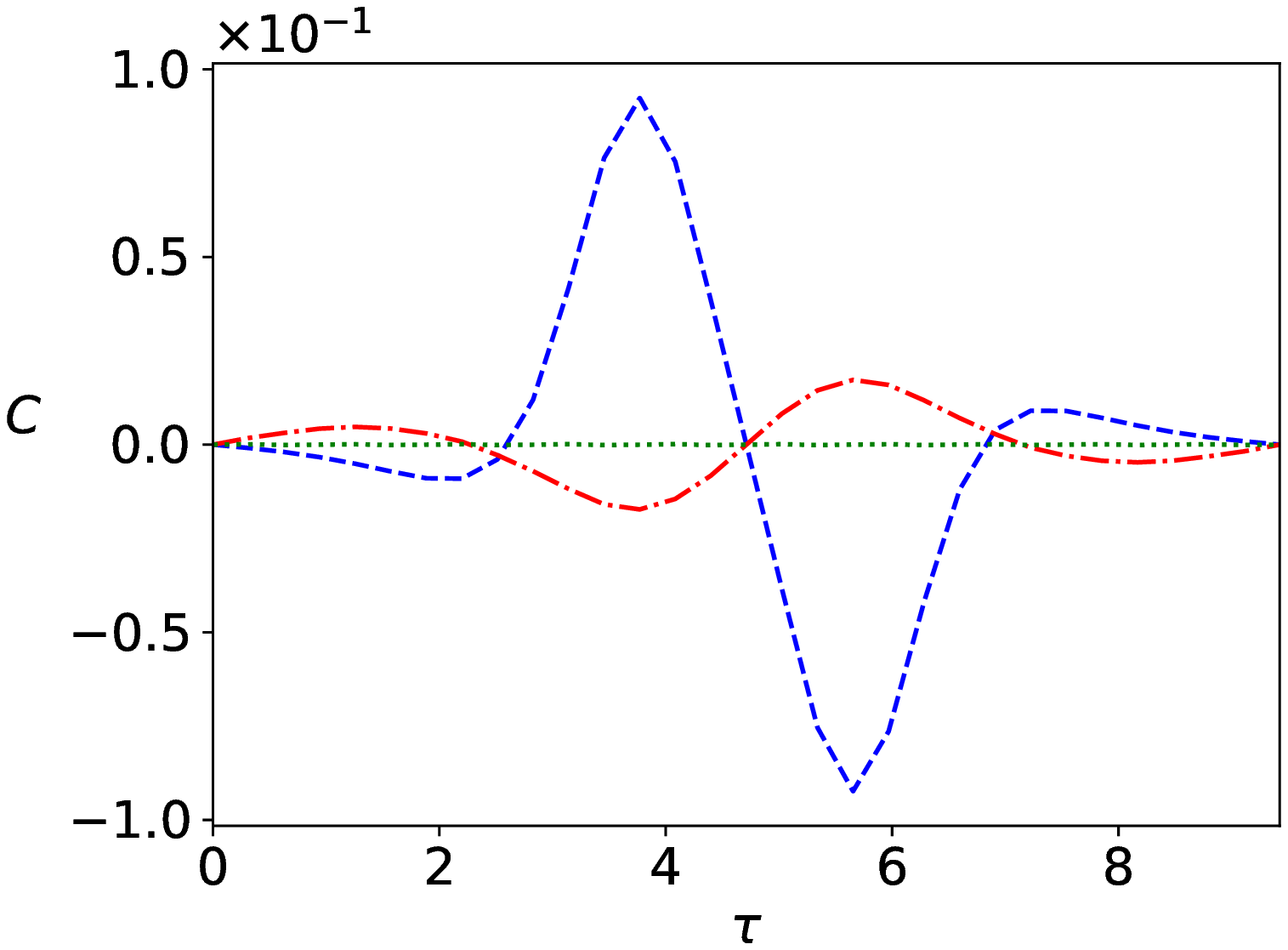}
  \includegraphics[width=0.49\textwidth]{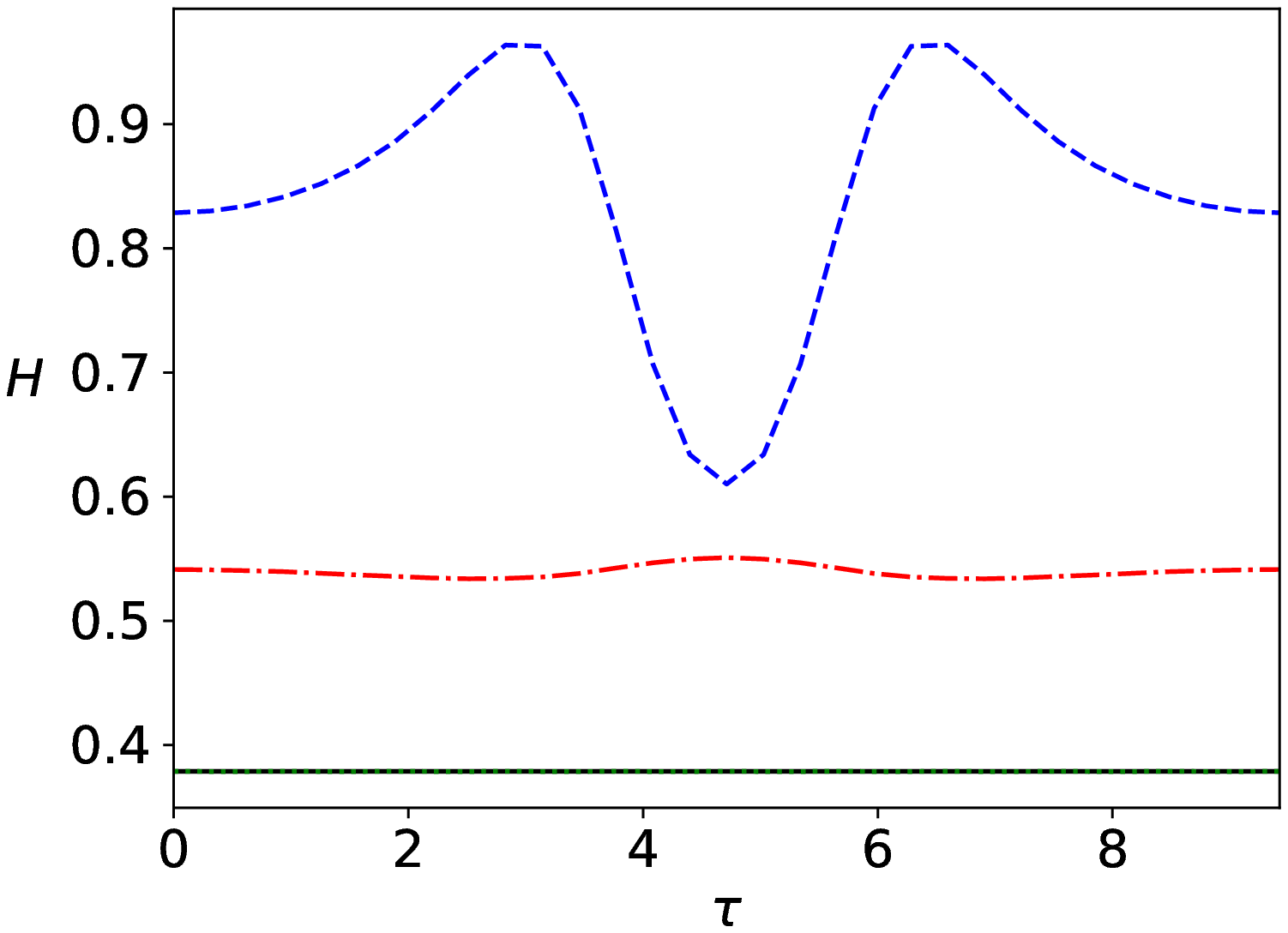}\\
  \includegraphics[width=0.49\textwidth]{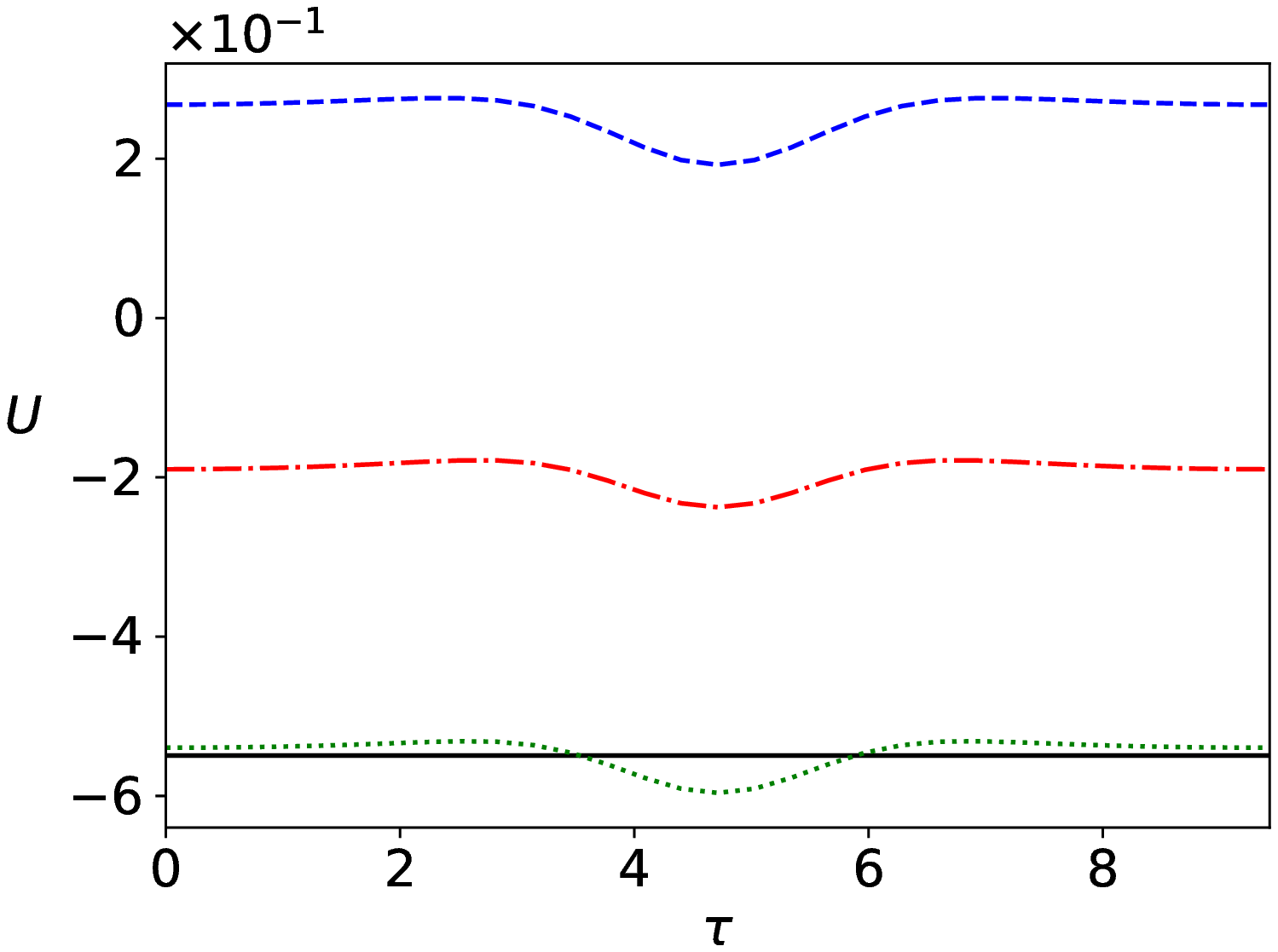}
  \includegraphics[width=0.49\textwidth]{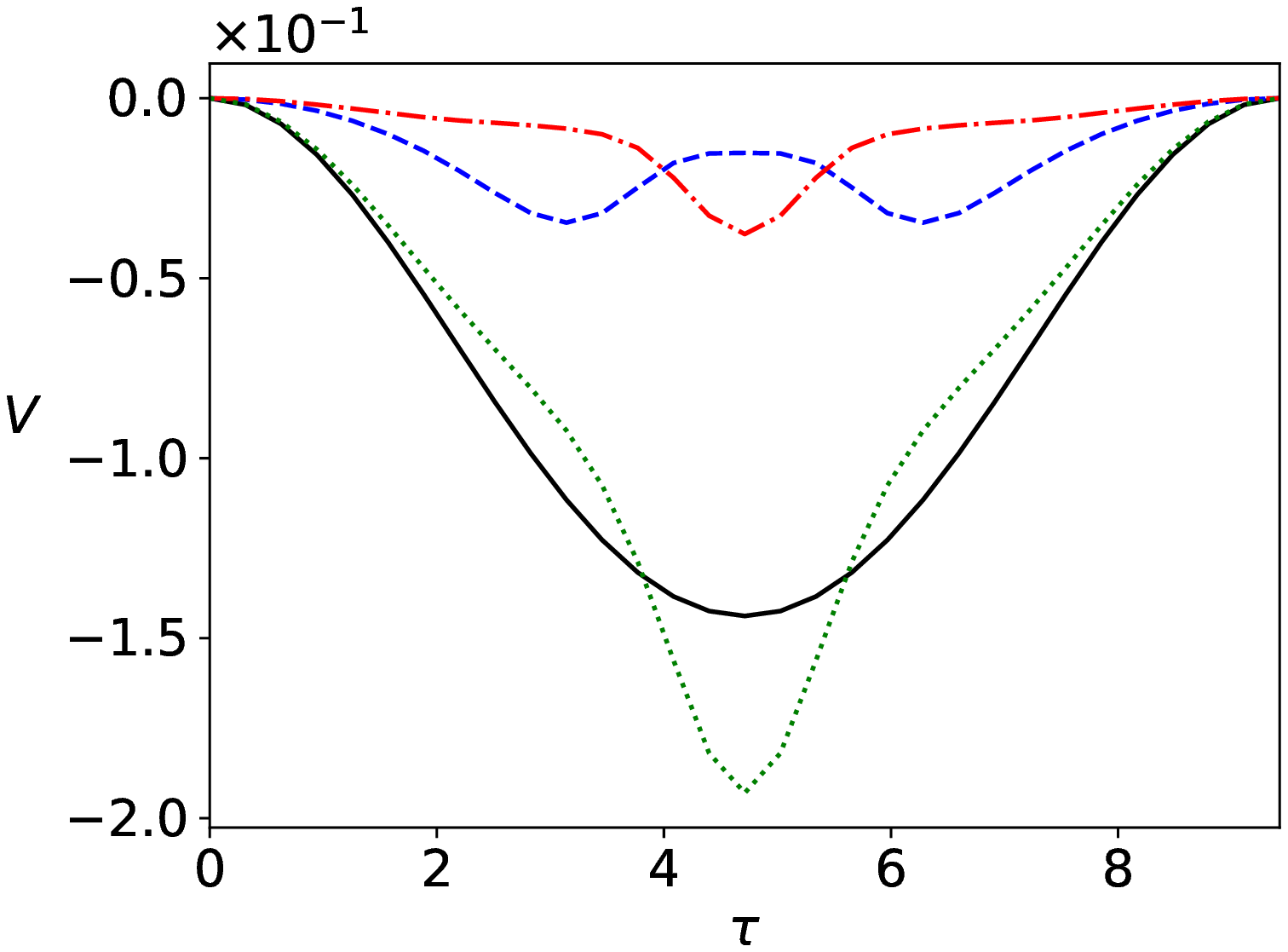}\\
  \includegraphics[width=0.49\textwidth]{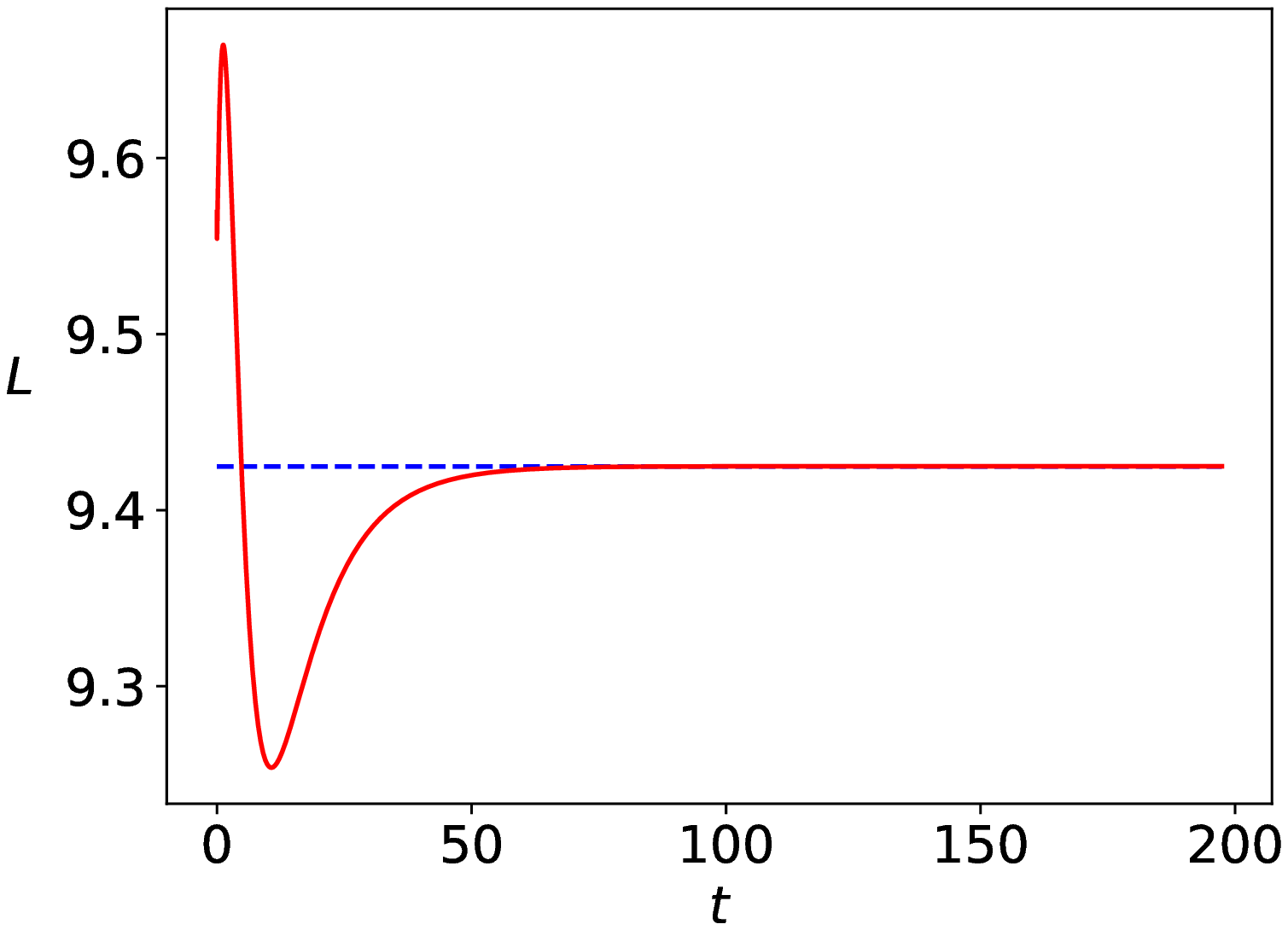}
  \includegraphics[width=0.49\textwidth]{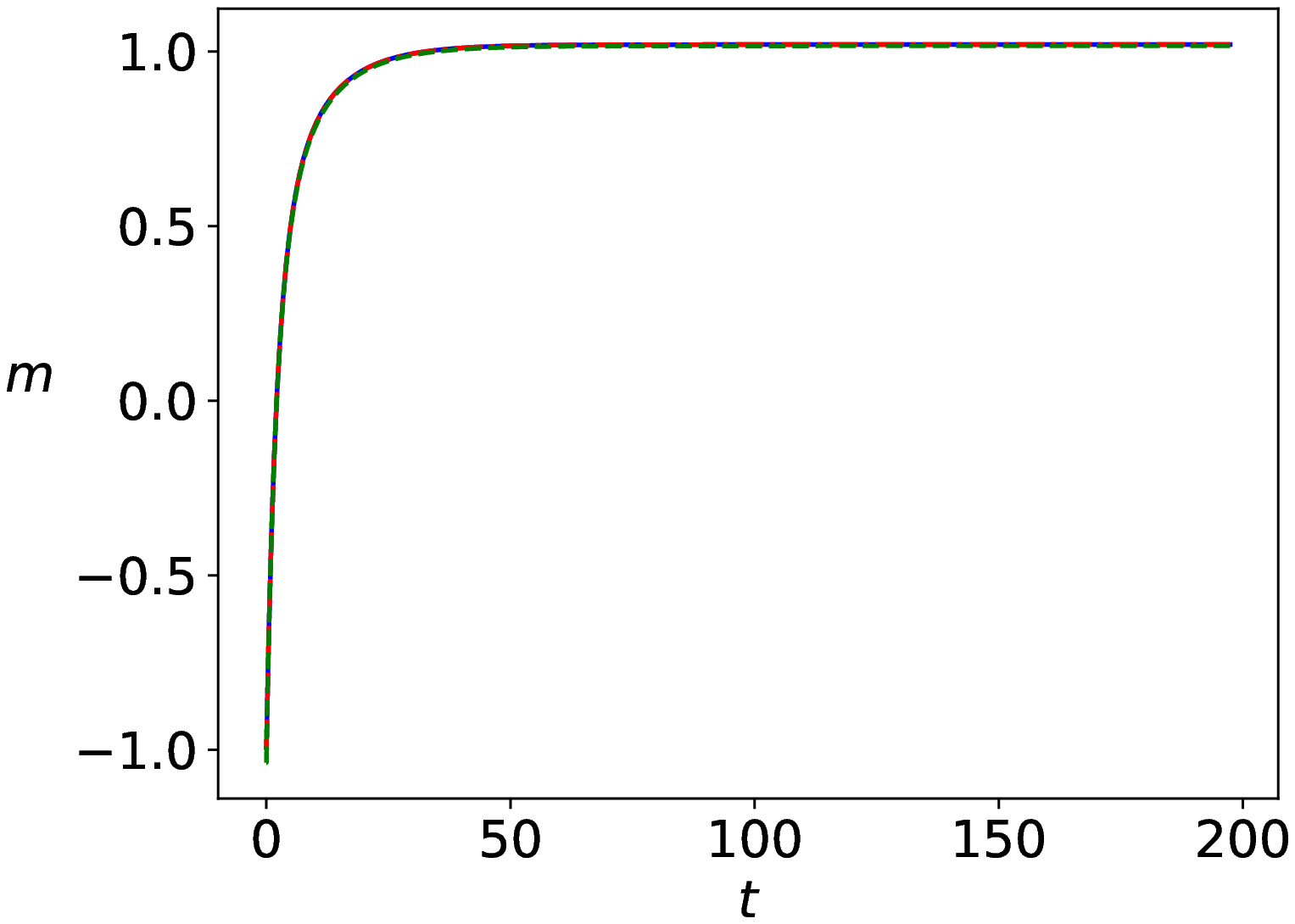}
  \caption{\small Flow with evolving metric. 
    The Bartnik data are constructed
    from the photon sphere in $M=1$ Schwarzschild perturbed by a
    Gaussian profile \eqref{e:pert_exp} with $A=0.1$,
    $\tau_0 = \overline{L}/2$ and $\sigma = \overline{L}/8$.
    The initial curve is taken to be a circle with radius $r_S = 4$ in 
    $M=1$ Schwarzschild.
    Solid black lines refer to the unperturbed data.
    The final ADM mass belonging to the perturbed Bartnik data is 
    $m_\mathrm{ADM} = 1.0200$ and the Hawking mass is
    $m_\mathrm{H} = 1.0162$.
    In the first five panels, the different curves correspond to 
    flow times $t=0$ (dashed blue), $t=4.9$ (dash-dotted red) and 
    $t=197.4$ (dotted green).
    Solid black lines refer to the \emph{unperturbed} data here.
  }
  \label{f:fly_pert_exp}
\end{figure}

Figure \ref{f:fly_pert_exp} demonstrates that
the flow converges; we have thus constructed the static metric extension
corresponding to the prescribed Bartnik data 
$([0,\overline{L}],\overline{\lambda},\overline{H})$.
The chosen amplitude $A=0.1$ of the perturbation is the maximum value
for which we were able to achieve a stable numerical evolution. 
We choose $M=1$ for the unperturbed target curve.
The final ADM mass is $m_\mathrm{ADM} = 1.0200$ and the final Hawking 
mass is $m_\mathrm{H} = 1.0162$.
This is in accordance with the generalised Penrose inequality~\eqref{genPI}.

\begin{figure}[t]
  \centering
 \includegraphics[width=0.7\textwidth]{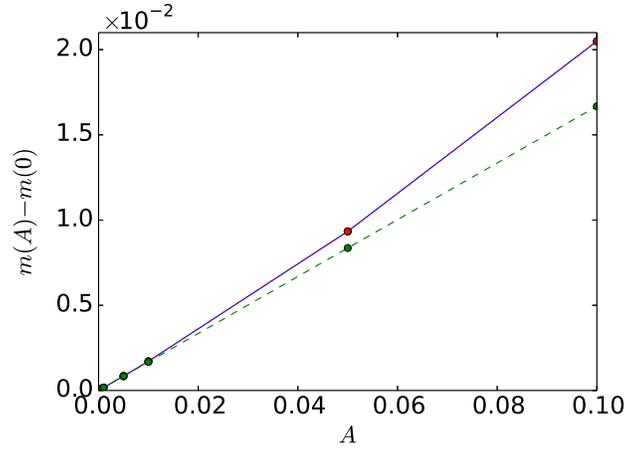}
  \caption{Dependence of the masses on the amplitude $A$ of the
    perturbation for the Bartnik data with Gaussian perturbation as in 
    Figure~\ref{f:fly_pert_exp}.
    Shown are the ADM mass (solid blue), the Hawking mass (dashed
    green), and the pseudo-Newtonian mass (dash-dotted red,
    indistinguishable from the solid blue line).
  }
  \label{f:mass_Amplitude}
\end{figure}

\clearpage
In Figure~\ref{f:mass_Amplitude}, we investigate the dependence of the
masses (ADM, Hawking, and pseudo-Newtonian) on the amplitude $A$ of the
perturbation.
All three masses appear to approach the unperturbed value continuously 
as $A\to 0$.
The ADM mass and the pseudo-Newtonian mass are identical within
numerical error, as expected from the theoretical analysis, 
see Section~\ref{s:masses}.

\clearpage

%% file: paper_sec6.tex
\section{Conclusion}
\label{s:discussion}

In this paper, we developed a new approach to construct static metric 
extensions as they arise in Bartnik's conjecture.
We restricted ourselves to axisymmetry and worked with the Weyl--Papapetrou 
formulation of the static axisymmetric vacuum (SAV) Einstein equations.
The metric extension problem becomes an elliptic free boundary value
problem in this setting, which we solved numerically using a geometric flow
\eqref{e:flow} coupled to the Weyl--Papapetrou system of equations \eqref{e:Ueqn}, \eqref{e:Veqns}.
As far as we know, this is the first time static metric extensions
have been constructed explicitly in general situations using numerical
methods.

It should be noted that we only considered axisymmetric Bartnik
data, and we only sought axisymmetric static metric extensions.
Even for axisymmetric data, non-axisymmetric static metric extensions might exist. Furthermore, we restricted our attention to reflection-symmetric Bartnik data and static metric extensions. Of course, even for reflection symmetric data, non-reflection symmetric static metric extensions may exist.

In a first step, we simplified the situation by fixing the metric to
some known background SAV solution to the Einstein equations and prescribing
Bartnik data corresponding to a given surface within that spacetime,
or rather to a given curve after performing the symmetry reduction.
We showed analytically that coordinate circles in Euclidean space and
in a fixed Schwarzschild background remain coordinate circles during our geometric
flow and approach the desired target circles in a globally stable manner.
Moreover, coordinate circles in Euclidean space are linearly stable against arbitrary perturbations.
We also presented arguments in favour of short-time existence of
solutions to \eqref{e:flow}. It would be very interesting to study short-time existence of the flow
more rigorously.

Numerically, our geometric flow behaved as expected in this situation. Interestingly, in Euclidean space, the Euclidean distance \eqref{e:distance}
between the flowing and the target curve appeared to decrease
monotonically in all cases we studied numerically.
Investigating this claim analytically would be an interesting topic
for future work.

Next, we studied the full elliptic free boundary value problem, which
involved solving the Weyl--Papapetrou form \eqref{e:Ueqn}, \eqref{e:Veqns} of the SAV Einstein
equations at each ``time'' step of the flow.
We specified Bartnik data corresponding to given surfaces in various known
SAV spacetimes, and in all cases where we obtained stable numerical
evolutions, the surface and static metric extension were found correctly. 
We also perturbed Bartnik data corresponding to (centred) round spheres in
Schwarzschild so that we did not know the corresponding metric 
extensions \textit{a priori}, and we were able to construct static metric 
extensions with ADM masses up to $2\%$ larger than in the unperturbed case. The ADM mass appeared to approach the unperturbed mass continuously in the
limit of vanishing amplitude of the perturbation.
The Hawking mass was always observed to be smaller than or equal to
the ADM mass (or, identically, the pseudo-Newtonian mass), in
agreement with the generalised Penrose inequality~\eqref{genPI}.
These results should be interesting in the light of ongoing analytical
work on the case of near-round spheres in Schwarzschild in~\cite{CederDiaz}.

Our analysis revealed that theoretically,
the flow has more stationary states than just the desired
position of a curve inducing the correct Bartnik data (even in a fixed
Euclidean background), although in all the simulations shown in this
paper, the flow did approach the desired stationary state.
Investigating if there really can be evolutions approaching
  one of the spurious stationary states and developing suitable work-around strategies would be an interesting topic for
future research.

In all situations where we investigated this, different initial data
to the coupled flow gave rise to the same asymptotic solution $(U,V)$ as $t\to\infty$.
If we had found an example where this is not the case then this would
disprove the uniqueness of static metric extensions for given Bartnik data.
It should be interesting to investigate this uniqueness question 
in more extreme situations.

We encountered numerical instabilities e.g.~when we specified Bartnik data
corresponding to a surface too close to the horizon in Schwarzschild
 or when the flowing curves became too strongly deformed (e.g. ellipses with
large eccentricity).
Sometimes, it was possible to cure these instabilities by adapting the parameter $\kappa$ in \eqref{e:flow} or by making the time steps
sufficiently small.
Another source of numerical instability was associated with our method
of solving the Poisson equation~\eqref{e:Ueqn} for the metric field $U$ described in
Section~\ref{s:nummethod.einstein} and arose when the radius of the
flowing curve became too small, $r\lesssim 1$ in Weyl--Papapetrou coordinates.
It was possible to somewhat alleviate this problem by the least squares method
(also described in Section~\ref{s:nummethod.einstein}).
More work is needed to obtain stable simulations in more extreme situations.

\enlargethispage{0.5cm}
 From an analytical perspective, it would be very interesting indeed
 to rigorously analyse the 
full coupled elliptic system with flowing boundary,
equations \eqref{e:Ueqn}, \eqref{e:Veqns} and \eqref{e:flow}, 
 and to thereby obtain theoretical results about existence of
 solutions to Bartnik's static metric extension conjecture in
 Weyl--Papapetrou form.

\newpage

%% file: paper.bbl
\begin{thebibliography}{26}

\bibitem{Anderson}
\begin{barticle}[author]
\bauthor{\bsnm{Anderson},~\bfnm{M.}\binits{M.}}
(\byear{2015}).
\btitle{Local existence and uniqueness for exterior static vacuum Einstein
  metrics}.
\bjournal{Proceedings of the American Mathematical Society}
\bvolume{143}
\bpages{3091--3096}.
\end{barticle}
\endbibitem

\bibitem{Anderson_2013}
\begin{barticle}[author]
\bauthor{\bsnm{Anderson},~\bfnm{Michael~T}\binits{M.~T.}} \AND
  \bauthor{\bsnm{Khuri},~\bfnm{Marcus~A}\binits{M.~A.}}
(\byear{2013}).
\btitle{On the Bartnik extension problem for the static vacuum Einstein
  equations}.
\bjournal{Classical and Quantum Gravity}
\bvolume{30}
\bpages{125005}.
\bdoi{10.1088/0264-9381/30/12/125005}
\end{barticle}
\endbibitem

\bibitem{ADM}
\begin{barticle}[author]
\bauthor{\bsnm{Arnowitt},~\bfnm{R.}\binits{R.}},
  \bauthor{\bsnm{Deser},~\bfnm{S.}\binits{S.}} \AND
  \bauthor{\bsnm{Misner},~\bfnm{Ch.}\binits{C.}}
(\byear{1961}).
\btitle{{Coordinate Invariance and Energy Expressions in General Relativity}}.
\bjournal{Phys. Rev.}
\bvolume{122}
\bpages{997--1006}.
\end{barticle}
\endbibitem

\bibitem{Bartnik86}
\begin{barticle}[author]
\bauthor{\bsnm{Bartnik},~\bfnm{Robert}\binits{R.}}
(\byear{1986}).
\btitle{{The Mass of an Asymptotically Flat Manifold}}.
\bjournal{Communications on Pure and Applied Mathematics}
\bvolume{39}
\bpages{661--693}.
\bdoi{10.1002/cpa.3160390505}
\end{barticle}
\endbibitem

\bibitem{Bartnik89}
\begin{barticle}[author]
\bauthor{\bsnm{Bartnik},~\bfnm{R.}\binits{R.}}
(\byear{1989}).
\btitle{New Definition of Quasilocal Mass}.
\bjournal{Phys.\ Rev.\ Lett.}
\bvolume{62}
\bpages{2346--2348}.
\end{barticle}
\endbibitem

\bibitem{Boyd2001}
\begin{bbook}[author]
\bauthor{\bsnm{Boyd},~\bfnm{J.~B.}\binits{J.~B.}}
(\byear{2001}).
\btitle{Chebyshev and Fourier Spectral Methods},
\bedition{2nd edn.} ed.
\bpublisher{Dover}.
\end{bbook}
\endbibitem

\bibitem{Cederbaum2011}
\begin{bphdthesis}[author]
\bauthor{\bsnm{Cederbaum},~\bfnm{C.}\binits{C.}}
(\byear{2011}).
\btitle{The {N}ewtonian Limit of Geometrostatics.}
\btype{PhD thesis},
\bpublisher{Freie Universit\"at Berlin}.
\end{bphdthesis}
\endbibitem

\bibitem{CederDiaz}
\begin{bmisc}[author]
\bauthor{\bsnm{Cederbaum},~\bfnm{C.}\binits{C.}} \AND
  \bauthor{\bsnm{Escobar~Diaz},~\bfnm{L.}\binits{L.}}
\bnote{Work in progress}.
\end{bmisc}
\endbibitem

\bibitem{Chrus2}
\begin{barticle}[author]
\bauthor{\bsnm{Chr\'usciel},~\bfnm{Piotr~T.}\binits{P.~T.}}
(\byear{1988}).
\btitle{{On the invariant mass conjecture in general relativity}}.
\bjournal{Communications in Mathematical Physics}
\bvolume{120}
\bpages{233--248}.
\end{barticle}
\endbibitem

\bibitem{Courant1928}
\begin{barticle}[author]
\bauthor{\bsnm{Courant},~\bfnm{R.}\binits{R.}},
  \bauthor{\bsnm{Friedrichs},~\bfnm{K.}\binits{K.}} \AND
  \bauthor{\bsnm{Lewy},~\bfnm{H.}\binits{H.}}
(\byear{1928}).
\btitle{\"Uber die partiellen {D}ifferenzengleichungen der mathematischen
  {P}hysik}.
\bjournal{Math.\ Ann.}
\bvolume{100}
\bpages{32--74}.
\end{barticle}
\endbibitem

\bibitem{Dittberner2017}
\begin{bphdthesis}[author]
\bauthor{\bsnm{Dittberner},~\bfnm{F.}\binits{F.}}
(\byear{2017}).
\btitle{Constrained Curve Flows.}
\btype{PhD thesis},
\bpublisher{Freie Universit\"at Berlin}.
\end{bphdthesis}
\endbibitem

\bibitem{Fornberg1996}
\begin{bbook}[author]
\bauthor{\bsnm{Fornberg},~\bfnm{B.}\binits{B.}}
(\byear{1996}).
\btitle{A Practical Guide to Pseudospectral Methods}.
\bseries{Cambridge Monographs on Applied and Computational Mathematics}.
\bpublisher{Cambridge University Press}.
\end{bbook}
\endbibitem

\bibitem{Gage1986}
\begin{bincollection}[author]
\bauthor{\bsnm{Gage},~\bfnm{M.~E.}\binits{M.~E.}}
(\byear{1986}).
\btitle{On an area preserving evolution equation for plane curves}.
In \bbooktitle{Nonlinear Problems in Geometry},
(\beditor{\bfnm{D.~M.}\binits{D.~M.}~\bsnm{DeTurck}}, ed.).
\bseries{Contemporary Mathematics}
\bvolume{51}
\bpages{51--62}.
\end{bincollection}
\endbibitem

\bibitem{Griffiths2009}
\begin{bbook}[author]
\bauthor{\bsnm{Griffiths},~\bfnm{J.~B.}\binits{J.~B.}} \AND
  \bauthor{\bsnm{Podolsk\'y},~\bfnm{J.}\binits{J.}}
(\byear{2012}).
\btitle{Exact space-times in {E}instein's general relativity}.
\bseries{Cambridge Monographs on Mathematical Physics}.
\bpublisher{Cambridge University Press}.
\end{bbook}
\endbibitem

\bibitem{Hawking}
\begin{barticle}[author]
\bauthor{\bsnm{Hawking},~\bfnm{Stephen~W.}\binits{S.~W.}}
(\byear{1968}).
\btitle{{Gravitational radiation in an expanding universe}}.
\bjournal{J. Math. Phys.}
\bvolume{9}
\bpages{598--604}.
\bdoi{10.1063/1.1664615}
\end{barticle}
\endbibitem

\bibitem{Huisken1987}
\begin{barticle}[author]
\bauthor{\bsnm{Huisken},~\bfnm{G.}\binits{G.}}
(\byear{1987}).
\btitle{The volume preserving mean curvature flow}.
\bjournal{J. Reine Angew. Math.}
\bvolume{382}
\bpages{35--48}.
\end{barticle}
\endbibitem

\bibitem{HI}
\begin{binproceedings}[author]
\bauthor{\bsnm{Huisken},~\bfnm{Gerhard}\binits{G.}} \AND
  \bauthor{\bsnm{Ilmanen},~\bfnm{Tom}\binits{T.}}
(\byear{2002}).
\btitle{{Energy Inequalities for Isolated Systems and Hypersurfaces Moving by
  their Curvature}}.
In \bbooktitle{{Proceedings of the 16th International Conference on General
  Relativity and Gravitation}}
(\beditor{\bfnm{Nigel~T.}\binits{N.~T.}~\bsnm{Bishop}} \AND
  \beditor{\bfnm{Sunil~D.}\binits{S.~D.}~\bsnm{Maharaj}}, eds.)
\bpages{162--173}.
\bpublisher{World Scientific}.
\end{binproceedings}
\endbibitem

\bibitem{Jauregui:2011fu}
\begin{barticle}[author]
\bauthor{\bsnm{Jauregui},~\bfnm{Jeffrey~L.}\binits{J.~L.}}
(\byear{2013}).
\btitle{{Fill-ins of nonnegative scalar curvature, static metrics, and
  quasi-local mass}}.
\bjournal{Pacific J. Math.}
\bvolume{261}
\bpages{417--444}.
\end{barticle}
\endbibitem

\bibitem{Komar}
\begin{barticle}[author]
\bauthor{\bsnm{Komar},~\bfnm{Arthur}\binits{A.}}
(\byear{1963}).
\btitle{{Positive-Definite Energy Density and Global Consequences for General
  Relativity}}.
\bjournal{Physical Review}
\bvolume{129}
\bpages{1873--1876}.
\end{barticle}
\endbibitem

\bibitem{Miao2003}
\begin{barticle}[author]
\bauthor{\bsnm{Miao},~\bfnm{P.}\binits{P.}}
(\byear{2003}).
\btitle{On the existence of static metric extensions in general relativity}.
\bjournal{Commun.\ Math.\ Phys.}
\bvolume{241}
\bpages{24--46}.
\end{barticle}
\endbibitem

\bibitem{Papapetrou}
\begin{barticle}[author]
\bauthor{\bsnm{Papapetrou},~\bfnm{A.}\binits{A.}}
(\byear{1948}).
\btitle{A static solution of the equations of the gravitational field for an
  arbitrary charge-distribution}.
\bjournal{Proc. R. Irish Acad. A}
\bvolume{52}.
\end{barticle}
\endbibitem

\bibitem{Pihan1998}
\begin{bphdthesis}[author]
\bauthor{\bsnm{Pihan},~\bfnm{D.~M.}\binits{D.~M.}}
(\byear{1998}).
\btitle{A length preserving geometric heat flow for curves.}
\btype{PhD thesis},
\bpublisher{University of Melbourne}.
\end{bphdthesis}
\endbibitem

\bibitem{Piubello}
\begin{bmastersthesis}[author]
\bauthor{\bsnm{Piubello},~\bfnm{A.}\binits{A.}}
(\byear{2017}).
\btitle{On Existence of Static Metric Extensions \`a la Miao.}
\btype{Master's thesis},
\bpublisher{University of T\"ubingen}.
\end{bmastersthesis}
\endbibitem

\bibitem{Schwarzschild}
\begin{barticle}[author]
\bauthor{\bsnm{Schwarzschild},~\bfnm{K.}\binits{K.}}
(\byear{1916}).
\btitle{{\"Uber das Gravitationsfeld eines Massenpunktes nach der Einsteinschen
  Theorie}}.
\bjournal{Sitzungsberichte der K\"oniglich Preu{\ss}ischen Akademie der
  Wissenschaften}
\bvolume{7}
\bpages{189--196}.
\end{barticle}
\endbibitem

\bibitem{TaylorBook}
\begin{bbook}[author]
\bauthor{\bsnm{Taylor},~\bfnm{M.~E.}\binits{M.~E.}}
(\byear{2011}).
\btitle{Partial Differential Equations III}.
\bseries{Applied Mathematical Sciences}
\bvolume{117}.
\bpublisher{Springer}, \baddress{New York}.
\end{bbook}
\endbibitem

\bibitem{Weyl1917}
\begin{barticle}[author]
\bauthor{\bsnm{Weyl},~\bfnm{H.}\binits{H.}}
(\byear{1917}).
\btitle{Zur {G}ravitationstheorie}.
\bjournal{Ann.\ Physik}
\bvolume{54}
\bpages{117--145}.
\end{barticle}
\endbibitem

\end{thebibliography}
